\documentclass[opre,nonblindrev]{informs3}

\DoubleSpacedXI % Made default 4/4/2014 at request
\usepackage{endnotes}
\let\footnote=\endnote

\usepackage{natbib}
 \bibpunct[, ]{(}{)}{,}{a}{}{,}%
 \def\bibfont{\small}%
 \def\BIBand{and}%
                 
\TheoremsNumberedThrough     % Preferred (Theorem 1, Lemma 1, Theorem 2)
\ECRepeatTheorems

\EquationsNumberedThrough    % Default: (1), (2), ...
\usepackage{bigstrut,bigdelim,multirow}
\usepackage{booktabs}
\usepackage{enumerate}
\usepackage[ruled,vlined,linesnumbered]{algorithm2e}
\usepackage{amsmath}
\usepackage{caption}
\usepackage[normalem]{ulem}

\usepackage{appendix}

\usepackage{subfigure}
\usepackage{float} 

\usepackage[colorlinks,linkcolor=blue,anchorcolor=blue,citecolor=blue]{hyperref}
\begin{document}
%%%%%%%%%%%%%%%%

\graphicspath{{figures/}}

\def\tinyl{\mbox{\tiny L}}
\def\tinyf{\mbox{\tiny F}}
\def\tinyu{\mbox{\tiny U}}
\def\tinyo{\circ}

% Outcomment only when entries are known. Otherwise leave as is and
%   default values will be used.
%\setcounter{page}{1}
%\VOLUME{00}%
%\NO{0}%
%\MONTH{Xxxxx}% (month or a similar seasonal id)
%\YEAR{0000}% e.g., 2005
%\FIRSTPAGE{000}%
%\LASTPAGE{000}%
%\SHORTYEAR{00}% shortened year (two-digit)
%\ISSUE{0000} %
%\LONGFIRSTPAGE{0001} %
%\DOI{10.1287/xxxx.0000.0000}%

% Author's names for the running heads
% Sample depending on the number of authors;
% \RUNAUTHOR{Jones}
% \RUNAUTHOR{Jones and Wilson}
% \RUNAUTHOR{Jones, Miller, and Wilson}
% \RUNAUTHOR{Jones et al.} % for four or more authors
% Enter authors following the given pattern:
\RUNAUTHOR{Qi, Jiang, and Shen}

% Title or shortened title suitable for running heads. 
\RUNTITLE{Sequential Competitive Facility Location}

% Full title. 
\TITLE{Sequential Competitive Facility Location: Exact and Approximate Algorithms}

% Block of authors and their affiliations starts here:
% NOTE: Authors with same affiliation, if the order of authors allows,
%   should be entered in ONE field, separated by a comma.
%   \EMAIL field can be repeated if more than one author
\ARTICLEAUTHORS{%
\AUTHOR{Mingyao Qi}
\AFF{Logistics and Transportation Division, Shenzhen International Graduate School, Tsinghua University, Shenzhen, China,  \EMAIL{qimy@sz.tsinghua.edu.cn}} %, \URL{}}
\AUTHOR{Ruiwei Jiang}
\AFF{Department of Industrial and Operations Engineering, University of Michigan, Ann Arbor, MI, \EMAIL{ruiwei@umich.edu}} %, \URL{}}% Enter all authors
\AUTHOR{Siqian Shen}
\AFF{Department of Industrial and Operations Engineering, University of Michigan, Ann Arbor, MI,  \EMAIL{siqian@umich.edu}}
} % end of the block

\ABSTRACT{%
We study a competitive facility location problem (CFLP), where two firms sequentially open new facilities within their budgets, in order to maximize their market shares of demand that follows a probabilistic choice model. This process is a Stackelberg game and admits a bilevel mixed-integer nonlinear program (MINLP) formulation. We derive an equivalent, single-level MINLP reformulation and exploit the problem structures to derive two valid inequalities, based on submodularity and concave overestimation, respectively. We use the two valid inequalities in a branch-and-cut algorithm to find globally optimal solutions. Then, we propose an approximation algorithm to find good-quality solutions with a constant approximation guarantee. We develop several extensions by considering general facility-opening costs, outside competitors, as well as diverse facility-planning decisions, and discuss solution approaches for each extension. We conduct numerical studies to demonstrate that the exact algorithm significantly accelerates the computation of CFLP on large-sized instances that have not been solved optimally or even heuristically by existing methods, and the approximation algorithm can quickly find high-quality solutions. We derive managerial insights based on sensitivity analysis of different settings that affect customers' probabilistic choices and the ensuing demand. 
}%

%\HISTORY{}

\maketitle

% Text of your paper here

\section{Introduction}
\label{sec:intro}

The competitive facility location problem (CFLP) involves decision games between two or multiple firms, who compete for customer demand of substitutable products or service in a shared market. CFLP arises in a wide variety of applications including opening new retail stores, locating park-and-ride car rental facilities, building charging stations for electric vehicles, etc. In addition, it extends classical location problems, e.g., $p$-median and maximum coverage, to a more complex decision-making environment, in which the market no longer assumes to have a spatial monopoly but to have co-existing competitors and a certain type of consumer patronage behavior.

\paragraph{Competition Type:} There are mainly three types of competition in CFLP: static, sequential, and dynamic \citep{Plastria2001}. In \emph{static} CFLP, a ``newcomer'' firm enters a market and knows \emph{a priori} the existing facilities set up by its competitors, such as their locations and levels of attractiveness. Following a certain customer behavior model which we will detail next, the firm decides where to locate new facilities to maximize its market share. Representative work on static CFLP includes \citet{BenatiHansen2002, Haase2014, Ljubic2018, mai2020multicut} and references therein. On the other hand, the sequential and dynamic CFLPs allow recourse actions of locating facilities by a firm once its competitor opens new facilities. For example, \citet{EiseltLaporte1997,Plastria2008,Kucukaydn2011,Kucukaydn2012,Kress2012,Drezner2015,Gentile2018} consider \emph{sequential} CFLP, in which a leader optimizes facility locations by taking into account a follower's potential location choices, made after the leader takes actions. Therefore, one could view static CFLP as the follower's problem in sequential CFLP, with fixed locations from the leader. In \emph{dynamic} CFLP, competing firms in a market make non-cooperative decisions simultaneously and iteratively, until a Nash equilibrium, if any, is reached \citep{Godinho2010}. In this paper, we focus on sequential CFLP, where (i) a firm can react to the competitor's actions by locating in remaining candidate sites (as opposed to static CFLP), and (ii) location decisions are relatively expensive and long-lasting, making dynamic relocation economically prohibitive (as opposed to dynamic CFLP). Notably, we focus on decisions of locating facilities made at the strategic-planning level, rather than short-term operational-level actions and therefore, we only consider long-term impacts of facility locations on aggregated customer demand trends and do not take into account daily sales strategies, dynamic pricing of certain products, and other operational decisions that can affect daily demand realizations. A typical application of our CFLP model is to locate chained business facilities such as hotels, supermarkets, and shopping malls, for which an investor needs to make a long-term facility deployment plan with foresight of its competitors' response. In Appendix \ref{sec:extension:attractiveness}, we extend our CFLP model to co-optimize locations and attractiveness levels of the facilities. 

\paragraph{Customer Behavior:} In a CFLP, customers are assumed to act as independent decision makers based on utilities they can receive from facilities. The utility typically depends on the distance to each facility, facility size, service price, etc. \citep[see, e.g.,][for a related empirical study]{okelly1999trade}. Then, a choice behavior model can be used to translate the utilities into how customers patronize facilities. The utility and choice models make key ingredients for CFLP and can largely determine the computational complexity as well as solution approaches. For example, a \textit{deterministic} choice model assumes that each customer purchases all the goods from a facility that has the highest utility, e.g., the closest facility. Accordingly, sequential CFLP with deterministic customer choice admits a mixed-integer \emph{linear} program (MILP) formulation \citep[see][]{Plastria2008, Roboredo2013, Alekseeva2015, Drezner2015, Gentile2018}, which can be efficiently solved by off-the-shelf solvers. In contrast, probabilistic choice models split customer demand across multiple facilities with certain probabilities. For example, in the well-celebrated multinomial logit (MNL) model~\citep{mcfadden1973conditional}, the probability of patronizing a facility is proportional to the natural exponential of the utility. The MNL model has been widely used in static CFLP \citep[see, e.g.,][]{BenatiHansen2002, Haase2014, Ljubic2018, mai2020multicut}, but receives much less consideration in sequential CFLP~\citep[see, e.g.,][]{Kucukaydn2012}, partly because it gives rise to a mixed-integer \emph{nonlinear} program (MINLP) and renders a significant computational challenge when presented in a bilevel formulation. 

In this paper, we investigate sequential CFLP with probabilistic customer choice following the MNL setting. Our goal is to derive exact and approximate formulations, and solution approaches that achieve high computational efficiency for solving instances with real-world problem sizes. In the ensuing bilevel model, both the upper- and lower-level problems (for the leader and the follower, respectively) are MINLPs given MNL-based customer demand. We recast the bilevel program as an equivalent single-level MINLP and derive two classes of valid inequalities to solve the reformulation exactly. We also propose an alternative algorithm to solve the model approximately and show that this algorithm admits a constant approximation guarantee. In Appendix~\ref{sec:extension:outside}, we extend the CFLP model to take into account outside competitors, in addition to the leader and the follower.

\paragraph{Connections with Network Interdiction:} Bilevel programs and their solution approaches have also been applied to network interdiction problems. Below we review the most relevant papers on network interdiction and their connections with this paper. The leader-follower Stackelberg game \citep{washburn1995two} is widely used in the modeling of network interdiction problems, where the leader seeks a set of moves that would minimize the maximum gain, or maximize the minimum loss of the follower over a network \citep[see, e.g.,][]{morton2007models, wood2010bilevel, shen-thesis, dimitrov2013interdiction, song-shen-2015}. With the leader's decisions influencing the follower's optimization problem, the sequential network interdiction games are often formulated and solved through bilevel optimization techniques. We refer to \citet{lim-smith-2007, smith2020survey} for thorough surveys of various interdiction problems, their formulations, and solution methods. \citet{colson2005bilevel,Bard98} reviewed the literature of bilevel programming models, algorithms, and applications, while \citet{kleinert2021survey, ralphs2015bilevel} conducted comprehensive surveys of (mixed) integer programming approaches for solving bilevel integer  programs. Recently, \citet{tahernejad2020branch} described a branch-and-cut algorithm for mixed-integer bilevel linear programming and \citet{bilevel-solver-ralphs} published an open-source solver that can handle bilevel programs with MILPs at both levels.  The sequential CFLP resembles the network interdiction problem in that (i) it can be viewed as a zero-sum game between a leader and a follower and (ii) the leader's decisions may ``interdict'' certain facilities/assets to be used by the follower (see constraints~\eqref{bl-f} in the (S-CFLP) model and, e.g., Section 6 of~\cite{kleinert2021survey}). Nonetheless, while most network interdiction problems consider linear objective functions, the sequential CFLP we study in this paper employs a \emph{nonlinear} objective function and admits a bilevel program with MINLPs at both levels. Consequently, although it is possible to compute an $\epsilon$-optimal solutions for the leader, the existing solution approaches have been applied to solving problems with rather small sizes due to the general hardness of bilevel MINLPs (see, e.g., Section 5.3 of~\cite{kleinert2021survey}). 

In Table \ref{table 1}, we summarize representative work of different sequential CFLP models based on deterministic/probabilistic choice models, their formulations and solution methods. Table \ref{table 1} also indicates whether each paper assumes continuous or discrete location space. If a bilevel formulation is adopted, we specify the formulation types in the upper/lower levels. We also indicate exact solution methods that guarantee global optimum. We refer the interested readers to Appendix \ref{sec:lit} for a detailed review of the existing works on static and sequential CFLP.
% Table generated by Excel2LaTeX from sheet 'literature'
\begin{table}[ht!]
  \centering
  \caption{Comparison with existing work on sequential CFLP}
   \label{table 1}
   \resizebox{\textwidth}{!}{
     \begin{tabular}{lllll}
    %\begin{tabular}{{p{11em}p{6em}p{10em}p{7em}}}
    \toprule
    Reference & Choice Model & Formulation (Upper/Lower) & Solution Approach & Location Space\\
    \midrule
    \citet{Drezner1998} & {Probabilistic} & {Bilevel (NLP/NLP)} & {Heuristic} & {Planar}\\
   \citet{Serra1994}  & Deterministic & MILP & Heuristic  & {Discrete}\\
   \citet{Fischer2002} & Deterministic   & Bilevel (MINLP/MILP) & Heuristic  & {Discrete}\\
   \citet{Plastria2008} & Deterministic  & MILP  & Exact; commercial solver & {Discrete}\\
   \citet{saiz2009branch} & {Probabilistic} & {Bilevel (NLP/NLP)} & {Exact; branch-and-bound} &  {Planar}\\ 
   \citet{Kucukaydn2011}  & Probabilistic & Bilevel (MINLP/NLP) & Exact; single-level MINLP reformulation & {Discrete}\\
  \citet{Kucukaydn2012}  & Probabilistic & Bilevel (MINLP/MINLP) & Heuristic & {Discrete} \\
   \citet{Roboredo2013} & Deterministic  & MILP  & Exact; branch-and-cut & {Discrete} \\
   \citet{Alekseeva2015} & Deterministic  & MILP  & Exact; iterative method & {Discrete} \\
   \citet{Drezner2015}  & Deterministic  & Bilevel (MINLP/MILP) & Heuristic & {Discrete} \\
   \citet{Gentile2018}  & Deterministic  & MILP  & Exact; branch-and-cut & {Discrete} \\
    This paper & Probabilistic & Bilevel (MINLP/MINLP) & Exact and approximate; branch-and-cut & {Discrete}\\
    \bottomrule
    \end{tabular}
   }
\end{table}

To the best of our knowledge, this paper provides the first exact solution approach, as well as the first approximation algorithm with an optimality guarantee, for solving sequential CFLP with probabilistic customer choice.
To achieve this, we interpret sequential CFLP as a robust optimization (RO) model, in which the follower acts adversarially in order to decrease the market share of the leader's. Nonetheless, as we shall see in Section \ref{sec:model}, the uncertainty set of this RO model depends on the leader's locations, leading to an intractable, decision-dependent RO model.

\paragraph{Main Contributions:}
We summarize the main contributions of the work in the following aspects. 
\begin{enumerate}[1.]
\item Without loss of optimality, we revise the objective function of sequential CFLP to make the uncertainty set of the equivalent RO-based reformulation decision-\emph{in}dependent. Accordingly, we recast the bilevel program as a single-level MINLP. This allows us to solve sequential CFLP to global optimum using a finitely-convergent branch-and-cut algorithm.
\item We derive two classes of valid inequalities to accelerate the branch-and-cut procedure. We further derive an approximate separation of these valid inequalities that only consumes a sorting procedure to compute, to significantly speed up the computation. 
\item We propose an approximation algorithm for solving sequential CFLP based on a mixed-integer second-order conic program (MISOCP), which can be readily solved in off-the-shelf optimization solvers. In addition, we derive a constant approximation guarantee on the ensuing market share.
\item Through extensive computational experiments, we demonstrate that our valid inequalities can significantly accelerate solving sequential CFLP. For example, our approach is able to solve instances with 100 candidate facilities and 2000 customer nodes in minutes, which has never been achieved in the CFLP literature  (even by heuristic approaches). In addition, our approximation algorithm can obtain good-quality solutions even more quickly. We also report results of varying parameters in the probabilistic choice model and how they impact the optimal location decisions of the leader and the follower. 
\item We extend the sequential CFLP model to incorporate additional features and constraints, including general facility-opening costs, outside competitors, attractiveness level, and utility change. We also show that our reformulation and branch-and-cut algorithm can still be applied to obtain exact or approximate solutions to these extensions. 
\end{enumerate}

\paragraph{Structure of the Paper:}
The remainder of the paper is organized as follows. In Section \ref{sec:model}, we formulate the bilevel model and its MINLP reformulation. In Section \ref{sec:cuts}, we derive valid inequalities, an approximate separation procedure used in a branch-and-cut algorithm, and also an alternative approximation algorithm. 
We present computational results based on instances with diverse sizes and complexity in Section \ref{sec:compu}. In Section \ref{sec:concl}, we conclude the work and propose future research directions. 

\noindent {\bf Notation}: For $a, b \in \mathbb R$, we define $a\vee b := \max\{a, b\}$. For set $S$, $|S|$ denotes its cardinality and $2^S$ denotes the collection of all its subsets. The notation $\|\cdot\|_2$ represents the $2$-norm in the real space, and $e$ represents an all-one vector with suitable dimension.

\section{Bilevel Model and Single-level Reformulation}
\label{sec:model}
In sequential CFLP, two firms (a leader and a follower) deploy facilities to provide substitutable commodities to customers located in a set $I$ of nodes. In each node $i\in I$, there is $h_i$ portion of the total customer demand to patronize these facilities, i.e., $\displaystyle\sum_{i \in I}h_i = 1$. The leader and the follower may already have existing facilities in the market, denoted by sets $J^{\tinyl}$ and $J^{\tinyf}$ respectively, with $J^{\tinyl} \cap J^{\tinyf} = \varnothing$. The new facilities may be deployed in a set $J$ of candidate sites such that $J \cap (J^{\tinyl} \cup J^{\tinyf}) = \varnothing$. The competition follows a Stackelberg game \citep{von1934marktform}, in which the leader first locates at most $p$ facilities to maximize the leader's market share while foreseeing that the follower will react and locate at most $r$ facilities in the remaining candidate sites to maximize the follower's market share, where $p+r \leq |J|$. Without loss of generality, we assume that two or more facilities do not co-locate at one candidate site, because mathematically we can always split a candidate site and the corresponding utilities if needed. 

We adopt the widely-used MNL model~\citep[see, e.g.,][]{mcfadden1973conditional,ben1985discrete} for probabilistic customer choice. Specifically, MNL assumes that a customer dwelling at node $i$ and patronizing a facility deployed at site $j$ receives utility $u_{ij} := \alpha_j - \beta d_{ij} + \epsilon_{ij}$, where $\epsilon_{ij}$ denotes a random noise, $-\beta < 0$ denotes the negative impact of traveling distance $d_{ij}$ between node $i$ and site $j$, and $\alpha_j$ denotes the attractiveness of facility $j$, which depends on various characteristics such as size, reputation, price levels, etc. In a seminal work \citep{mcfadden1973conditional}, the author shows that if the random noises $\epsilon_{ij}$ are independent and identically follow the standard Gumbel distribution then the probability $P_{ij}$ of the customer patronizing facility $j$ follows 
$P_{ij} = \frac{\exp\{\alpha_j - \beta d_{ij}\}}{\sum_{k \in J^0}\exp\{\alpha_k - \beta d_{ik}\}}$, 
where $J^0$ denotes the set of facilities deployed by either the leader or the follower. For notational brevity, we denote $w_{ij}:=\exp\{\alpha_j - \beta d_{ij}\}$ for all $i \in I$ and $j \in J\cup J^{\tinyl}\cup J^{\tinyf}$, and let $U^{\tinyl}_i := \sum_{j \in {J^{\tinyl}}} w_{ij}$ and $U^{\tinyf}_i := \sum_{j \in {J^{\tinyf}}} w_{ij}$ denote utility of the pre-existing facilities already open by the leader and the follower, respectively. For new facilities, define binary variables $x_j$ and $y_j$, for all $j \in J$, to indicate whether or not the leader/follower deploys a facility at site $j$, respectively. The leader's market share is given by \small
\begin{equation}
\label{eq:l+}
L^{+}(x, y) := \sum_{i\in I} h_i \left(\frac{U^{\tinyl}_i+\sum_{j\in J}w_{ij} x_j}{U^{\tinyl}_i+U^{\tinyf}_i+\sum_{j\in J}w_{ij}(x_{j}+y_{j})}\right).
\end{equation}
\normalsize
Accordingly, we formulate sequential CFLP with probabilistic customer choice as a bilevel program: 
\begin{subequations}\small
\label{bl}
\begin{align}
(\mbox{\bf S-CFLP}) \ \ \ \max_{x} \ & \ L^{+}(x, y^*)    \label{bl-a} \\
      \rm{s.t.} \ & \ \sum_{j\in J}x_{j}\leq p, \label{bl-b} \\
         & \ x_{j}\in\{0,1\}, \quad \forall j\in J, \label{bl-c} 
\end{align}
\begin{align}
\mbox{where} \qquad y^* \in \argmax \ & \ \sum_{i\in I}h_{i} \left(\frac{U^{\tinyf}_i+\sum_{j\in J}w_{ij} y_j}{U^{\tinyl}_i+U^{\tinyf}_i+\sum_{j\in J}w_{ij}(x_{j}+y_{j})}\right) \label{bl-d}\\
      \rm{s.t.} \ & \ \sum_{j\in J}y_{j}\leq r, \label{bl-e} \\
         & \ y_{j}\leq 1-x_{j}, \quad \forall j\in J,  \label{bl-f} \\
         & \ y_{j}\in\{0,1\}, \quad \forall j\in J. \label{bl-g}
\end{align}
\end{subequations}
The objective function \eqref{bl-a} of the upper-level problem aims to maximize the leader's market share, and constraints \eqref{bl-b}--\eqref{bl-c} ensure that the leader deploys no more than $p$ new facilities. In the lower-level problem, the objective function \eqref{bl-d} derives an optimal solution $y^*$ that maximizes the follower's market share given $x$. Constraints \eqref{bl-e} and \eqref{bl-g} ensure that the follower opens up to $r$ new facilities, and constraints \eqref{bl-f} prohibit co-location. Both upper- and lower-level problems are MINLPs, giving rise to significant computational challenges. In fact, the bilevel MINLP formulation of (S-CFLP) suggests that it is $\Sigma^p_2$-hard (see, e.g.,~\cite{jeroslow1985polynomial}). The following hardness result indicates that it is not even possible to find a good approximation of (S-CFLP) unless $\mbox{P} = \mbox{NP}$, for which we present a detailed proof in Appendix \ref{apx-thm-hardness}.
\begin{theorem}[Adapted from Theorem 3 of~\cite{krause2008robust}] \label{thm:hardness}
There does not exist a polynomial-time, constant approximation algorithm for (S-CFLP) unless $\mbox{P} = \mbox{NP}$. Specifically, let $z^*$ represent the optimal value of (S-CFLP). If there exists a constant $c_0 > 0$ and an algorithm, which runs in time polynomial in $|J|$, $p$ and guarantees to find a solution $x$ such that $L^+(x, y^*) \geq c_0 z^*$, then $\mbox{P} = \mbox{NP}$.
\end{theorem}

Next, we take a RO perspective to recasting (S-CFLP) into a computable form. We start by noticing that the leader's and the follower's objective functions \eqref{bl-a} and \eqref{bl-d} sum up to $1$. Intuitively, this is because a customer patronizes either the leader or the follower. (In this section, we focus on the case of having two firms only. In reality, however, customers may patronize outside competitors. We will extend the baseline (S-CFLP) model to take this into account in Appendix~\ref{sec:extension:outside}.) This suggests that formulation \eqref{bl} can be viewed as a RO model, in which the follower acts adversarially to decrease the market share of the leader's. From this perspective, (S-CFLP) is equivalent to a RO model:%\small
\begin{equation}
\label{eq:max-min}
\max_{x \in \mathcal X} \ \min_{y \in \mathcal{Y}(x)} L^{+}(x, y),
\end{equation}
where $\mathcal X := \{ x \in \{0, 1\}^{|J|} : e^{\top}x \leq p \}$ denotes the leader's budget for opening facilities and $\mathcal{Y}(x) := \left\{y \in \{0, 1\}^{|J|} : \mbox{\eqref{bl-e}--\eqref{bl-g}}\right\}$ denotes the follower's, which is interpreted as an uncertainty set in RO. Unfortunately, this uncertainty set is decision-\emph{dependent}, preventing us from applying standard reformulation techniques. For example, one might suggest to solve formulation \eqref{eq:max-min} by rewriting the inner formulation in a hypographic form and applying delayed constraint generation (DCG). Specifically, one can rewrite \eqref{eq:max-min} as $\displaystyle\max_{x \in \mathcal{X}, \theta^+}\left\{\theta^+: \ \theta^+ \leq L^+(x, y), \ \forall y \in \mathcal{Y}(x)\right\}$ and iteratively incorporates inequalities (cuts) $\theta^+ \leq L^+(x, y)$ only when they are violated. Unfortunately, this is not applicable, because the validity of these constraints depends on the values of $x$. For example, suppose that we incorporate a cut $\theta^+ \leq L^+(x, \widehat{y})$, where $\widehat{y} \in  \mathcal{Y}(\widehat{x})$ for an $\widehat{x} \in \mathcal{X}$. Then, this cut fails to be valid whenever $\widehat{y} \notin  \mathcal{Y}(\overline{x})$ for a different $\overline{x} \in \mathcal{X}$ because it would (incorrectly) undervalue the objective function at $\overline{x}$, or equivalently, $\displaystyle\min_{y \in  \mathcal{Y}(\overline{x})}L^+(\overline{x}, y)$. Another possibility is to take the dual of the inner minimization formulation and produce a single-level formulation. Unfortunately, this does not apply to \eqref{eq:max-min} either, because $\mathcal{Y}(x)$ involves binary restrictions \eqref{bl-g} and hence strong duality fails to hold for the inner formulation. For the same reason, the exactness of this formulation would be lost if we replace the inner formulation with its Karush-Kuhn-Tucker conditions. Also, note that constraints \eqref{bl-g} may not be na{\"{i}}vely relaxed because $L^+(x, y)$ is convex in $y$.

To make the uncertainty set decision-\emph{in}dependent and relax the co-location constraints~\eqref{bl-f} without loss of optimality, the network interdiction literature suggests adding a penalty term $\sum_{j \in J}M_j y_j x_j$ to the objective function $L^+(x, y)$, where $M_j$ represent sufficiently large positive numbers (see, e.g., Section 3.1.1 of~\cite{smith2020survey} and Section 6.2 of~\cite{kleinert2021survey}). The choice of the coefficients $M_j$ is crucial because large $M_j$ significantly weaken $L^+(x, y)$ and the continuous relaxation of formulation~\eqref{eq:max-min}. We adopt an alternative approach without weakening $L^+(x, y)$. We specify the result in the following theorem and present a proof in Appendix~\ref{apx-thm-submodularity}.
\begin{theorem}
\label{lemma:l1}
Define $\theta^+, \theta: \{0,1\}^{|J|} \rightarrow \mathbb{R}$ such that $\displaystyle \theta^+(x) := \min_{y \in \mathcal{Y}(x)} L^+(x, y)$ and $\displaystyle\theta(x) := \min_{y \in \mathcal{Y}} L(x, y)$, where $\mathcal{Y} := \{y: \mbox{\eqref{bl-e}, \eqref{bl-g}} \}$ and 
\begin{equation}
\label{eq:l}
L(x, y):=\sum_{i\in I}h_{i}\left(\frac{U^{\tinyl}_i+\sum_{j\in J}w_{ij} x_j}{U^{\tinyl}_i+U^{\tinyf}_i+\sum_{j\in J}w_{ij}(x_{j}\vee y_{j})}\right).
\end{equation}
Then, it holds that $\theta^+(x)=\theta(x)$ for all $x\in\mathcal{X}$.
\end{theorem}

Theorem \ref{lemma:l1} recasts the bilevel model ({S-CFLP}) as the following single-level MINLP: %\small
\begin{subequations}
\label{sl}
\begin{align}
   \max_{x \in \mathcal{X}, \theta} \ & \ \theta \\
   \mbox{s.t.} \ & \ \theta\leq L(x,y), \quad \forall y \in \mathcal{Y}. \label{sl-b}
\end{align}
\end{subequations} 
Although model \eqref{sl} incorporates an exponential number of constraints due to the cardinality of set $\mathcal{Y}$, it can be solved through DCG. Specifically, we relax constraints \eqref{sl-b} and iteratively add them back if needed. In each iteration, we obtain an incumbent solution $(\hat{x}, \hat{\theta})$ from the relaxed formulation. Then, we solve the following separation problem
\begin{equation}
\label{sp}
    \min_{y} \Bigl\{L(\hat{x},y) \ : \ y \in \mathcal{Y} \Bigr\}
\end{equation}
to decide if this solution violates any of constraints \eqref{sl-b}. If not, then $(\hat{x}, \hat{\theta})$ is an optimal solution to \eqref{sl}; otherwise, we find a $\hat{y} \in \mathcal{Y}$ such that constraint $\theta \leq L(x, \hat{y})$ is violated, i.e., $\hat{\theta} > L(\hat{x}, \hat{y})$. We append this violated constraint in the iteratively-solved relaxed formulation of \eqref{sl} to cut off the incumbent solution. Since $\mathcal{Y}$ is finite, DCG terminates with a global optimal solution in a finite number of iterations. We summarize full details of the DCG approach in Algorithm \ref{algo-b&c} at the end of Section \ref{sec:cuts}, after completing the derivation of  the valid inequalities and the approximate separation.

\section{Valid Inequalities, Approximate Separation, and Approximation Algorithm}
\label{sec:cuts}
There are two challenges on applying DCG. First, the violated constraints we incorporate, $\theta \leq L(x, \hat{y})$, are nonlinear. As a consequence, in every iteration we need to solve the relaxed formulation of model \eqref{sl} as a MINLP. Moreover, as we shall see in Section \ref{sec:cuts-oa}, function $L(x, \hat{y})$ is non-concave in $x$ in its current form presented in \eqref{eq:l}. That is, the relaxed formulation remains a non-convex NLP even if we further relax its integer restrictions. To address these challenges, we derive two classes of \emph{linear} valid inequalities in Sections \ref{sec:cuts-sc} and \ref{sec:cuts-oa}, respectively. Together, they generate a tight MILP relaxation of the nonlinear, non-convex formulation, which can be readily solved by off-the-shelf solvers. Second, the separation problem \eqref{sp} by itself is a MINLP. To address this, we recall that \eqref{sp} is equivalent to a static CFLP model and solve it via the state-of-the-art approach from~\citet{Ljubic2018}. As a further improvement, in Section \ref{sec:cuts-approx}, we derive an approximate, but much faster, approach to solving \eqref{sp} and generating valid inequalities via a single-round sorting. The two valid inequalities and approximate separation procedures are used in a branch-and-cut algorithm described in Section \ref{sec:cuts:bcf}, and we derive an approximation algorithm for (S-CFLP) with a constant approximation guarantee in Section \ref{sec:heuristic}. 

\subsection{Submodular Inequalities} 
\label{sec:cuts-sc}
We start by recalling the following definition of submodular functions.
\begin{definition}[Submodular Functions]
A function $f:2^J \rightarrow \mathbb{R}$ is submodular if %\small
$$
f(S\cup\{j\}) - f(S) \geq f(R\cup\{j\}) - f(R)
$$
for all subsets $S \subseteq R \subseteq J$ and all element $j \in J \setminus R$.
\Halmos
\end{definition}
Intuitively, $f$ is submodular if the marginal gain of incorporating any additional element $j$ is non-increasing in the subset $S$. We show that, for any fixed $y \in \mathcal{Y}$, the function $L(x, y)$ is submodular with respect to the index set $X$ of $x$, i.e., $X := \{j \in J: x_j = 1\}$. This observation enables us to represent the nonlinear constraint $\theta \leq L(x, y)$ as a set of linear inequalities.

To state the results formally, we denote $Y := \{j \in J: y_j = 1\}$ as the index set of $y$. In addition, we define set functions $L_Y$: $2^J\rightarrow \mathbb{R_+}$ and $f_{i,Y}$ : $2^J\rightarrow \mathbb{R_+}$ such that $\displaystyle L_Y(X) := \sum_{i\in I}h_{i}f_{i,Y}(X)$ and 
\begin{equation}
\label{eq:f-i-J'}
  f_{i,Y}(X):=\frac{U^{\tinyl}_i+\sum_{j\in X}w_{ij}}{U^{\tinyl}_i+U^{\tinyf}_i+\sum_{j\in X\cup Y}w_{ij}}
\end{equation}
for all $i \in I$. Intuitively, $f_{i,Y}(X)$ evaluates the percentage of demand from customer $i$ that patronizes the leader's facilities and $L_Y$ evaluates the total market share of the leader, if the leader and the follower deploy facilities in sets $X$ and $Y$, respectively. Hence, it holds that $L(x,y) = L_{Y}(X)$.
\begin{proposition} \label{p2}
For any $Y \subseteq J$, $L_{Y}$ is submodular.
\end{proposition}
A detailed proof of Proposition \ref{p2} is in Appendix \ref{apx-p2}. Since $L_Y$ is submodular, we follow~\citet{nemhauser1981maximizing} to rewrite the constraint $\theta \leq L(x, y) \equiv L_Y(X)$ as a set of linear inequalities.
\begin{proposition}[Adapted from Theorem 6 of~\citet{nemhauser1981maximizing}] \label{p3}
For any $y \in \mathcal{Y}$, the constraint $\theta \leq L(x,y)$ is equivalent to the following linear constraints:
\begin{align}
\theta \le L_Y (S) - \sum_{k \in S} \rho_Y (J\setminus \{k\}; k)(1 - x_k) + \sum_{k \in J\setminus S} \rho_Y (S;k)x_k, \quad \forall{S \subseteq J}, \label{7}
\end{align}
where $\rho_Y (S; k) := L_Y (S \cup k) - L_Y (S) $ for all $S \subseteq J$ and $k \in J \setminus S$.
\end{proposition}
Constraints \eqref{7} involve an exponential number of inequalities. Hence, in DCG, a straightforward replacement of $\theta \leq L(x,y)$ with \eqref{7} drastically increases the formulation size. Instead, we can replace $\theta \leq L(x,y)$ with the \emph{most violated} inequality among \eqref{7}, which serves the same purpose of cutting off the incumbent solution. Specifically, for given $(\hat{x}, \hat{\theta})$, inequalities \eqref{7} hold valid \emph{iff} %if and only if 
\begin{align}
\hat{\theta} \leq \min_{S \subseteq J} \left\{L_Y (S) - \sum_{k \in S} \rho_Y (J\setminus \{k\}; k)(1-\hat{x}_k) + \sum_{k \in J\setminus S} \rho_Y (S;k) \hat{x}_k\right\},  \label{8}
\end{align}
and to find the most violated inequality it suffices to solve the combinatorial optimization problem on the right-hand side of \eqref{8}. In what follows, we show that this task can be accomplished efficiently, in polynomial time. 

Suppose that $\hat{x} \in \{0, 1\}^{|J|}$. This can take place when a relaxation of model \eqref{sl} happens to produce a binary-valued solution or at a leaf node of the branch-and-bound tree for solving \eqref{sl}. For this case,~\citet{Ljubic2018} showed that the index set $\widehat{X}$ of $\hat{x}$, i.e., $\widehat{X} := \{j \in J: \hat{x}_j = 1\}$ is optimal to problem \eqref{8}. That is, the most violated inequality among \eqref{7} is the one with $S = \widehat{X}$.

More generally, we consider $\hat{x} \in [0, 1]^{|J|}$, i.e., $\hat{x}$ can be either fractional- or binary-valued. This can take place at any node of the branch-and-bound tree, where we relax (a part of or all of) the binary restrictions on variables $x$. The next proposition shows that problem \eqref{8} has a submodular objective function and so it admits a polynomial-time solution~\citep[see][]{edmonds1970submodular,topkis1978minimizing}.
\begin{proposition} \label{p4}
For any $x \in [0,1]^{|J|}$ and $Y \subseteq J$, define $H:2^J \to \mathbb{R}$ such that $\displaystyle H(S) := L_Y (S) - \sum_{k \in S} \rho_Y (J\setminus{k}; k)(1-x_k) + \sum_{k \in J \setminus S} \rho_Y (S;k) x_k$. Then, $H$ is submodular.
\end{proposition}
We present a proof of Proposition \ref{p4} in Appendix \ref{apx-p4}. Proposition \ref{p4} indicates that we can find the most violated inequality among \eqref{7} in polynomial time. 

\subsection{Bulge Inequalities} \label{sec:cuts-oa}
A second alternative of the nonlinear constraints $\theta \leq L(x, y)$ is the supporting hyperplanes of the hypograph of function $L(x, y)$, also known as the outer approximation method~\citep[see, e.g.,][]{duran1986outer,Ljubic2018}. To this end, we extend the domain of $L(x, y)$ by defining
$$
\tilde{L}(x, y) := \sum_{i\in I}h_{i}\left(\frac{U^{\tinyl}_i+\sum_{j\in J}w_{ij}x_j}{U^{\tinyl}_i+U^{\tinyf}_i+\sum_{j\in J}w_{ij}\bigl[(1-y_j)x_j+y_j\bigr]}\right),
$$
where we replace the $x_j \vee y_j$ in $L(x, y)$ with $(1-y_j)x_j+y_j$. Note that $L(x, y)$ coincides with $\tilde{L}(x, y)$ whenever $(x, y)$ are binary-valued, but $\tilde{L}(x, y)$ is well-defined on $[0, 1]^{2|J|}$. For fixed $y \in \mathcal{Y}$, we can replace $\theta \leq L(x, y)$ with a supporting hyperplane if $\tilde{L}(x, y)$ is concave in $x$. Unfortunately, this fails to hold as evidenced by the following example.
\begin{example}[Non-concavity of $\tilde{L}$] \label{exam:non-concave}
Suppose that $I = \{1\}$ and $h_1 = 1$, i.e., there is one single customer node. In this example, we ignore the index $i$ for notational simplicity. In addition, suppose that $J = \{1, 2\}$, $U^{\tinyl}_i = U^{\tinyf}_i = 0$, $w_1 = w_2 = 1$, and $y_1 = 1 - y_2 = 0$. Then, $\tilde{L}(x, y) = \frac{x_1 + x_2}{x_1 + 1}$ and its Hessian reads
$$
-\frac{1}{(x_1 + 1)^2}
\begin{bmatrix}
\frac{2(1 - x_2)}{x_1 + 1} \quad & \quad 1\\
1 \quad & \quad 0
\end{bmatrix},
$$
which is not negative semidefinite on $(0, 1)^{|J|}$. Therefore, $\tilde{L}(x, y)$ is not concave in $x$. In particular, restricting $\tilde{L}(x, y)$ on the line $x_1 + x_2 = 1$ yields $\frac{1}{x_1 + 1}$, which is in fact a convex function.
\Halmos
\end{example}
The above example suggests that we should ``bulge up'' $\tilde{L}(x, y)$ in order to obtain a concave hypograph. To this end, we replace the linear term $x_j$ in the numerator of the $\tilde{L}(x, y)$ definition with a larger, quadratic term. This yields a desired concave function as shown next, for which we provide a detailed proof in Appendix \ref{apx-lemma-l2}.
\begin{proposition}
\label{lemma:l2}
For fixed $y \in \{0, 1\}^{|J|}$, define $\widehat{L}: [0, 1]^{|J|} \rightarrow \mathbb{R}_+$ such that
\begin{equation}
\label{eq:L-bar}
\widehat L (x, y):=\sum_{i\in I}h_{i}\left(\frac{U^{\tinyl}_i+\sum_{j\in J}w_{ij}{\bigl[- y_j x_j^2+(1+ y_j)x_j\bigr]}}{U^{\tinyl}_i+U^{\tinyf}_i+\sum_{j\in J}w_{ij}\bigl[(1-y_j)x_j+y_j\bigr]}\right).
\end{equation}
Then, $\widehat L (x,y)$ is concave in $x$. In addition, $\widehat L (x,y) = L(x, y)$ for all $x \in \{0, 1\}^{|J|}$.
\end{proposition}

\begin{figure}[ht!]
\centering  
\subfigure[$\tilde{L}(x,y)$]{
\label{figure 2a}
\includegraphics[width=0.45\textwidth]{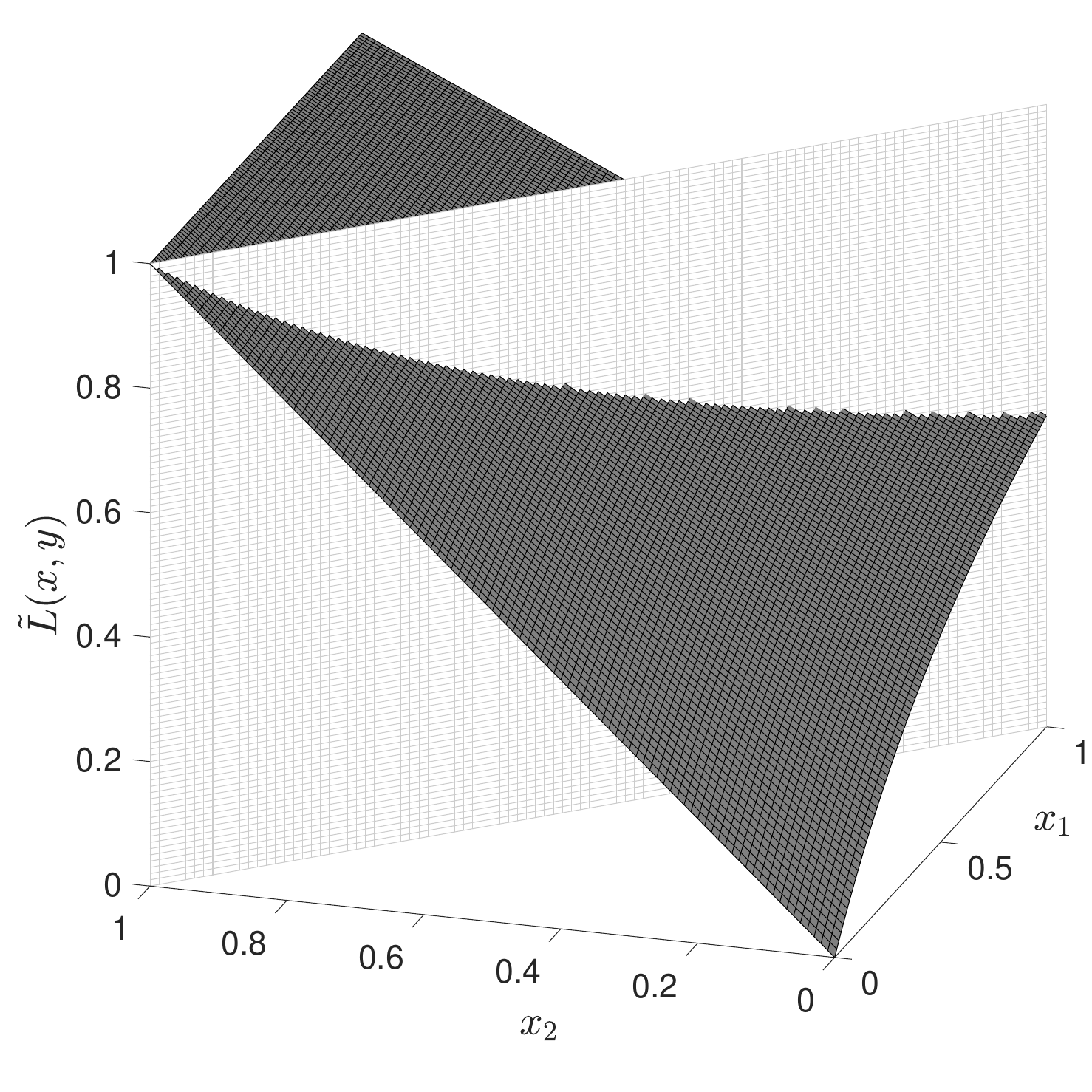}}
\subfigure[$\widehat L(x,y)$]{
\label{figure 2b}
\includegraphics[width=0.45\textwidth]{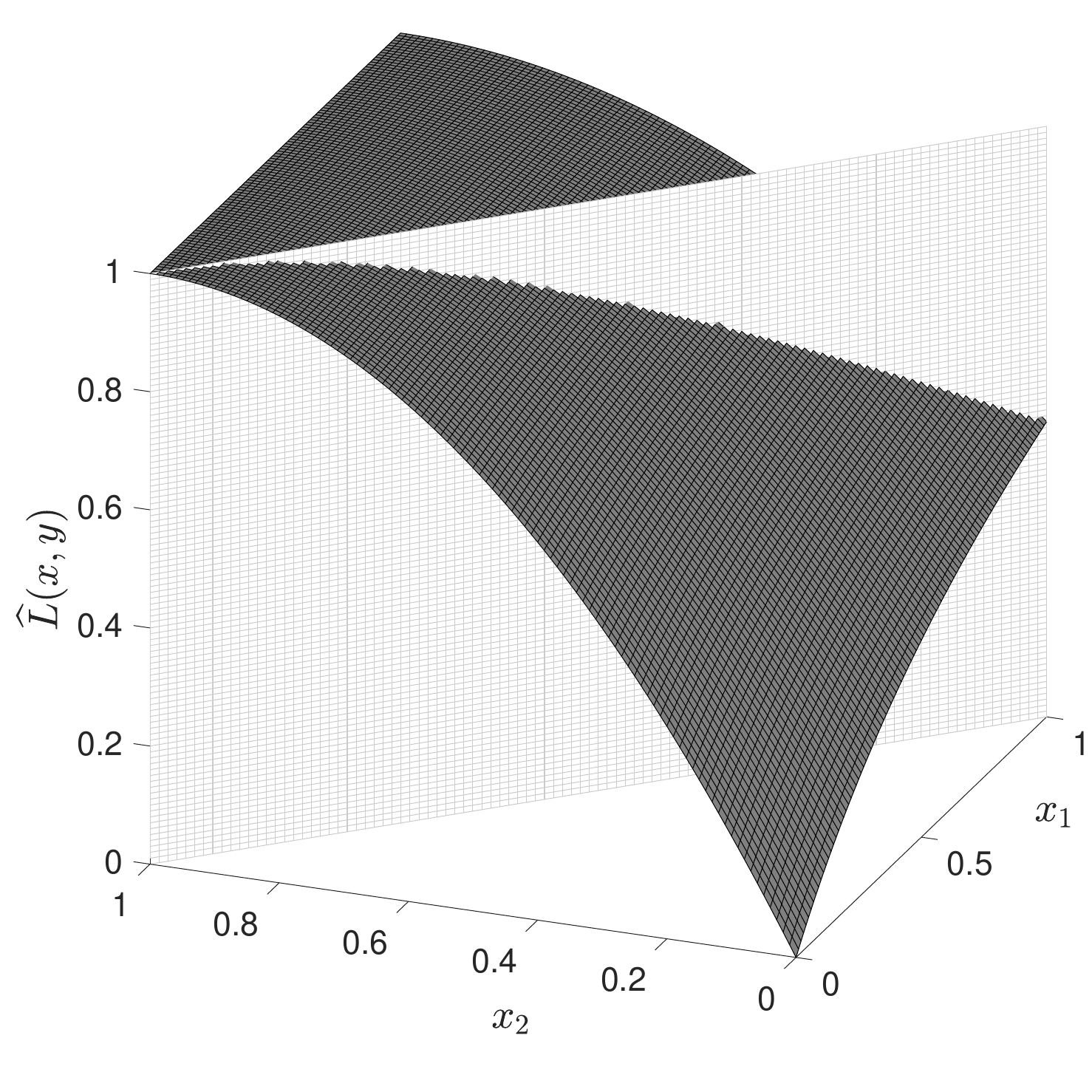}}
\caption{Illustration of functions $\tilde{L}$ and $\widehat L$ (the black surfaces) and their restrictions on the line $x_1 + x_2 = 1$ (the intersecting curves of the gray and black surfaces) with $y = [0,1]^{\top}$} \label{figure-concave}
\end{figure}

Proposition \ref{lemma:l2} indicates that $\widehat{L}(x, y)$ is a concave representation of function $L(x, y)$. That is, $\widehat{L}$ bulge up $L$ to make it concave while retaining the exactness at any binary-valued $x$. We illustrate this observation below.
\begin{example}[Concavity of $\widehat{L}$]
Continuing from Example \ref{exam:non-concave}, we construct $\widehat{L}(x, y) = \frac{x_1 - x_2^2 + 2x_2}{x_1 + 1}$ for the same $y = [0, 1]^{\top}$. Then, its Hessian reads
$$
-\frac{2}{(x_1 + 1)^3}
\begin{bmatrix}
(1-x_2)^2 \quad & \quad (x_1+1)(1-x_2)\\
(x_1+1)(1-x_2) \quad & \quad (x_1 + 1)^2
\end{bmatrix} \ = \ -\frac{2}{(x_1 + 1)^3}
\begin{bmatrix}
1-x_2 \\
x_1+1
\end{bmatrix}
\begin{bmatrix}
1-x_2 \\
x_1+1
\end{bmatrix}^{\top},
$$
which is negative semidefinite. Therefore, $\widehat{L}(x, y)$ is concave in $x$. In particular, restricting $\widehat{L}(x, y)$ on the line $x_1 + x_2 = 1$ yields $3 - (x_1 + 1) - \frac{1}{x_1+1}$, which is a concave function. In Figure~\ref{figure-concave}, we depict functions $\tilde{L}$ and $\widehat{L}$ (the black surfaces) as well as their restrictions on the line $x_1 + x_2 = 1$ (the intersecting curves of the gray and black surfaces).
\Halmos
\end{example}
Thanks to the concave representation, we can replace the \emph{non}-convex constraints $\theta \leq L(x, y)$ with \emph{convex} ones $\theta \leq \widehat{L}(x, y)$ when solving model \eqref{sl} in DCG. In Appendix \ref{apx:prop:soc}, we show that $\theta \leq \widehat{L}(x, y)$ are not only convex but also second-order conic representable in Proposition \ref{prop:soc} and present a detailed proof there. In implementation, we replace constraints $\theta \leq L(x, y)$ with a supporting hyperplane of $\widehat{L}(x, y)$ as linear cuts instead of using second-order conic constraints described in Proposition \ref{prop:soc} as cutting planes.  
Specifically, for given $(\hat{x}, \hat{y})$ in DCG, we incorporate the linear inequality
\begin{align}
& \theta \leq \widehat L (\hat x, \hat y) + \sum_{j \in J} {g_j {(\hat x, \hat y)} (x_j - \hat x_j)}, \label{ineq-bi}
\end{align}
where, for all $j \in J$,
\begin{equation*}
g_j (\hat{x},\hat{y}) := \frac{\partial \widehat L (x,\hat y)}{\partial x_j}\Big|_{x = \hat{x}} = \sum_{i \in I} {h_i \left(\frac{-w_{ik} (1-\hat y_k)Q}{P^2} + \frac{w_{ik}(-2\hat y_k \hat{x}_k+1+\hat y_k)}{P}\right)},
\end{equation*}
$P = U^{\tinyl}_i+U^{\tinyf}_i+\sum_{j \in J} {w_{ij}\bigl[(1-\hat{y}_j)\hat{x}_j+\hat{y}_j\bigr]}$, and $Q = U^{\tinyl}_i+\sum_{j \in J} {w_{ij} \bigl[-\hat{y}_j \hat{x}_j^2 +(1+\hat{y}_j) \hat{x}_j\bigr]}$. In what follows, we call \eqref{ineq-bi} the bulge inequalities. 

\subsection{Approximate Separation} \label{sec:cuts-approx}
We derive an approximate approach to solving the separation problem \eqref{sp}. Although this problem can be viewed as static CFLP, for which~\citet{Ljubic2018} has provided an exact solution approach, our goal is to solve it significantly faster in order to accelerate DCG, for which we demonstrate the resulting speedup numerically in Section \ref{sec:compu}. We start by identifying the following optimality conditions for Formulation (S-CFLP) and Model \eqref{sp}, respectively. 
\begin{lemma} \label{lem:oc}
There exists an optimal solution $x^*$ to (S-CFLP) such that $e^{\top}x^* = p$. In addition, for any fixed $\hat{x} \in \mathcal{X}$, there exists an optimal solution $y^*$ to the separation problem \eqref{sp} such that $e^{\top} y^* = r$.
\end{lemma}
A detailed proof of Lemma \ref{lem:oc} is presented in Appendix \ref{apx:lem:oc}. 
Following Lemma \ref{lem:oc}, the separation problem \eqref{sp} is equivalent to $\displaystyle \min_{y \in \overline{\mathcal{Y}}}L(\hat{x}, y)$ for given $\hat{x} \in \mathcal{X}$, where $\overline{\mathcal{Y}} := \{y \in \{0,1\}^J: e^{\top}y = r\}$. In what follows, we approximate the nonlinear function $L(\hat{x}, y)$ from above using a linear one. This leads to a relaxation of the separation problem \eqref{sp} that can be solved by a single-round sorting. 
\begin{proposition} \label{prop:approx}
For fixed $\hat{x} \in [0,1]^{|J|}$, define constants $\displaystyle a_i(\hat{x}) := U^{\tinyl}_i + \sum_{j \in J}w_{ij}\hat{x}_j$, $\displaystyle w^{\tinyl}_i(\hat{x}) := \min_{y \in \overline{\mathcal{Y}}} \Bigl\{U^{\tinyf}_i+\sum_{j\in J}w_{i j}(1-\hat{x}_j)y_j\Bigr\}$, and $\displaystyle w^{\tinyu}_i(\hat{x}) :=\max_{y \in \overline{\mathcal{Y}}} \Bigl\{U^{\tinyf}_i + \sum_{j\in J} w_{ij}(1-\hat{x}_j)y_j\Bigr\}$ for all $i \in I$. In addition, define
\begin{align*}
\alpha(\hat{x}) := \sum_{i\in I}h_i\left[\frac{a_i(\hat{x})\left(a_i(\hat{x})+w^{\tinyu}_i(\hat{x})+w^{\tinyl}_i(\hat{x})-U^{\tinyf}_i\right)}{\left(a_i(\hat{x})+w^{\tinyu}_i(\hat{x})\right)\left(a_i(\hat{x})+w^{\tinyl}_i(\hat{x})\right)}\right],
\end{align*}
and vector $\beta(\hat{x}) := [\beta_1(\hat{x}), \ldots, \beta_{|J|}(\hat{x})]^{\top}$ with
\begin{align*}
\beta_j(\hat{x}) := \sum_{i\in I}h_i \left[\frac{a_i(\hat{x})w_{ij}(1-\hat{x}_j)}{\left(a_i(\hat{x})+w^{\tinyu}_i(\hat{x})\right)\left(a_i(\hat{x})+w^{\tinyl}_i(\hat{x})\right)}\right]
\end{align*}
for all $j \in J$. Then, it holds that
\begin{align*}
    L(\hat{x}, y) \leq \alpha(\hat{x}) - \beta(\hat{x})^{\top}y, \quad \forall y \in \overline{\mathcal{Y}}.
\end{align*}
In addition, the problem $\displaystyle \min_{y \in \overline{\mathcal{Y}}}\bigl\{\alpha(\hat{x}) - \beta(\hat{x})^{\top}y\bigr\}$ admits an optimal solution $\hat{y} = \{[1], \ldots, [r]\}$, where $\{[j]: j \in J\}$ is a sorting of the set $J$ such that $\beta_{[1]}(\hat{x}) \geq \beta_{[2]}(\hat{x}) \geq \cdots \geq \beta_{[|J|]}(\hat{x})$.
\end{proposition}
We present the proof of Proposition \ref{prop:approx} in Appendix \ref{apx-prop-approx}. 
Notice that for given $\hat{x} \in [0,1]^{|J|}$, the values of $\alpha(\hat{x})$ and $\{\beta_j(\hat{x}): j \in J\}$ can all be computed in closed-form. In addition, Proposition \ref{prop:approx} suggests an approximate separation. Specifically, suppose that an incumbent solution $(\hat{x}, \hat{\theta})$, obtained from solving a relaxation of Model \eqref{sl}, satisfies  $\hat{\theta} > \alpha(\hat{x}) - \beta(\hat{x})^{\top} \hat{y}$, where $\hat{y}$ is described in Proposition~\ref{prop:approx}. Then, this solution violates the constraint $\theta \leq L(x, \hat{y})$, which can be added back to the relaxed formulation of \eqref{sl}. This approximate separation reduces the effort of solving a MINLP to a single-round sorting. As we shall report in Section \ref{sec:compu}, this leads to a substantial numerical speedup.

\subsection{A Branch-and-Cut Framework}
\label{sec:cuts:bcf}
In Algorithm \ref{algo-b&c}, we summarize the DCG algorithm for solving (S-CFLP), or equivalently its reformulation \eqref{sl}, in a branch-and-cut framework. This framework uses a single branching tree and maintains a set $\mathcal{F}$ of formulations with respect to the active tree nodes, i.e., the tree nodes that have not yet produced an integral solution or whose optimal value is not dominated by the best lower bound $\theta_{\text{LB}}$ found so far.  Note that in line \ref{alg:step} of Algorithm~\ref{algo-b&c}, we add the valid inequality globally to all active tree nodes.
\begin{algorithm}[ht!]
\SetAlgoLined
{\bf Initialization}: create a set $\mathcal{F}$ of formulations and insert the formulation $\displaystyle\max_{x, \theta}\left\{\theta: e^{\top}x = p, x \in [0,1]^{|J|}, \theta \in [0, 1]\right\}$, a continuous relaxation of \eqref{sl}, into $\mathcal{F}$\;
Set \texttt{BestSol}$\leftarrow \emptyset$ and $\theta_{\mbox{\tiny LB}} \leftarrow 0$\;
\While{$\mathcal{F}$ is non-empty}{
    Remove a formulation from $\mathcal{F}$\;
    Solve the formulation and obtain an incumbent solution $(\hat{x}, \hat{\theta})$ with optimal value $\hat{\theta}$\;
    Solve either the approximate separation problem $\displaystyle\min_{y \in \overline{\mathcal{Y}}}\{\alpha(\hat{x}) - \beta(\hat{x})^{\top}y\}$ or the exact separation problem \eqref{sp} and obtain an optimal solution $\hat{y}$\;
    \uIf{$\hat{\theta} > L(\hat{x}, \hat{y})$}
    {
        Add a cut $\theta \leq L(x, \hat{y})$ based on either the submodular inequality \eqref{7} or the bulge inequality \eqref{ineq-bi} to strengthen the formulation and all formulations in $\mathcal{F}$\;
        \label{alg:step}
        Insert the strengthened formulation back into $\mathcal{F}$\;
    }
    \uElseIf{$\hat{\theta} > \theta_{\mbox{\tiny LB}}$ and $\hat{x}$ is integral}
    {
        Set \texttt{BestSol}$\leftarrow \hat{x}$ and $\theta_{\mbox{\tiny LB}} \leftarrow \hat{\theta}$\;
    }
    \uElseIf{$\hat{\theta} > \theta_{\mbox{\tiny LB}}$ and $\hat{x}$ is fractional}{
        Branch on $\hat{x}$ and insert the two resulting formulations into $\mathcal{F}$\;
    }
}
\caption{A Branch-and-Cut Framework for Solving (S-CFLP)} \label{algo-b&c}
\end{algorithm}
Since there are a finite number of cuts $\theta \leq L(x, \hat{y})$ to add and a finite number of candidate solutions of $\hat{x}$, we make the following claim.
\begin{proposition}
Algorithm \ref{algo-b&c} terminates in a finite number of steps with a global optimal solution to (S-CFLP).
\end{proposition}

\subsection{Approximation Algorithm} \label{sec:heuristic}
In view of the computational challenges of solving (S-CFLP) to global optimum, we propose an approximate algorithm to more quickly obtain  good-quality solutions. Specifically, the proposed algorithm solves a single-level MISOCP, which can be readily solved by off-the-shelf solvers directly, waiving the need to design specialized solution methods like Algorithm \ref{algo-b&c}. Notably, this algorithm admits a constant approximation guarantee, as detailed in the following Theorem \ref{thm:heuristic}, of which we present a detailed proof in Appendix \ref{apx-thm-heuristic}. 
\begin{theorem} \label{thm:heuristic}
Let $z^*$ represent the optimal value of (S-CFLP), $x^{\mbox{\tiny H}}$ represent an optimal location decision to the following MISOCP
\begin{subequations}
\label{heu}
\begin{align}
\min_{\substack{x \in \mathcal{X}, \ \mu \geq 0,\\ s \geq 0, \ t \geq 0}}
\ & \ \sum_{i \in I} h_i U_i^{\tinyf} s_i + r\lambda + e^{\top}\mu \\
\rm{s.t.} \ & \ \lambda+\mu_j \geq \sum_{i \in I} h_i w_{ij}t_{ij}, \quad \forall j \in J,\\
& \begin{Vmatrix}
\begin{bmatrix} 
2  \\
U_i^{\tinyl} + \sum_{k \in J}w_{ik}x_k-s_i 
\end{bmatrix}
\end{Vmatrix} \leq s_i + U_i^{\tinyl} + \sum_{k \in J}w_{ik}x_k, \quad \forall i \in I, \label{heu-1c-ref}\\
& \begin{Vmatrix}
\begin{bmatrix} 
2(1-x_j)  \\
U_i^{\tinyl} + \sum_{k \in J} w_{ik}x_k-t_{ij} 
\end{bmatrix}
\end{Vmatrix} \leq t_{ij} + U_i^{\tinyl} + \sum_{k \in J} w_{ik}x_k, \quad \forall i \in I, j \in J, \label{heu-1d-ref}
\end{align}
\end{subequations}
and $z^{\mbox{\tiny H}}$ represents the objective function value of $x^{\mbox{\tiny H}}$ in (S-CFLP), i.e., $\displaystyle z^{\mbox{\tiny H}} := \min_{y \in \mathcal{Y}} L(x^{\mbox{\tiny H}}, y)$. Then, it holds that
$$
\frac{4\gamma_M \gamma_m}{(\gamma_M + \gamma_m)^2} \ \leq \ \frac{z^{\mbox{\tiny H}}}{z^*} \ \leq \ 1,
$$
where
\begin{equation*}
\gamma_m := \min_{i \in I}\left\{\frac{1}{1+\frac{\max_{y \in \overline{\mathcal{Y}}}\left\{U_i^{\tinyf}+\sum_{j \in J}w_{ij}y_j\right\}}{\min_{x \in \overline{\mathcal{X}}}\left\{U_i^{\tinyl}+\sum_{j \in J}w_{ij}x_j\right\}}}\right\}, \quad \gamma_M := \max_{i \in I}\left\{\frac{1}{1+\frac{\min_{y \in \overline{\mathcal{Y}}}\left\{U_i^{\tinyf}+\sum_{j \in J}w_{ij}y_j\right\}}{\max_{x \in \overline{\mathcal{X}}}\left\{U_i^{\tinyl}+\sum_{j \in J}w_{ij}x_j\right\}}}\right\}
\end{equation*}
with $\overline{\mathcal{X}} := \left\{x \in \{0,1\}^{|J|}: e^{\top}x = p\right\}$ and $\overline{\mathcal{Y}} = \left\{y \in \{0,1\}^{|J|}: e^{\top}y = r\right\}$.
\end{theorem}
\begin{remark}
Constants $\gamma_M$ and $\gamma_m$ can be computed in closed-form. Hence, Theorem~\ref{thm:heuristic} presents a constant approximation algorithm for (S-CFLP) via solving a MISOCP. In theory, it is impossible to improve this result (i.e., to a polynomial-time approximation algorithm) unless $\mbox{P} = \mbox{NP}$, in view of the inapproximability conclusion of Theorem \ref{thm:hardness}. Nevertheless, The MISOCP is practically tractable thanks to promising performance of state-of-the-art solvers. We can further improve the efficacy of solving formulation \eqref{heu} by adding valid inequalities similar to those presented in Section \ref{sec:cuts}. We present the details of these inequalities in Appendix \ref{apx-vi-misocp}.
\end{remark}

\section{Computational Results}
\label{sec:compu}
We test a variety of sequential CFLP instances to validate the efficacy of our approaches on obtaining optimal or approximate solutions to (S-CFLP). We describe parameter configurations and experimental design in Section \ref{sec:instance} and conduct result analysis to demonstrate that
\begin{enumerate}[(i)]
\item the submodular inequalities \eqref{7} and bulge inequalities \eqref{ineq-bi} significantly accelerate the branch-and-cut framework for solving (S-CFLP) in Section \ref{sec:perform};
\item the approximate separation provides a further speed-up in Appendix~\ref{sec:compu:appr};
\item the approximation algorithm can quickly find a good-quality solution in Appendix~\ref{sec:compu:socp};
\item the patterns of the leader's and follower's optimal locations highly depends on settings of the choice model, generating insights on winning market share in Section \ref{sec:compu:sens}. 
\end{enumerate}

\subsection{Experimental Design} \label{sec:instance}
All the code is written in C++ and we solve all mixed-integer programs in CPLEX 12.6 using the default configurations of the solver.  All numerical tests are run on a PC with Intel CORE (TM) i7-8550 1.8GHz CPU, 16G RAM running 64-bit Windows 10. When implementing Algorithm \ref{algo-b&c}, we generate valid inequalities derived in Sections \ref{sec:cuts-sc} and \ref{sec:cuts-oa} via the \texttt{lazy callback} function in CPLEX, which is called upon when a binary-valued solution is found at a branching node. (In a separate implementation (not reported in this paper), we tried separating the submodular inequalities~\eqref{7} for incumbent solutions $\hat{x}$ that are fractional-valued, but the computational efficacy is inferior to that obtained by separating~\eqref{7} for binary-valued solutions only.) All the code and instances are available at \url{https://github.com/MingyaoQi/S-CFLP.git}.

To the best of our knowledge, there are no immediately available benchmark instances having comparable sizes to what we aim to solve in the existing literature. Therefore, we generate our instances following the settings in the sequential CFLP literature including  \citet{Kucukaydn2012,Haase2014,Gentile2018}. In specific, we consider a $[0, 50] \times [0, 50]$ square on a planar surface, in which the locations of demand points and candidate facility sites are randomly generated with integer coordinates. We employ the Euclidean metric to compute distances and use a MNL model with the default $\beta := 0.1$ to generate probabilistic utilities. In the default case, we set $\alpha_j := 0$ for all $j \in J$ to generate homogeneous attractiveness among facilities. In Section~\ref{sec:compu:sens}, we vary the values of $\beta$ and test heterogeneous $\alpha_j$-values to observe their impacts on the optimal locations selected by the leader and the follower. We assume that $J^{\tinyl} = J^{\tinyf} = \emptyset$ as the case where both the leader and the follower do not own any pre-existing facilities. 

\begin{table}[ht!]
\fontsize{9}{15}\selectfont
  \centering
  \caption{Parameter settings and instant sizes in existing literature}
  \label{table 2}
   \resizebox{0.5\textwidth}{!}{
    \begin{tabular}{lrrrr}
    \toprule
    Reference & $\mid I\mid$ &$\mid J\mid$ &$p$&$r$ \\
    \midrule
   \citet{Kucukaydn2011}  & 30    & 5     & $\ast$     & 0 \\
    \citet{Kucukaydn2012} & 100   & 20      & $\ast$     & $\ast$ \\
   \citet{Roboredo2013} & 100   & 50        & 4     & 4 \\
   \citet{Alekseeva2015} & 100   & 100       & 20    & 20 \\
    \citet{Gentile2018} & 225   & 225      & 5     & 5 \\
    \bottomrule
    \end{tabular}%
    }
\end{table}%
Before introducing our results, we briefly review the scale of experiments conducted by previous studies on sequential CFLP in Table \ref{table 2}, where entries $\ast$ indicate that the corresponding paper also optimized $p$ or $r$, in addition to optimizing locations. Note that except for \citet{Kucukaydn2011}, where only five candidate sites are considered, none of the previous studies can guarantee solution optimum. In this paper, we take up the challenge of solving problems with up to 100 facility sites and 2000 customer nodes to global optimum. 

\subsection{Strength of Valid Inequalities} \label{sec:perform}
We perform three implementations of Algorithm \ref{algo-b&c}, by adding the submodular inequalities \eqref{7} only, the bulge inequalities \eqref{ineq-bi} only, and both, denoted by SC, BI, and SCBI, respectively. Since no exact solution approach exists for (S-CFLP) in the literature, to benchmark these implementations, we solve each instance by an enumeration method, which finds optimal leader locations by evaluating the objective function value $\min_{y \in \mathcal{Y}}L(x, y)$ for all $x \in \mathcal{X}$. We notice that this formulation, as well as the separation problem~\eqref{sp}, is equivalent to a static CFLP model and then we can solve it via the state-of-the-art approach from~\citet{Ljubic2018} (specifically, using their outer approximation cuts), both in the enumeration method and in Algorithm~\ref{algo-b&c}. 

\begin{figure}[ht!]
\centering
\includegraphics[width=0.5\textwidth]{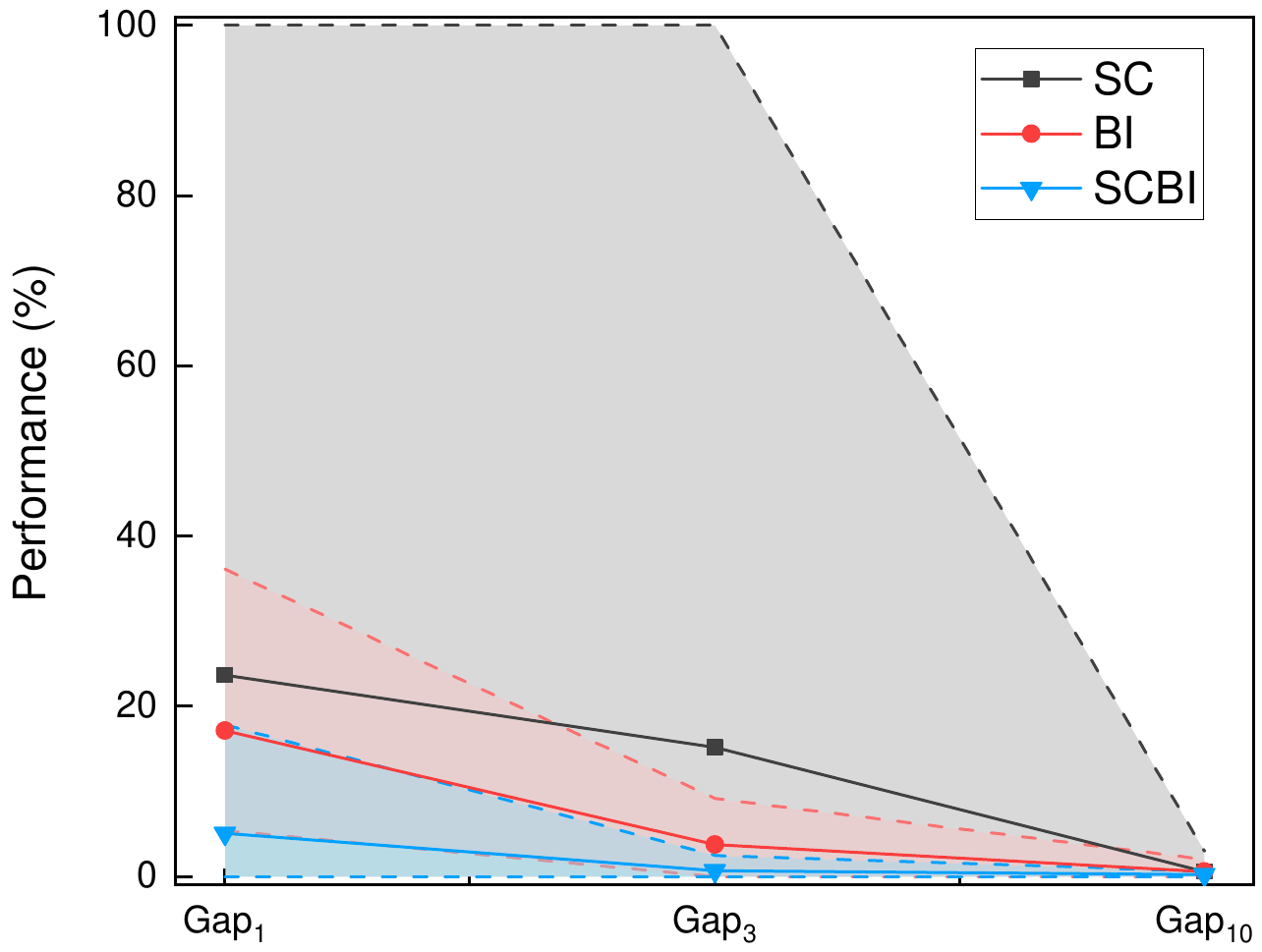}
\caption{Optimality gap improvements by different valid inequalities after adding the first, third, and tenth cut in the branch-and-cut framework (optimality gap equals the relative difference between the final optimal value and the objective value of the best integer solution found so far).} \label{fig:gaps}
\end{figure}

In the first group of instances, the number of candidate facilities (i.e., $|J|$) ranges from 20 to 100 and the number of customers (i.e., $|I|$) is set to be equal to the number of candidate facilities. The values of $p$ and $r$ vary as 2 or 3. Each instance is indicated by its scale $|I|$-$|J|$-$p$-$r$. We report the computational results in Table \ref{table 3} and Figure \ref{fig:gaps}. 
In Table \ref{table 3}, ``\texttt{Time(s)}'' report the CPU seconds it takes to solve each instance to optimum and we highlight the shortest solution time for each instance in bold. We report ``LIMIT'' if an instance is not solved to optimum within 4 hours. In addition, we report the total number of valid inequalities added until an optimal solution is found and the size of the final branch-and-bound tree in columns ``\texttt{\#Cuts}'' and ``\texttt{\#Nodes},'' respectively. We record the gaps between the final optimal objective value and the objective value of the best integer solution found after adding the first, third and tenth cut in CPLEX, denoted by ``\texttt{Gap$_1$}'', ``\texttt{Gap$_3$}'', and ``\texttt{Gap$_{10}$}'', respectively. Note that for each round of callback, at most one cut is added in SC and in BI, and at most two cuts are added in SCBI. In Figure \ref{fig:gaps}, we report the average of the gaps (solid line), as well as the second smallest and second largest gaps (shaded area) among all instances (see the detailed results of these gaps in Table \ref{tab:fig-gaps} in Appendix \ref{apx-fig-gaps}).

% Table generated by Excel2LaTeX from sheet 'Table3'
\begin{table}[ht!]
\fontsize{9}{11}\selectfont
  \centering
  \caption{Effectiveness of valid inequalities used in  Algorithm \ref{algo-b&c} for optimizing (S-CFLP) instances}
    \resizebox{\textwidth}{!}{
    \begin{tabular}{crrrrrrrrrr}
    \hline
    \multirow{2}[2]{*}{Instance} & \multicolumn{3}{c}{SC} & \multicolumn{3}{c}{BI} & \multicolumn{3}{c}{SCBI} & \multicolumn{1}{c}{Enumeration} \bigstrut[t]\\
          & \multicolumn{1}{r}{\texttt{Time(s)}} & \multicolumn{1}{r}{\texttt{\#Cuts}} & \multicolumn{1}{r}{\texttt{\#Nodes}} & \multicolumn{1}{r}{\texttt{Time(s)}} & \multicolumn{1}{r}{\texttt{\#Cuts}} & \multicolumn{1}{r}{\texttt{\#Nodes}} & \multicolumn{1}{r}{\texttt{Time(s)}} & \multicolumn{1}{r}{\texttt{\#Cuts}} & \multicolumn{1}{r}{\texttt{\#Nodes}} & \multicolumn{1}{r}{\texttt{Time(s)}} \\
    \hline
    \multicolumn{1}{l}{20-20-2-2} & 0.94  & 38    & 129   & 2.14  & 33    & 88    & \textbf{0.75} & 50    & 110   & 3.47 \bigstrut[t]\\
    \multicolumn{1}{l}{20-20-3-2} & 1.66  & 76    & 450   & 3.09  & 59    & 361   & \textbf{0.97} & 83    & 311   & 19.58 \\
    \multicolumn{1}{l}{20-20-2-3} & 2.00  & 58    & 151   & 2.42  & 44    & 117   & \textbf{1.30} & 73    & 120   & 5.48 \\
    \multicolumn{1}{l}{40-40-2-2} & 13.23 & 220   & 839   & 4.06  & 56    & 576   & \textbf{3.75} & 114   & 514   & 46.42 \\
    \multicolumn{1}{l}{40-40-3-2} & 68.09 & 1192 & 6339 & 15.28 & 196   & 2155  & \textbf{11.16} & 369   & 2323  & 496.77 \\
    \multicolumn{1}{l}{40-40-2-3} & 64.44 & 335   & 964   & 20.88 & 60    & 448   & \textbf{11.02} & 115   & 538   & 146.52 \\
    \multicolumn{1}{l}{60-60-2-2} & 69.95 & 514   & 1694 & 29.36 & 116   & 1172  & \textbf{9.67} & 124   & 1527  & 220.33 \\
    \multicolumn{1}{l}{60-60-3-2} & 777.56 & 5631 & 33312 & 79.88 & 334   & 4006  & \textbf{39.33} & 549   & 11678 & 3630.78 \\
    \multicolumn{1}{l}{60-60-2-3} & 640.11 & 755   & 2290 & 211.22 & 121   & 1268  & \textbf{94.92} & 156   & 1462  & 1399.27 \\
    \multicolumn{1}{l}{80-80-2-2} & 353.55 & 1240 & 3601 & 65.49 & 122   & 2119  & \textbf{25.75} & 175   & 3134  & 817.75 \\
    \multicolumn{1}{l}{80-80-3-2} & 13655.10 & 15236 & 93558 & 147.78 & 345   & 10278 & \textbf{146.78} & 941   & 24219 & LIMIT \\
    \multicolumn{1}{l}{80-80-2-3} & 5181.42 & 1538 & 4480 & 384.99 & 142   & 2421  & \textbf{228.08} & 207   & 3387  & 6989.31 \\
    \multicolumn{1}{l}{100-100-2-2} & 636.63 & 1573 & 5741 & 57.97 & 89    & 2834  & \textbf{44.95} & 176   & 3656  & 2087.59 \\
    \multicolumn{1}{l}{100-100-3-2} & 13418.00 & 22628 & 155384 & 233.02 & 323   & 7972  & \textbf{190.53} & 772   & 33046 & LIMIT \\
    \multicolumn{1}{l}{100-100-2-3} & 5469.91 & 2143 & 6943 & 384.00 & 105   & 3124  & \textbf{273.86} & 194   & 4066  & LIMIT \\
    \hline
    \textbf{Average} & 2690.17 & 3545 & 21058 & 109.44 & 143   & 2596  & \textbf{72.19} & 273   & 6006  & N/A \\
    \hline
    \end{tabular}%
    }
  \label{table 3}%
\end{table}%

From Table \ref{table 3} and Figure \ref{fig:gaps}, we observe the following about the submodular and bulge inequalities.
\begin{enumerate}
\item Both inequalities strengthen the formulation significantly. For example, all three implementations outperform the benchmark method in all instances. In 10 rounds of cuts, all implementations are able to prove an optimality gap of below 3\% in all instances (below 1\% on average). In particular, the incorporation of both inequalities (i.e., the SCBI implementation) exhibits the best strength. For example, SCBI proves an average optimality gap of 0.65\% and 0.20\% across all instances in three and ten cuts, respectively. The small gaps suggest that Algorithm~\ref{algo-b&c} has the potential of solving much larger-scale S-CFLP instances -- even though it may fail to prove optimality, the best integer solutions found by the time limit may be of high quality. 

\item The submodular inequalities are less effective than the bulge inequalities. (Note that, however, this is not the case in Table \ref{table 4}, which we will discuss later.) For example, BI solves each instance to optimum within 600 seconds. In contrast, SC spends 2690.17 seconds on average to solve an instance, which is roughly 25 times that of BI. In addition, SC involves more cuts and a larger branch-and-bound tree, which are roughly 25 and 8 times those of BI, respectively.

\item SCBI performs even better than BI. For example, the average solution time by SCBI is about $34\%$ less than that by BI. This suggests that the submodular and bulge inequalities complement each other quite well.
\end{enumerate}

Next, we increase $p$ and $r$ while fixing $|I| = |J|$ at 20 and 30, respectively. The search space of the (S-CFLP) problem has a cardinality $C(|J|,p) \times C(|J|-p,r)=\dbinom{|J|}{p} \times \dbinom{|J|-p}{r}$. Thus, (S-CFLP)'s search space increases exponentially with $p$ and $r$ when $p\le |J|/2$ and $r\le(|J|-p)/2$ (see a proof of the exponential increase in Appendix~\ref{apx-prop:exponential}). This makes (S-CFLP) extremely difficult to solve. We report the computational results in Tables \ref{table 4}--\ref{table 5}. The instances that were not solved to optimum within 1 hour are marked as ``LIMIT,'' and if an implementation was not able to find any feasible solution in ten rounds of cuts, we report ``N/A'' for $\texttt{Gap}_{10}$. The average metrics in Tables \ref{table 4}--\ref{table 5} (reported in the bottom rows) are calculated among the instances that were solved to optimum by all three implementations.

\begin{table}[ht!]
\fontsize{9}{11}\selectfont
  \centering
  \caption{Computational results of instances with $|I|=|J|=20$ and varying $p$-, $r$-values}
  \label{table 4}
  \resizebox{0.87\textwidth}{!}{
    \begin{tabular}{lrrrrrrrrrr}
    \hline
    \multirow{2}[2]{*}{Instance} & \multicolumn{3}{c}{SC} & \multicolumn{3}{c}{BI} & \multicolumn{3}{c}{SCBI} & \multicolumn{1}{c}{{Enumeration}} \bigstrut[t]\\
          & \multicolumn{1}{c}{\texttt{Time(s)}} & \multicolumn{1}{c}{\texttt{Gap$_{10}$}} & \multicolumn{1}{c}{\texttt{\#Cuts}} & \multicolumn{1}{c}{\texttt{Time(s)}} & \multicolumn{1}{c}{\texttt{Gap$_{10}$}} & \multicolumn{1}{c}{\texttt{\#Cuts}} & \multicolumn{1}{c}{\texttt{Time(s)}} & \multicolumn{1}{c}{\texttt{Gap$_{10}$}} & \multicolumn{1}{c}{\texttt{\#Cuts}} & \multicolumn{1}{c}{\texttt{Time(s)}} \\
    \hline
    \multicolumn{1}{l}{20-20-2-2} & 0.97  & 0.00\% & 38    & {1.24}  & 0.00\% & 33    & \textbf{0.56} & 0.00\% & 50 & {6.09} \bigstrut[t]\\
    \multicolumn{1}{l}{20-20-4-2} & 2.70  & 0.00\% & 140   & 2.03  & 0.00\% & 95    & \textbf{1.16} & 0.00\% & 117 & {126.38}\\
    \multicolumn{1}{l}{20-20-6-2} & 6.61  & 1.49\% & 350   & \textbf{1.69} & 0.18\% & 90    & 1.75  & 0.23\% & 201 & {547.59}\\
    \multicolumn{1}{l}{20-20-8-2} & 3.58  & 2.58\% & 194   & 2.45  & 0.00\% & 150   & \textbf{1.39} & 3.04\% & 170 & {1677.17}\\
    \multicolumn{1}{l}{20-20-10-2} & 2.67  & 0.00\% & 132   & 1.91  & 0.00\% & 121   & \textbf{1.44} & 0.00\% & 182 & {2327.67}\\
    \multicolumn{1}{l}{20-20-2-4} & 3.09  & 0.00\% & 70    & 1.67  & 0.00\% & 45    & \textbf{1.34} & 0.00\% & 74 & {7.42}\\
    \multicolumn{1}{l}{20-20-4-4} & 9.14  & 4.23\% & 299   & 4.19  & 0.18\% & 175   & \textbf{3.88} & 4.20\% & 310 & {133.50}\\
    \multicolumn{1}{l}{20-20-6-4} & 13.27  & 2.36\% & 504   & 6.83  & 0.00\% & 353   & \textbf{6.05} & 3.81\% & 578 & {812.38}\\
    \multicolumn{1}{l}{20-20-8-4} & 17.11  & 0.00\% & 711   & 12.70 & 0.00\% & 781   & \textbf{9.75} & 0.00\% & 953 & {2119.78}\\
    \multicolumn{1}{l}{20-20-10-4} & \textbf{17.45} & 1.93\% & 881   & 21.53 & 0.56\% & 1206  & 34.64 & 1.93\% & 1552 & {2729.45}\\
    \multicolumn{1}{l}{20-20-2-6} & 2.73  & 2.08\% & 75    & 5.70  & 0.00\% & 46    & \textbf{1.75} & 0.00\% & 82 & {6.53}\\
    \multicolumn{1}{l}{20-20-4-6} & 11.36  & 5.64\% & 395   & 17.39 & 0.00\% & 292   & \textbf{8.52} & 6.49\% & 497 & {142.66}\\
    \multicolumn{1}{l}{20-20-6-6} & \textbf{24.50} & 6.46\% & 981   & 28.88 & 0.00\% & 1218  & 27.06 & 4.04\% & 1441 & {844.77}\\
    \multicolumn{1}{l}{20-20-8-6} & \textbf{43.13} & 2.52\% & 1909  & 80.33 & 0.34\% & 3453  & 57.06 & 2.52\% & 3464 & {2177.16}\\
    \multicolumn{1}{l}{20-20-10-6} & \textbf{50.69} & 0.93\% & 2285  & 137.94 & 0.00\% & 5885  & 69.58 & 1.07\% & 4399 & {2782.27}\\
    \multicolumn{1}{l}{20-20-2-8} & 2.64  & 1.68\% & 81    & 2.06  & 0.00\% & 58    & \textbf{1.58} & 1.68\% & 123 & {5.61}\\
    \multicolumn{1}{l}{20-20-4-8} & 12.70  & 9.36\% & 519   & 18.14 & 0.21\% & 620   & \textbf{8.30} & 9.36\% & 821 & {118.39}\\
    \multicolumn{1}{l}{20-20-6-8} & 43.05  & 5.28\% & 1867  & 87.31 & 0.00\% & 3308  & \textbf{33.94} & 5.40\% & 3330 & {766.03}\\
    \multicolumn{1}{l}{20-20-8-8} & \textbf{68.17} & 2.54\% & 3067  & 303.06 & 0.00\% & 9930  & 69.72 & 3.46\% & 5679 & {1992.34}\\
    \multicolumn{1}{l}{20-20-10-8} & \textbf{61.42} & 0.90\% & 2526  & 1370.64 & 0.00\% & 28118 & 81.05 & 0.00\% & 4857 & {2608.61}\\
    \multicolumn{1}{l}{20-20-2-10} & 2.24  & 0.00\% & 89    & 1.69  & 0.00\% & 77    & \textbf{1.66} & 0.00\% & 152 & {5.45}\\
    \multicolumn{1}{l}{20-20-4-10} & 15.17  & 3.83\% & 728   & 19.59 & 0.18\% & 1065  & \textbf{11.14} & 0.18\% & 1236 & {103.17}\\
    \multicolumn{1}{l}{20-20-6-10} & 45.11  & 7.26\% & 2133  & 134.98 & 0.00\% & 6677  & \textbf{40.56} & 4.27\% & 3962 & {700.84}\\
    \multicolumn{1}{l}{20-20-8-10} & \textbf{79.94} & 4.95\% & 3216  & 1619.81 & 0.00\% & 28190 & 92.58 & 4.95\% & 6145 & {2030.36}\\
    \multicolumn{1}{l}{20-20-10-10} & \textbf{74.02} & 0.00\% & 2964  & 8493.41  & 0.00\% & 33664 & 87.84 & 0.00\% & 5669 & {2495.48}\\
    \hline
    \textbf{Average} & \textbf{24.54} & 2.64\% & 1046  & {495.09}  & 0.07\% & 5026  & 26.17 & 2.27\% & 1842 & {1090.68} \\
    \hline
    \end{tabular}%
    }
\end{table}%

\begin{table}[ht!]
\fontsize{9}{11}\selectfont
  \centering
  \caption{Computational results of instances with $|I|=|J|=30$ and varying $p$-, $r$-values}
  \label{table 5}
  \resizebox{0.9\textwidth}{!}{
    \begin{tabular}{lrrrrrrrrr}
    \hline
    \multicolumn{1}{c}{\multirow{2}[2]{*}{Instance}} & \multicolumn{3}{c}{SC} & \multicolumn{3}{c}{BI} & \multicolumn{3}{c}{SCBI} \bigstrut[t]\\
          & \multicolumn{1}{c}{\texttt{Time(s)}} & \multicolumn{1}{c}{\texttt{Gap$_{10}$}} & \multicolumn{1}{c}{\texttt{\#Cuts}} & \multicolumn{1}{c}{\texttt{Time(s)}} & \multicolumn{1}{c}{\texttt{Gap$_{10}$}} & \multicolumn{1}{c}{\texttt{\#Cuts}} & \multicolumn{1}{c}{\texttt{Time(s)}} & \multicolumn{1}{c}{\texttt{Gap$_{10}$}} & \multicolumn{1}{c}{\texttt{\#Cuts}} \\
    \hline
    30-30-3-3 & 34.11 & 0.99\% & 567   & 8.33  & 0.17\% & 141   & \textbf{6.06} & 3.10\% & 214 \bigstrut[t]\\
    30-30-6-3 & 450.97 & 2.45\% & 5149  & 21.84 & 0.13\% & 609   & \textbf{19.16} & 1.84\% & 1021 \\
    30-30-9-3 & LIMIT & 2.50\% & 11842 & 36.02 & 0.58\% & 1361  & \textbf{34.19} & 1.62\% & 2047 \\
    30-30-12-3 & 2198.05 & 2.45\% & 3723  & 27.42 & 5.39\% & 1245  & \textbf{20.20} & 2.45\% & 1331 \\
    30-30-15-3 & 88.42 & 2.07\% & 1130  & 28.59 & 7.19\% & 1510  & \textbf{20.49} & 2.07\% & 1409 \\
    30-30-3-6 & 143.83 & 0.17\% & 825   & \textbf{33.52} & 0.17\% & 245   & 34.53 & 0.17\% & 493 \\
    30-30-6-6 & 3259.80 & 1.45\% & 12203 & 140.97 & 0.17\% & 2446  & \textbf{130.28} & 1.45\% & 4425 \\
    30-30-9-6 & LIMIT & N/A   & 14238 & \textbf{269.45} & 0.67\% & 6131  & 400.83 & 0.69\% & 9737 \\
    30-30-12-6 & LIMIT & N/A   & 14453 & \textbf{1837.03} & 4.30\% & 18120 & 1981.98 & 0.25\% & 17798 \\
    30-30-3-9 & 92.06 & 2.36\% & 897   & 38.48 & 6.34\% & 372   & \textbf{31.13} & 4.09\% & 697 \\
    30-30-6-9 & 2098.56 & 0.37\% & 12108 & 314.88 & 0.00\% & 5487  & \textbf{304.11} & 0.37\% & 8989 \\
    30-30-3-12 & 55.72 & 1.09\% & 919   & 25.36 & 0.23\% & 493   & \textbf{20.95} & 3.88\% & 898 \\
    30-30-6-12 & 2199.31 & 1.76\% & 14509 & 1026.81 & 1.52\% & 14695 & \textbf{965.25} & 1.76\% & 17479 \\
    30-30-3-15 & 39.92 & 4.01\% & 968   & 22.73 & 0.88\% & 715   & \textbf{21.39} & 4.01\% & 1156 \\
    30-30-6-15 & 2609.08 & 2.33\% & 17089 & LIMIT & N/A   & N/A   & \textbf{2234.25} & 2.33\% & 25750 \\
    \hline
    \textbf{Average} & {969.16} & {1.74\%} & {4818}  & {153.54} & {2.02\%} & {2542}  & {\textbf{143.05}} & {2.29\%} & {3465}\\
    \hline
    \end{tabular}%
    }
\end{table}%

Shown in Table \ref{table 4}, BI is less efficient than SC, different from the results in Tables \ref{table 3} and \ref{table 5}. From Tables \ref{table 3}--\ref{table 5}, we observe that SCBI remains the most competitive implementation among the three alternatives. For example, the solution time of SCBI is the shortest in a majority of instances, and in the instances this is not the case, SCBI performs comparably with the best implementation. In contrast, the solution time of SC and BI increases quickly as $p$ and $r$ increase (see BI in Table \ref{table 4} and SC in Table \ref{table 5}), reaching the time limit in several instances. This indicates that the submodular and bulge inequalities still complement each other well in these more challenging instances. For this reason, we adopt SCBI as the benchmark approach in all subsequent experiments. 

\subsection{Sensitivity Analysis}
\label{sec:compu:sens}
We analyze solution time and results of solving (S-CFLP) under various parameter settings. The results will provide insights on winning market share and how to conduct parameter selection for firms with leader and follower roles in sequential CFLP. We present the results for varying coefficients $\beta$ and $\alpha_j$ of the choice model in Sections~\ref{sec:compu:beta} and \ref{sec:compu:alpha}, respectively, and for varying customer sizes $|I|$ in Appendix~\ref{apx-sec:vary-customer-size-I}.

\subsubsection{Varying $\beta$ in the Choice Model} \label{sec:compu:beta}
\begin{figure}[ht!]
\centering
\includegraphics[width=\textwidth]{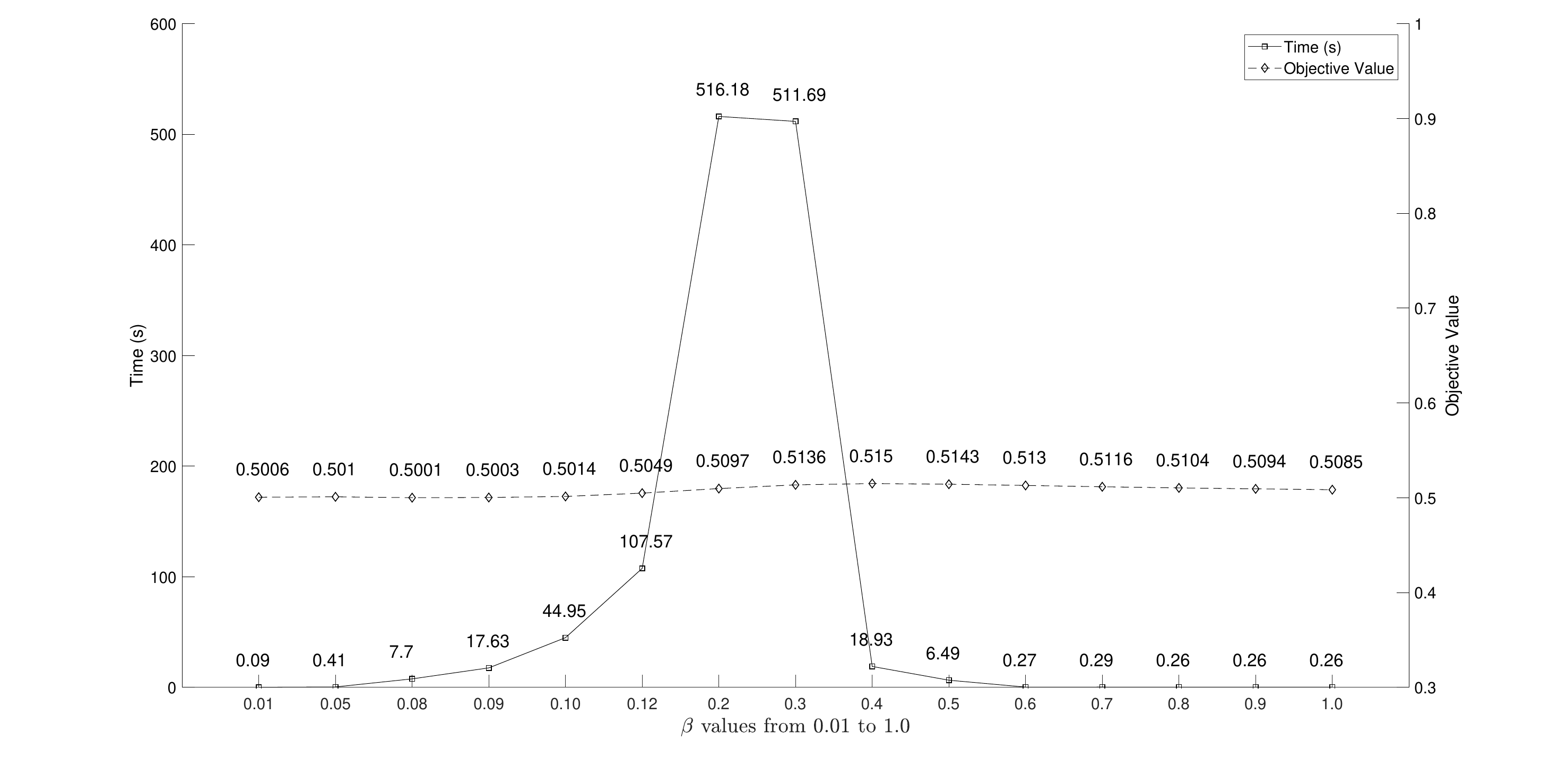}
\caption{Solution time and optimal objective values given varying  $\beta = 0.01$ to $1.0$.} 
\label{figure 3}
\end{figure}

\begin{figure}[ht!]
\centering
\subfigure[{$\beta=0.05$}]{
\includegraphics[width=0.45\linewidth]{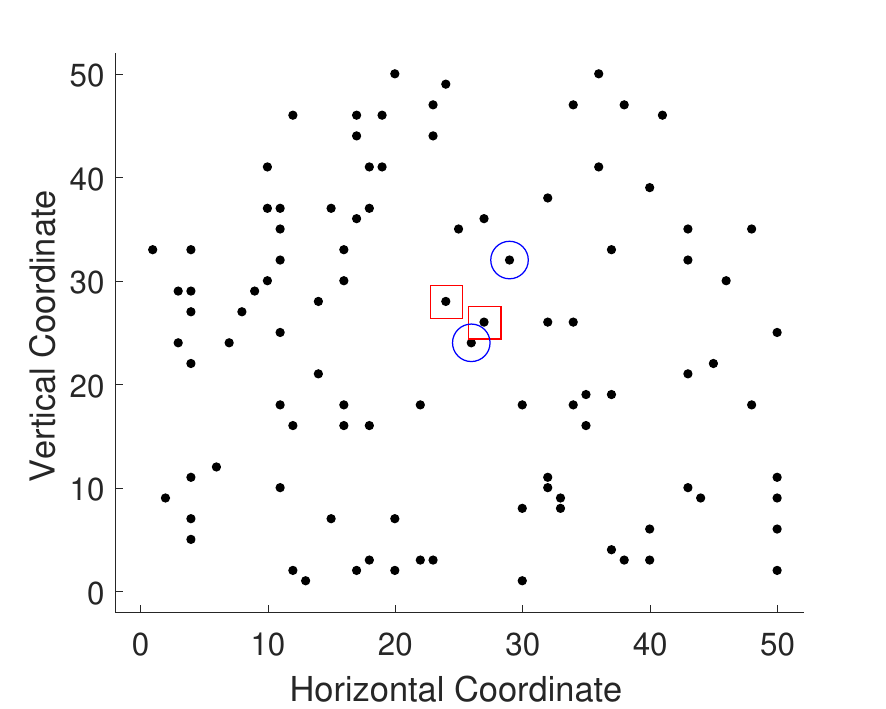}
}
\subfigure[{$\beta=0.08$}]{
\includegraphics[width=0.45\linewidth]{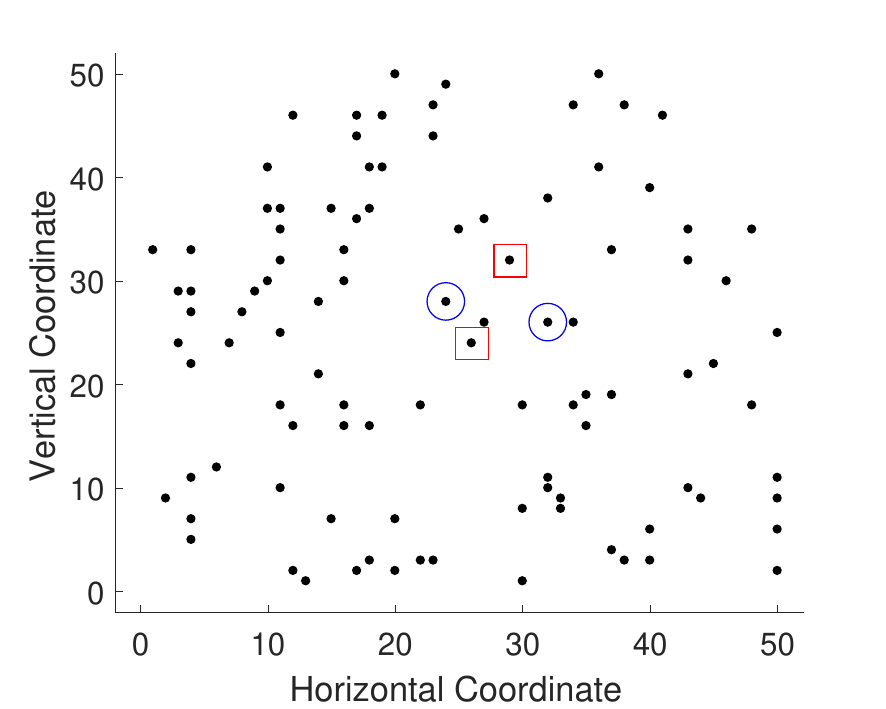}
}
\subfigure[{$\beta=0.1$}]{
\includegraphics[width=0.45\linewidth]{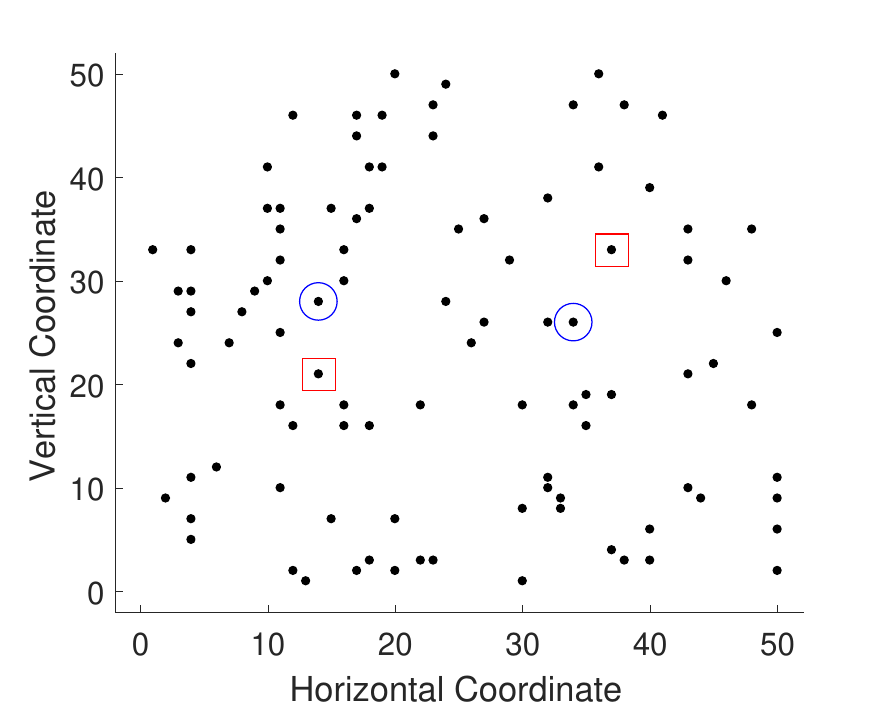}
}
\subfigure[{$\beta=0.2$}]{
\includegraphics[width=0.45\linewidth]{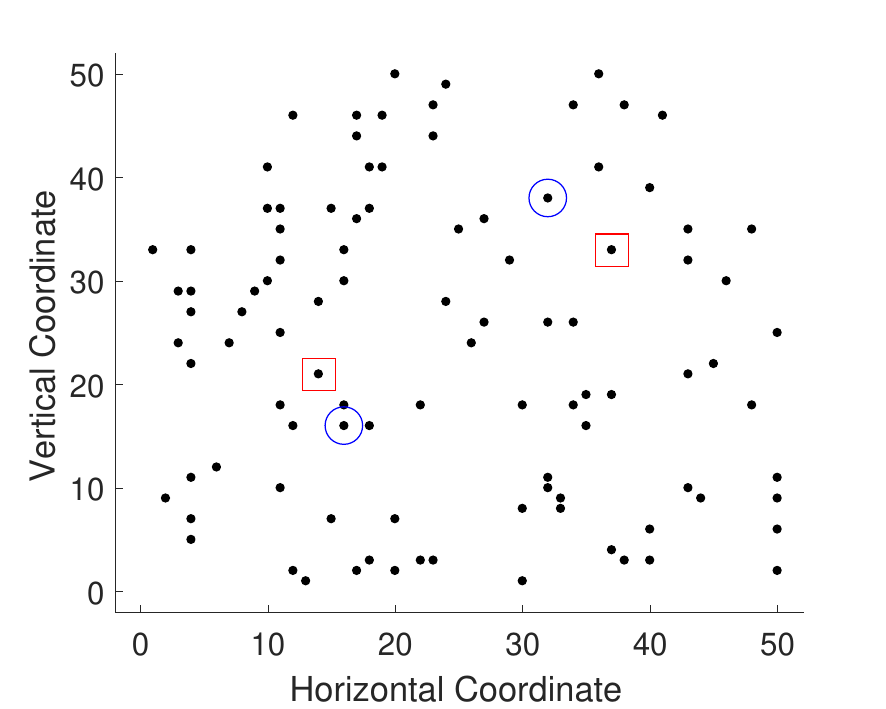}}
\subfigure{
\includegraphics[width=0.45\linewidth]{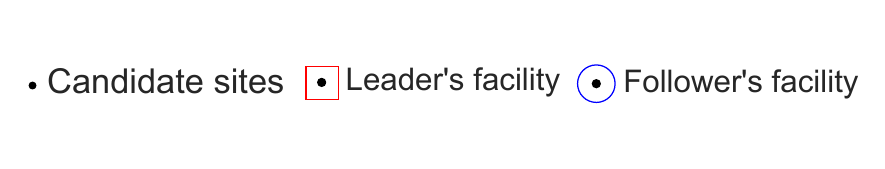}
}
\caption{Leader's and follower's optimal location decisions given varying $\beta$.}
\label{figure 4}
\end{figure}

Recall that in the MNL model, parameter $\beta$ measures the impact of traveling distance $d_{ij}$ on the utility of a customer $i$ patronizing facility $j$, i.e., $w_{ij}=\mbox{exp}\{\alpha_j-\beta d_{ij}\}$. We examine how the optimal objective value and solution time of (S-CFLP) change under various choices of $\beta$. In specific, we implement SCBI on the instance 100-100-2-2 with $\beta$ ranging from 0.01 to 1.0. We report the resulting optimal objective values and solution time in Figure \ref{figure 3} and the optimal locations of some representative instances in Figure \ref{figure 4}.

From Figure \ref{figure 3}, the solution time is short for most $\beta$-choices, but it becomes extremely long when $\beta$ ranges between 0.2 and 0.3. The optimal objective value (as the leader's total market share) is quite stable around 50\% for different $\beta$ choices. That is, given that the leader and the follower each can open at most two facilities (i.e., $p=r=2$), the optimal solutions to (S-CFLP) make sure that they split the market almost equally, regardless whether or not customers are more or less willing to travel for patronizing. 

Figure \ref{figure 4} depicts optimal location choices by the leader and the follower for $\beta=0.05$, $0.08$, $0.1$, and $0.2$. First, note that the optimal locations are clustered when $\beta$ is small (e.g., $\beta = 0.05$) and then spread out when $\beta$ increases. Indeed, a small $\beta$ implies lower spatial impedance effect, by which the facilities tend to be located at the center of the region to attract customers from all directions. As $\beta$ increases, customers become more likely to patronize nearby facilities according to the MNL model. As a result, it becomes optimal (for both leader and follower) to spread out the facilities in order to avoid self-competition and cover as many customers as possible. When $\beta$ is sufficiently large, the location results remain the same because customer $i$ would almost only visit its nearest facility $j$ yielding the dominant utility $w_{ij}$.
In that case, since customers only patronize locally, the follower harvests more customers by locating its facilities farther away from the leader's. This observation is particularly relevant when customers traveling for shopping becomes inconvenient (e.g., in challenging weather) or risky (e.g., during a pandemic). 
We also observe that the follower tends to locate its facilities near the leader's, demonstrating the economies of agglomeration. 

\subsubsection{Heterogeneous $\alpha_j$ and Impacts} \label{sec:compu:alpha}
We select Figures \ref{figure 4}(a) and \ref{figure 4}(c) as two representative cases with $\beta = 0.05$ and $\beta=0.1$, respectively, and vary the $\alpha_j$-values for different locations $j \in J$ to examine how heterogeneous attractiveness levels will affect the location choices under the same $\beta$ (i.e., customers' location preferences given by the distance factor remain the same).

\begin{figure}[ht!]
\centering
\subfigure[{$\alpha^{\textrm{Left}} = 0, \ \alpha^{\textrm{Right}}=0$}]{
\includegraphics[width=0.45\linewidth]{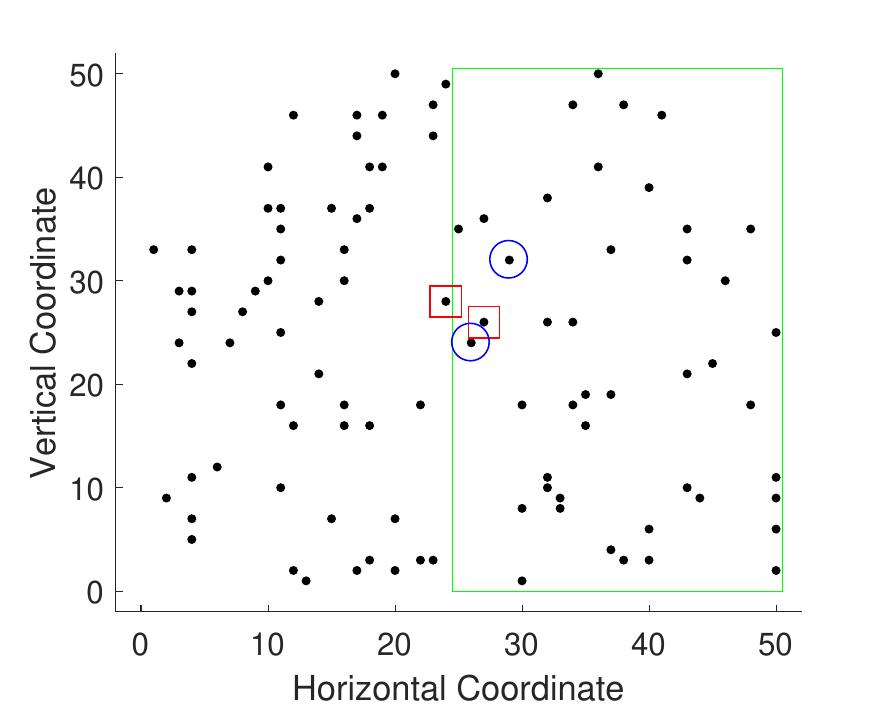}
}
\subfigure[{$\alpha^{\textrm{Left}} = 0, \ \alpha^{\textrm{Right}}=0.01$}]{
\includegraphics[width=0.45\linewidth]{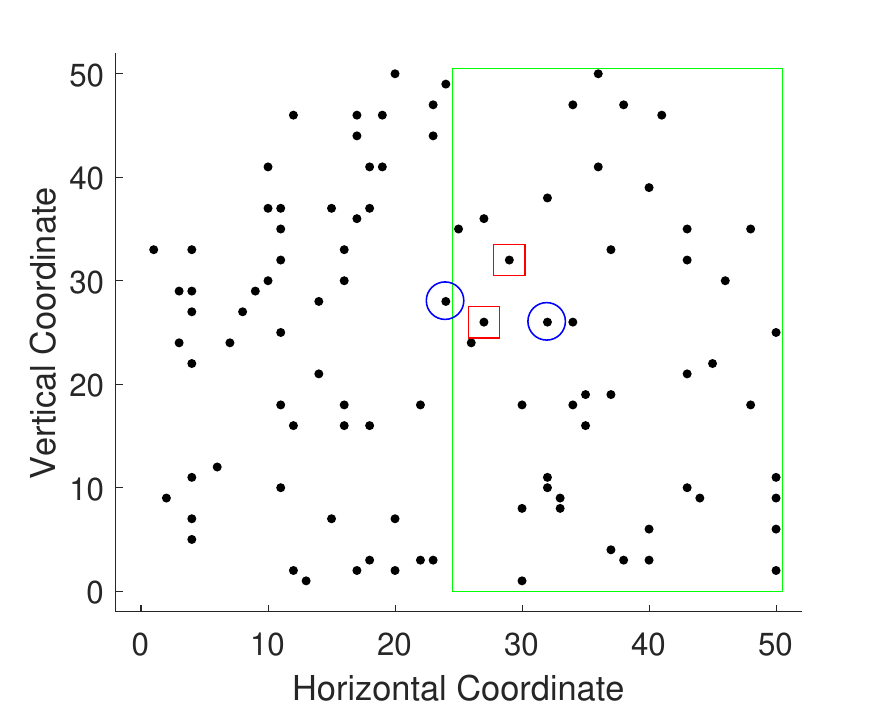}
}
\subfigure[{$\alpha^{\textrm{Left}} = 0, \ \alpha^{\textrm{Right}}=0.02$}]{
\includegraphics[width=0.45\linewidth]{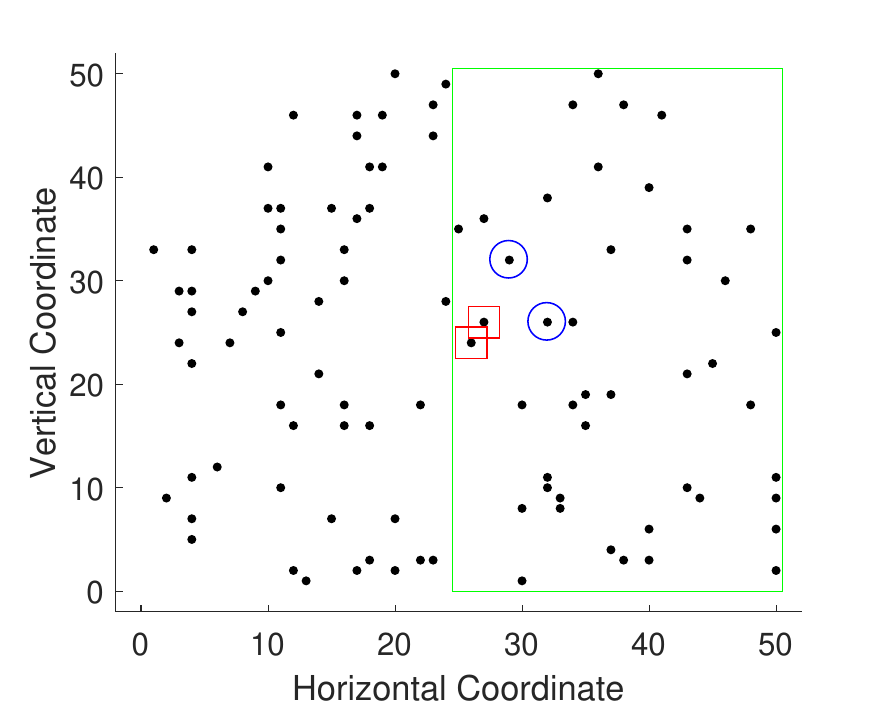}
}
\subfigure[{$\alpha^{\textrm{Left}} = 0, \ \alpha^{\textrm{Right}}=0.03$}]{
\includegraphics[width=0.45\linewidth]{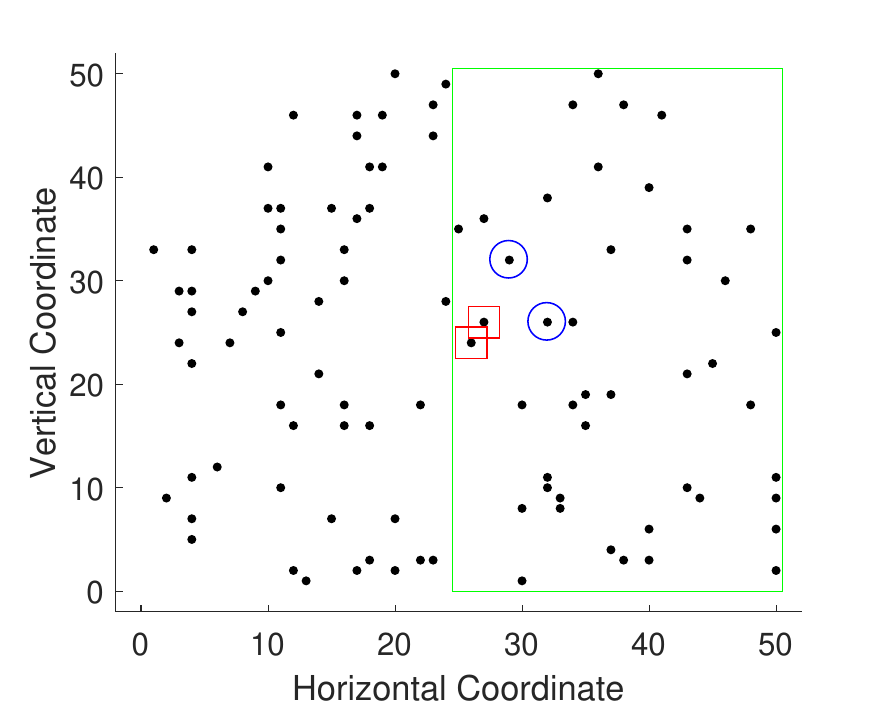}}
\subfigure{
\includegraphics[width=0.45\linewidth]{figure-legend.pdf}
}
\caption{{Leader's and follower's optimal locations given $\beta=0.05$ and heterogeneous $\alpha_j$ for left/right locations.} 
}
\label{figure 5}
\end{figure}

\begin{figure}[ht!]
\centering
\subfigure[{$\alpha^{\textrm{Left}} = 0, \ \alpha^{\textrm{Right}}=0$}]{
\includegraphics[width=0.45\linewidth]{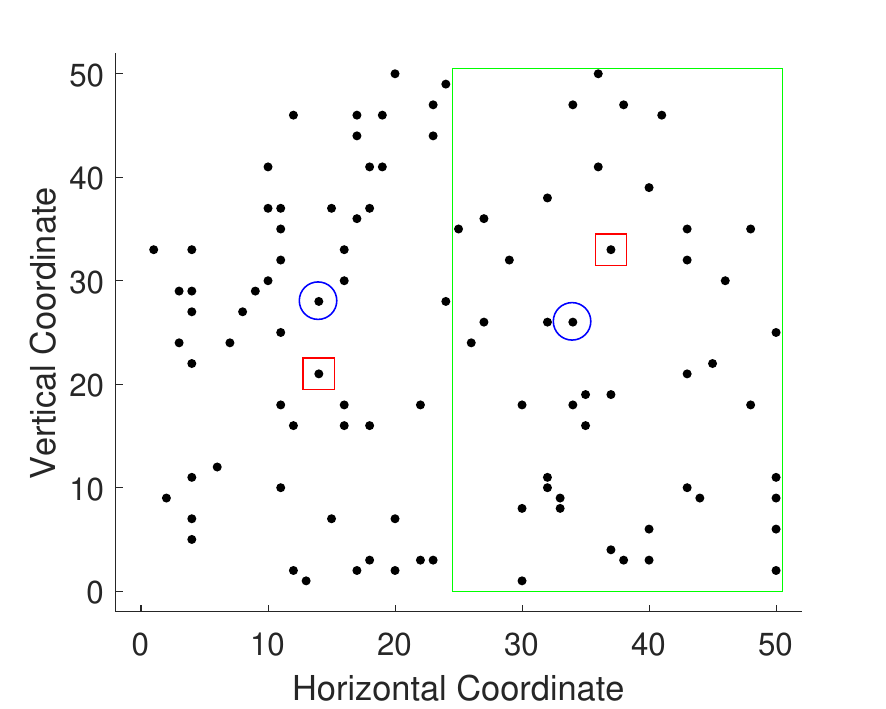}
}
\subfigure[{$\alpha^{\textrm{Left}} = 0, \ \alpha^{\textrm{Right}}=0.06$}]{
\includegraphics[width=0.45\linewidth]{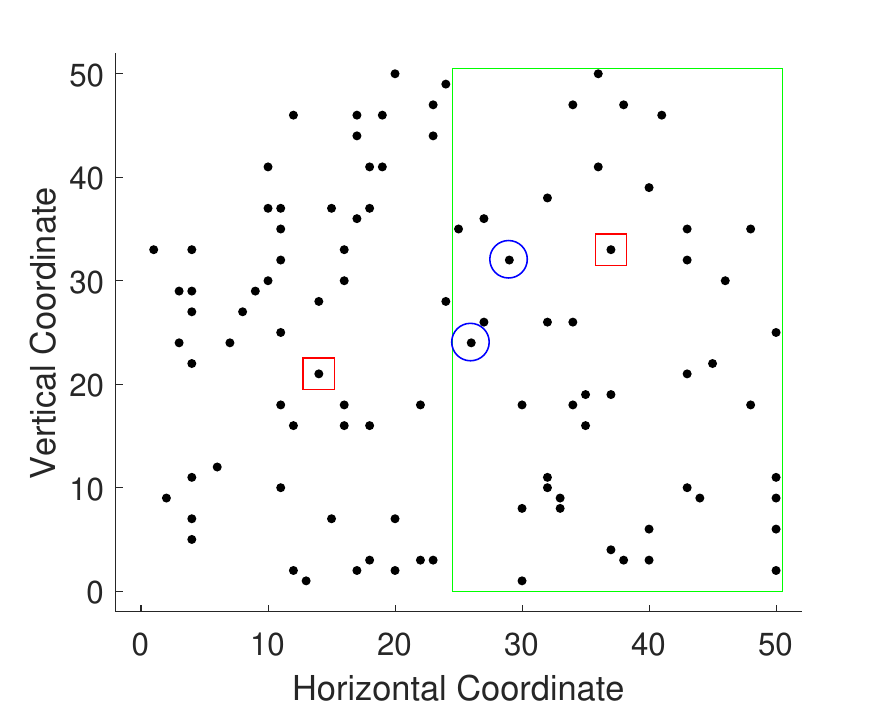}
}
\subfigure[{$\alpha^{\textrm{Left}} = 0, \ \alpha^{\textrm{Right}}=0.08$}]{
\includegraphics[width=0.45\linewidth]{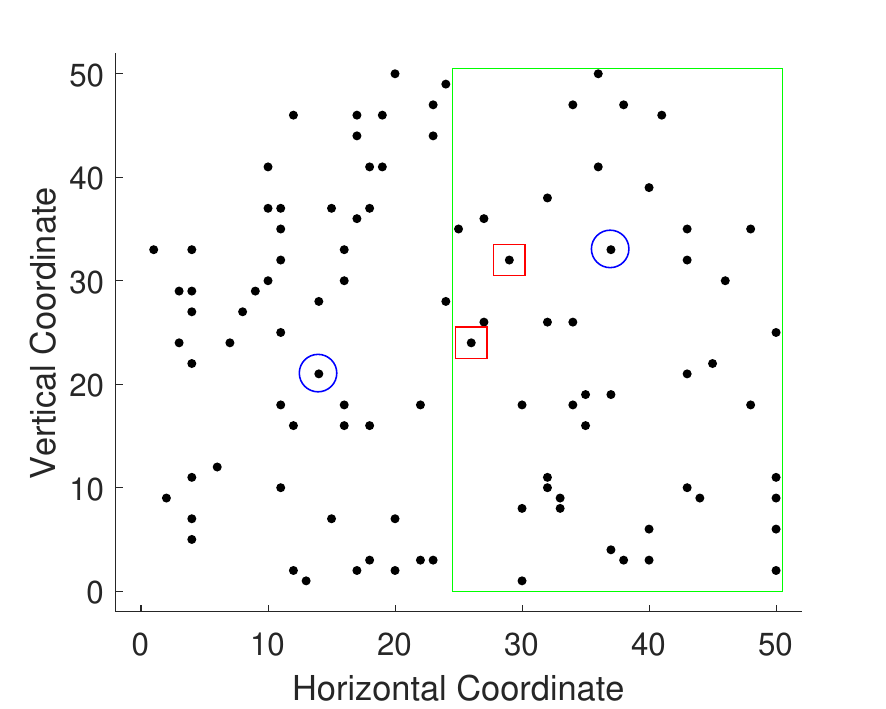}
}
\subfigure[{$\alpha^{\textrm{Left}} = 0, \ \alpha^{\textrm{Right}}=0.10$}]{
\includegraphics[width=0.45\linewidth]{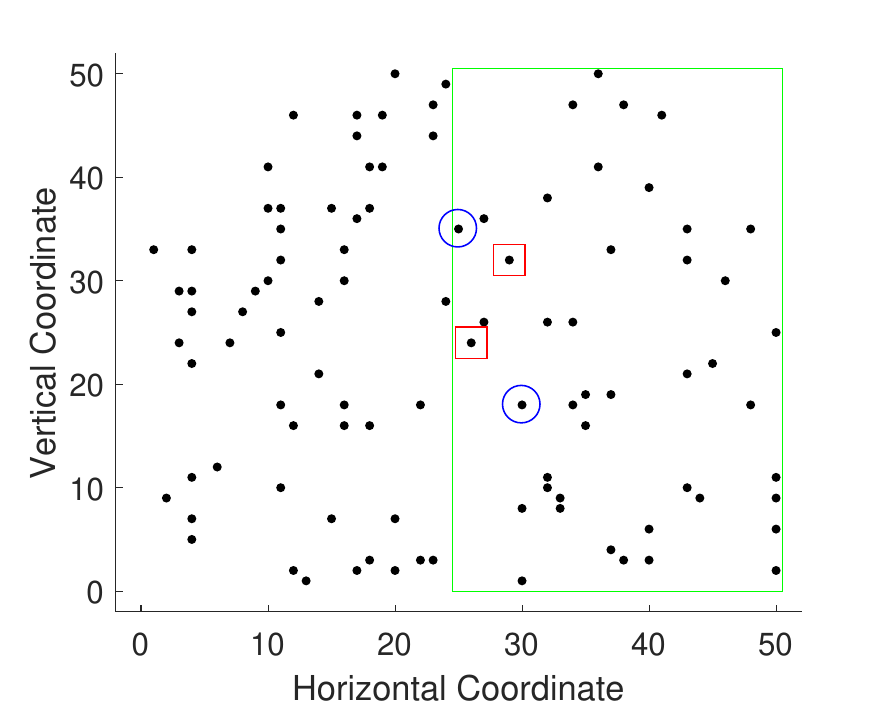}}
\subfigure{\includegraphics[width=0.45\linewidth]{figure-legend.pdf}
}
\caption{Leader's and follower's optimal locations given $\beta=0.1$ and heterogeneous $\alpha_j$ for left/right locations.}
\label{figure 6}
\end{figure}

Specifically, we differentiate the $\alpha_j$-values for candidate sites located on the left and right regions by keeping $\alpha_j = 0$ for all locations on the left and enlarging $\alpha_j$ for all the locations on the right (circled by a green rectangle in Figures \ref{figure 5} and \ref{figure 6}). In the two figures, we denote $\alpha_j$-values for all the left- and right-hand-side candidate locations as $\alpha^{\textrm{Left}}$ and $\alpha^{\textrm{Right}}$, respectively, and provide their specific values in each sub-figure's caption. 

In Figure \ref{figure 5}, both the leader and the follower locate their facilities around the center area in order to attract customers from all directions (recall that $\beta = 0.05$ indicates a low spatial impedance effect), but as $\alpha^{\textrm{Right}}$ increases, they both move facilities slightly towards the right into the green rectangle, reflecting the higher attractiveness level therein. Perhaps more interestingly, in Figure \ref{figure 6} where $\beta =0.1$, the customers intend to shop more locally and as a result, when $\alpha^{\textrm{Right}}$ increases from 0 to 0.06, we observe that the follower moves its facilities into the green rectangle while the leader's optimal locations remain unchanged (see Figure~\ref{figure 6}(b)). As we continue increasing $\alpha^{\textrm{Right}}$ to $0.08$, the leader moves its facilities to the right while the follower moves one of its facilities back to the left to attract local customers there, which are not covered by the leader anymore (see Figure~\ref{figure 6}(c)). Finally, when $\alpha^{\textrm{Right}}$ increases to $0.1$ in Figure~\ref{figure 6}(d), both the leader and the follower locate within the green rectangle due to the higher attractiveness level. Nevertheless, their facilities spread out more than those in Figure \ref{figure 5}(d), in order to attract more customers locally.

\section{Conclusion and Future Research}
\label{sec:concl}

This paper provides exact and approximation algorithms for solving sequential CFLP with probabilistic customer choice. We adopt integer programming and reformulation techniques to develop an exact, branch-and-cut algorithm. In addition, we derive an approximation algorithm with a constant guarantee on the ensuing market share. Extensive computational studies generate insights on winning market share and demonstrate the effectiveness of our approaches. 
For future research, the baseline (S-CFLP) model~\eqref{bl} can be further strengthened by incorporating additional designing features of the facilities, such as their capacities and the ``spillover'' effect (cf.~\citet{dan2019competitive}), which refers to that a facility spills its excessive demand to the nearby facilities. In addition, it is interesting to explore the use of other utility functions and customer choice models than the MNL model considered in this paper. An example is to explore how the customer demand of a facility may be impacted by the attractiveness levels of its nearby facilities. Another interesting direction is to consider uncertainty and asymmetric information in decision making. For example, the follower's strategy and budget may be unveiled to the leader. For all the aforementioned extensions, when they can still be formulated as a bilevel program (but with mixed-integer nonlinear structures in both levels), it is of high interest to investigate ways of designing algorithms based on cutting planes and other integer-programming related approaches to obtain high-quality solutions more efficiently. 

% Acknowledgments here
\ACKNOWLEDGMENT{Dr. Mingyao Qi is partially supported by the National Natural Science Foundation of China under grant No.\ 71772100 to work on this project. All the authors thank the Department Editor, Associate Editor, and three reviewers for making helpful suggestions and comments to improve the work.}

% References here (outcomment the appropriate case)
% CASE 1: BiBTeX used to constantly update the references
%   (while the paper is being written).
%\bibliographystyle{informs2014} % outcomment this and next line in Case 1
%\bibliography{<your bib file(s)>} % if more than one, comma separated
% CASE 2: BiBTeX used to generate mypaper.bbl (to be further fine tuned)

 % outcomment this line in Case 2

%If you don't use BiBTex, you can manually itemize references as shown below.

%\newpage
%\bibliographystyle{informs2014} % outcomment this and next line in Case 1
%\bibliography{cflp} % if more than one, comma separated

\clearpage
\newpage
\setcounter{page}{1}
%%%%%
%%% APPENDIX
%%%

\begin{center}
\Large{\bf Online Appendices of the Paper ``Sequential Competitive Facility Location: Exact and Approximate Algorithms''}~\\
~\\

\large{\bf Mingyao Qi, Ruiwei Jiang, and Siqian Shen}
\end{center}

\begin{appendices}
\section{(S-CFLP) Model Extensions}
\label{sec:extension}
This section considers extensions of the baseline (S-CFLP) model in~\eqref{bl} to take into account additional features and constraints. Section~\ref{sec:extension:setup} models heterogeneous costs for setting up facilities, 
Section~\ref{sec:extension:attractiveness} co-optimizes facility location and the choice of attractiveness levels, 
Section~\ref{sec:extension:outside} incorporates outside firms to share the market with the two competitors, and Section~\ref{sec:extension:change} considers potential changes in utility of the pre-existing facilities once new facilities are set up.

\subsection{Heterogeneous Setup Costs}\label{sec:extension:setup}
(S-CFLP) restricts the numbers of new facilities the leader and the follower can deploy through cardinality constraints~\eqref{bl-b} and~\eqref{bl-e}, respectively. Implicitly, these constraints assume that each facility incurs a homogeneous cost to set up. To model heterogeneous setup costs, we replace~\eqref{bl-b} and~\eqref{bl-e} with general 0-1 knapsack constraints 
\begin{align*}
\sum_{j \in J} {c^{\textrm{L}}_j x_j} \leq C^{\textrm{L}}, \quad \sum_{j \in J} {c^{\textrm{F}}_j y_j} \leq C^\textrm{F},
\end{align*}
where $c^{\textrm{L}}_j$ and $c^{\textrm{F}}_j$ are leader's and follower's costs of opening a facility at location $j$, and $C^{\textrm{L}}$ and $C^{\textrm{F}}$ are their total budgets, respectively. It can be observed that the solution method and valid inequalities described in Section~\ref{sec:cuts} remain applicable in this extension.

\subsection{Attractiveness Level}\label{sec:extension:attractiveness}
(S-CFLP) assumes that the attractiveness level $\alpha_j$ of each facility is fixed and known, making it impossible for a competitor to strategically increase its market share by adjusting, e.g., the price level at a facility. To co-optimize locations and attractiveness levels, we replicate each site $j$, which now consists of a set $N_j$ of potential facilities to build and all these facilities share the same distances $d_{ij}$ to demand nodes. These (replicated) facilities only differ in the attractiveness level, denoted by $\alpha_{jn}$, and setup costs, denoted by $c^{\textrm{L}}_{jn}$ (for the leader) and $c^{\textrm{F}}_{jn}$ (for the follower), for all $n \in N_j$. Accordingly, we extend the decision variables $(x_j, y_j)$ to be $(x_{jn}, y_{jn})$ to reflect the attractiveness level choice. For example, $x_{jn} = 1$ if and only if the leader deploys a facility at site $j$ with attractiveness level $\alpha_{jn}$. Additionally, we denote $w_{ijn} := \exp\{\alpha_{jn} - \beta d_{ij}\}$. Then, the leader's market share
\begin{equation*}
L(x, y) = \sum_{i\in I} h_i \left(\frac{U^{\tinyl}_i+\sum_{j\in J}\sum_{n \in N_j}w_{ijn} x_{jn}}{U^{\tinyl}_i+U^{\tinyf}_i+\sum_{j\in J}\sum_{n \in N_j}w_{ijn}(x_{jn}\vee y_{jn})}\right)
\end{equation*}
follows from the MNL model. This leads to the following extended model:
\begin{align*}
(\text{\bf S-CFLP-$\alpha$}) \quad \max_x \min_y \ & \ L(x, y) \\
\mbox{s.t.} \ & \ \sum_{n \in N_j} x_{jn} \leq 1, \quad \sum_{n \in N_j} y_{jn} \leq 1, \quad \forall j \in J, \\
& \ \sum_{j \in J}\sum_{n \in N_j} {c^{\textrm{L}}_{jn} x_{jn}} \leq C^{\textrm{L}}, \quad \sum_{j \in J}\sum_{n \in N_j} {c^{\textrm{F}}_{jn} y_{jn}} \leq C^\textrm{F}, \\
& \ x_{jn}, y_{jn} \in \{0,1\} \quad \forall j \in J, \ \forall n \in N_j.
\end{align*}
The first set of constraints ensure that the leader (respectively, the follower) commits to at most one attractiveness level if she builds a facility at a site. We notice that (S-CFLP-$\alpha$) admits the same RO formulation as (S-CFLP) and, as a result, (S-CFLP-$\alpha$) can be solved by a similar branch-and-cut framework as Algorithm~\ref{algo-b&c}.

\subsection{Outside Competitors}\label{sec:extension:outside}
In reality, a customer may choose to patronize options outside of the two competitors. An example is the Amazon home delivery in the retail business. When outside options are taken into account, there are two ways to apply (S-CFLP). First, we can estimate the market share occupied by the outside options and consider (S-CFLP) within the remaining share for the two competitors only. Although this allows us to directly apply (S-CFLP), it overlooks the impact of the outside options on how customers choose among the players (i.e., leader, follower, and outside competitors).

The second way addresses this issue by incorporating the utility of the outside options into the MNL model. Specifically, suppose that each customer from demand location $i$ receives a utility $U^{\tinyo}_i$ by patronizing the outside options. Then, it follows from the MNL model that the leader's market share now becomes
\begin{equation*}
L^+(x, y) = \sum_{i\in I} h_i \left(\frac{U^{\tinyl}_i + \sum_{j \in J}w_{ij}x_j}{U^{\tinyl}_i+U^{\tinyf}_i+\sum_{j \in J}w_{ij}(x_{j} + y_{j}) + U^{\tinyo}_i}\right)
\end{equation*}
and the follower's market share equals
\begin{equation*}
F^+(x, y) = \sum_{i\in I} h_i \left(\frac{U^{\tinyf}_i + \sum_{j \in J}w_{ij}y_j}{U^{\tinyl}_i+U^{\tinyf}_i+\sum_{j \in J}w_{ij}(x_{j} + y_{j}) + U^{\tinyo}_i}\right).
\end{equation*}
Then, (S-CFLP) can be extended as follows to consider outside options:
\begin{align*}
(\mbox{\bf S-CFLP-O}) \qquad z^{\tinyo} := \max_{x \in \mathcal{X}} \ & \ L^{+}(x, y^*), \quad \mbox{where} \ y^* \in \argmax_{y \in \mathcal{Y}(x)} \ F^+(x, y).
\end{align*}
Unfortunately, since $L^+(x, y) + F^+(x, y) < 1$, (S-CFLP-O) no longer admits the RO reformulation as in~\eqref{eq:max-min}. Nevertheless, one can show that
\begin{equation*}
\max_{x \in \mathcal{X}} \min_{y \in \mathcal{Y}(x)} L^+(x, y) \leq z^{\tinyo} \leq \max_{x \in \mathcal{X}} \min_{y \in \mathcal{Y}(x)} \big\{1-F^+(x, y)\big\}.
\end{equation*}
That is, the bilevel (S-CFLP-O) model can be approximated by two RO formulations, from both above and below respectively (see a proof in Appendix~\ref{apx-ext-outside}). Furthermore, like in Theorem~\ref{lemma:l1}, the decision dependency in these RO formulations can be relaxed (i.e., replacing $\mathcal{Y}(x)$ with $\mathcal{Y}$) without loss of optimality (see Appendix~\ref{apx-ext-outside}). As a result, these RO formulations can be solved by a similar branch-and-cut framework as Algorithm~\ref{algo-b&c}.

\subsection{Utility Change}\label{sec:extension:change}
The opening of new facilities may change the utility of a customer patronizing a pre-existing facility. An example is the economy of agglomeration, which suggests that each of a cluster of facilities can attract more customers than a stand-alone facility. To capture the utility change, we assume that opening new facilities changes $U^{\tinyl}_i$, the utility received by the customers located in node $i$ patronizing the leader's pre-existing facilities, to be
\begin{equation*}
U^{\tinyl}_i + \sum_{j \in J} v^{\tinyl}_{ij}(x_j \vee y_j),
\end{equation*}
where parameter $v^{\tinyl}_{ij}$ evaluates the sensitivity of $U^{\tinyl}_i$ on the opening of a facility in location $j$ by either the leader or the follower. Similarly, we assume that opening new facilities changes $U^{\tinyf}_i$, the utility of patronizing from the follower's pre-existing facilities, to be $U^{\tinyf}_i + \sum_{j \in J} v^{\tinyf}_{ij}(x_j \vee y_j)$. Accordingly, the leader's market share follows from the MNL model:
\begin{equation*}
L^{\textrm{U}}(x, y) = \sum_{i\in I} h_i \left(\frac{U^{\tinyl}_i + \sum_{j \in J} v^{\tinyl}_{ij}(x_j \vee y_j) + \sum_{j \in J}w_{ij}x_j}{U^{\tinyl}_i+U^{\tinyf}_i+\sum_{j \in J}(v^{\tinyl}_{ij} + v^{\tinyf}_{ij} + w_{ij})(x_{j} \vee y_{j})}\right).
\end{equation*}
This extends (S-CFLP) to 
\begin{equation*}
\text{\bf (S-CFLP-U)} \quad \max_{x \in \mathcal{X}} \min_{y \in \mathcal{Y}(x)} L^{\textrm{U}}(x, y).
\end{equation*}
It can be shown that the decision-dependency $\mathcal{Y}(x)$ can once again be relaxed, i.e., we can replace $\mathcal{Y}(x)$ with $\mathcal{Y}$ without loss of optimality as in Theorem~\ref{lemma:l1} (see a proof in Appendix~\ref{apx-prop:extension:change}). In addition, the objective function $L^{\textrm{U}}(x, y)$ in (S-CFLP-U) admits a second-order conic representation (see Appendix~\ref{apx-prop:extension:change}). As a result, (S-CFLP-U) can be solved by a similar branch-and-cut framework as Algorithm~\ref{algo-b&c}.

\section{A Detailed Literature Review}
\label{sec:lit}

In this section, we review the most relevant literature on static CFLP and sequential CFLP. 

\paragraph{Static CFLP:} The problem was first proposed by \citet{slater1975underdevelopment} and also known as the $(r\mid X_p)$-medianoid problem \citep{hakimi1983}, in which a decision maker locates $r$ facilities when the $p$ locations of its competitors', denoted by $X_p$, are given. \citet{Plastria2001} provided a survey of the models and solution approaches, including heuristic methods, for static CFLP. Assuming a deterministic choice model and that both customers and facilities are located at discrete points of a network, one can model static CFLP as a MILP \citep[see, e.g.,][]{revelle1986, revelle1995}. In contrast, \citet{luce1959} and \citet{huff1964} proposed a probabilistic choice model, which characterized customer behavior given utilities of multiple facilities. \citet{BenatiHansen2002} adopted this probabilistic choice model in static CFLP and formulated it as a MINLP. They exploited both concavity and submodularity of the objective function to construct an exact and a heuristic algorithm, respectively. In the exact algorithm, they recast this model as a MILP, which was further strengthened in \citet{haase2009discrete}, \citet{aros2013p}, and \citet{zhang2012impact} to improve the computational performance. Recently, \citet{Haase2014} compared the aforementioned methods via extensive computational studies and empirical analysis, and \citet{Freire2016} conducted numerical studies on diverse instances, showing that the state-of-the-art methods can handle small- to medium-sized problem instances. \citet{Ljubic2018} proposed a branch-and-cut algorithm using outer approximation (OA) inequalities and submodularity inequalities for solving static CFLP with random utilities. Their method outperformed the state-of-the-art approaches with two to three orders of magnitude. \citet{mai2020multicut} proposed multicut OA inequalities for groups of demand points and implemented them as cutting planes. Different from most existing static CFLP models, \citet{dan2019competitive} considered random utility based on both travel time and queuing delay at facilities, leading to a MINLP formulation. They derived a piecewise linear approximation and a heuristic method for solving this model. In this paper, we consider sequential CFLP, which gives rise to a more challenging  bilevel program with MINLPs at both levels. Additionally, in terms of methodology, our valid inequalities are different from those derived in \citet{Ljubic2018, mai2020multicut}. Specifically, \citet{Ljubic2018} derived OA inequalities from the concave objective function of the static CFLP model, but our sequential CFLP model undermines such concavity. Nevertheless, we are able to ``bulge up'' our objective function to restore concavity while retaining exactness. This yields a new class of valid inequalities that have not been developed in the existing literature and can significantly speed up the computation of the MINLP reformulation of the bilevel sequential CFLP derived in Section \ref{sec:model}. 

\paragraph{Sequential CFLP:} The problem is also known as \emph{CFLP with leader-follower game} in the literature and can be modeled as a bilevel program that involves the leader's location model in the upper level and a static CFLP in the lower level (i.e., the follower's location model). Assuming deterministic choice and discrete location space, one can reformulate the bilevel program as a single-level MILP \citep{Plastria2008,Roboredo2013,Alekseeva2015,Gentile2018}, which involves a polynomial number of variables but an exponential number of constraints. Instances with 100 candidate facility sites and 100 customers can be optimally solved through a branch-and-cut algorithm using commercial solvers \citep[see][]{Gentile2018}. \citet{Drezner1998} were the first to study sequential CFLP with probabilistic choice, for which they developed a heuristic algorithm (without optimality-gap guarantees). \citet{saiz2009branch} extended this work and applied a branch-and-bound algorithm to seek exact solutions, but assumed a planar (continuous) location space. Different from \citet{Drezner1998} and \citet{saiz2009branch}, \citet{Kucukaydn2011} optimized the leader's facility locations with \emph{fixed} locations from the follower who only optimizes facility attractiveness, e.g., facility sizes. In their model, the lower-level problem became a convex NLP. The authors then applied the Karush-Kuhn-Tucker (KKT) optimality conditions to seek the follower's optimal decisions, and the bilevel program became a single-level MINLP after adding the KKT conditions in the upper level. To the best of our knowledge, \citet{Kucukaydn2012} is the only study that considered the same sequential CFLP as ours, which assumes a probabilistic choice model and discrete location space. \citet{Kucukaydn2012} developed three heuristics, as well as an $\epsilon$-optimal method by iteratively fixing the leader's decisions. Small-sized instances with only 16 candidate facilities were solved using these heuristic approaches without optimality guarantee.

\section{Proofs} 
\label{apx-proofs}

\subsection{Proof of Theorem \ref{thm:hardness}} \label{apx-thm-hardness}
\proof{Proof:}
By Theorem \ref{lemma:l1}, which we shall prove later, (S-CFLP) is equivalent to $\displaystyle \max_{x \in \mathcal{X}} \min_{y \in \mathcal{Y}} L(x, y)$, where $\mathcal{X} = \{x \in \{0, 1\}^{|J|}: e^{\top}x \leq p\}$, $\mathcal{Y} = \{y \in \{0, 1\}^{|J|}: e^{\top}y \leq r\}$, and $L(x, y)$ is defined in \eqref{eq:l}. In addition, for any fixed $y \in \mathcal{Y}$, $L(x, y)$ is submodular in the index set $X$ of variables $x$ (defined as $X:=\{j \in J: x_j =1\}$) by Proposition \ref{p2} presented later. Hence, (S-CFLP) is equivalent to a robust submodular maximization model as defined in formulation (2) of~\cite{krause2008robust}. The conclusion follows from Theorem 3 of~\cite{krause2008robust}.
\Halmos
\endproof

\subsection{Proof of Theorem \ref{lemma:l1}} \label{apx-thm-submodularity}
\proof{Proof:}
Define $F_{0}:=\mathcal{Y}$ and $F_{x}:=\{y\in\{0,1\}^{|J|}: y_{j}\leq 1-x_{j},\ \forall j\in J\}$. Then, for any $x\in\mathcal{X}$, we have $\displaystyle\theta^{+}(x)=\min_{y\in F_{0}\cap F_{x}}L^{+}(x,y)$ and $ \theta(x)=\min\{\theta_{1}(x),\theta_{2}(x)\}$, where $\displaystyle\theta_{1}(x)=\min_{y\in F_{0}\cap F_{x}}L(x,y)$ and $ \displaystyle\theta_{2}(x)=\min_{y\in F_{0}\backslash F_{x}}L(x,y)$. In addition, for any $y\in F_{0}\cap F_{x}$, we have $L^{+}(x,y)=L(x,y)$ because $x_{j}+y_{j}=x_{j}\vee y_{j}$ for all $j \in J$. Hence, $\theta^{+}(x)=\theta_{1}(x)$. It remains to show that $\theta_{1}(x)\leq\theta_{2}(x)$ for all $x \in \mathcal{X}$ and then the equivalence can be drawn between $\theta^+(x)$ and $\theta(x)$.

To this end, for any $x \in \mathcal{X}$ and $\overline{y} \in F_{0}\backslash F_{x}$, we construct a $\tilde{y} \in F_0\cap F_x$ such that $L(x, \tilde{y}) \leq L(x, \overline{y})$. Since $\overline{y} \in F_{0}\backslash F_{x}$, there exists a nonempty subset $K\subseteq J$ such that (i) $\overline{y}_{j} > 1-x_{j}$, i.e., $x_{j}=\overline{y}_{j}=1$ for all $j\in K$ and (ii) $\overline{y}_{j} \leq 1-x_{j}$ for all $j\in J\backslash K$. We claim that there exists a subset $M \subseteq J\backslash K$ with $|M|=|K|$ and $x_{j}=\overline{y}_{j}=0$ for all $j \in M$. To see this, we denote $J_{mn} := \{j \in J: x_j = m, \overline{y}_j = n\}$ for $m, n \in \{0, 1\}$. Then, it holds that
\begin{align*}
    |J_{00}| & = |J| - |J_{11}| - |J_{01}| - |J_{10}| \\
    & \geq |J| - |K| - (r - |K|) - (p - |K|) \\
    & = |K| + (|J| - r - p) \ \geq \ |K|,
\end{align*}
where the first inequality is because $\sum_{j \in J} \overline{y}_j \leq r$ and $\sum_{j \in J} x_j \leq p$ and the last inequality holds because $p+r\leq|J|$. Then, the existence of $M$ follows from the pigeonhole principle. Now define a $\tilde{y} \in \{0, 1\}^{|J|}$ such that $\tilde{y}_{j}=\overline{y}_{j}$ for all $j\in J\backslash (K\cup M)$, $\tilde{y}_{j}=0$ for all $j\in K$, and $\tilde{y}_{j}=1$ for all $j\in M$. Then, $\tilde{y}\in F_{0}\cap F_{x}$ by construction. In addition, for each $i\in I$, we have
\begin{eqnarray*}
& & \sum_{j\in J}w_{ij}(x_{j}\vee\overline{y}_{j})   \\
&= & \sum_{j\in K}w_{ij}(x_{j}\vee\overline{y}_{j})+\sum_{j\in M}w_{ij}(x_{j}\vee\overline{y}_{j})+\sum_{j\in J\backslash (K\cup M)}w_{ij}(x_{j}\vee\overline{y}_{j})  \\
&= &  \sum_{j\in K}w_{ij}\times 1+\sum_{j\in M}w_{ij}\times 0+\sum_{j\in J\backslash (K\cup M)}w_{ij}(x_{j}\vee\overline{y}_{j}) \\
&\leq &  \sum_{j\in K}w_{ij}\times 1+\sum_{j\in M}w_{ij}\times 1+\sum_{j\in J\backslash (K\cup M)}w_{ij}(x_{j}\vee\overline{y}_{j}) \\
&= & \sum_{j\in K}w_{ij}(x_{j}\vee\tilde{y}_{j})+\sum_{j\in M}w_{ij}(x_{j}\vee\tilde{y}_{j})+\sum_{j\in J\backslash (K\cup M)}w_{ij}(x_{j}\vee\tilde{y}_{j})  \\
&= & \sum_{j\in J}w_{ij}(x_{j}\vee\tilde{y}_{j}).
\end{eqnarray*}
As a result, $L(x,\overline{y}) \geq L(x,\tilde y)$ and the proof is completed.
\Halmos
\endproof

\subsection{Proof of Proposition \ref{p2}} \label{apx-p2}
\proof{Proof:}
For all $i\in I$, $X, Y\subseteq J$, and $k\in J \backslash X$, it follows from \eqref{eq:f-i-J'} that
\begin{align}
& f_{i,Y}(X\cup\{k\})-f_{i,Y}(X) \nonumber \\
= & \frac{U^{\tinyl}_i+\sum_{j\in X\cup\{k\}}w_{ij}}{U^{\tinyl}_i+U^{\tinyf}_i+\sum_{j\in X\cup Y\cup\{k\}}w_{ij}}-\frac{U^{\tinyl}_i+\sum_{j\in X}w_{ij}}{U^{\tinyl}_i+U^{\tinyf}_i+\sum_{j\in X\cup Y}w_{ij}} \nonumber \\[0.1cm]
= & \begin{cases} \frac{w_{ik}}{U^{\tinyl}_i+U^{\tinyf}_i+\Sigma_{j\in X \cup Y}w_{ij}} & \mbox{if } k\in Y \\ 
\frac{w_{ik}\left(U^{\tinyf}_i+\Sigma_{j\in Y\backslash X}w_{ij}\right)}{(U^{\tinyl}_i+U^{\tinyf}_i+\Sigma_{j\in X\cup Y}w_{ij}+w_{ik})(U^{\tinyl}_i+U^{\tinyf}_i+\Sigma_{j\in X\cup Y}w_{ij})} & \mbox{if } k\not\in Y. \end{cases}
\end{align}
It follows that $f_{i,Y}(X\cup\{k\})-f_{i,Y}(X)$ is non-increasing in $X$, i.e., $ f_{i,Y}(X\cup \{k\})-f_{i,Y}(X)\geq f_{i,Y}(X'\cup\{k\})-f_{i,Y}(X')$ for all $X\subseteq X'\subseteq J\backslash \{k\}$. Therefore, $f_{i,Y}$ is submodular and so is $L_{Y}$ because $L_{Y}$ is a linear combination of $f_{i,Y}$. This completes the proof. 
\Halmos
\endproof

\subsection{Proof of Proposition \ref{p4}} \label{apx-p4}
\proof{Proof:}
For notational brevity, we omit the subscript $Y$ in this proof. Pick any subsets $S \subseteq R \subseteq J$ and any element $j \in J \setminus R$. By definition, we have 
\begin{align*}
H(S \cup \{j\}) = 
L(S \cup \{j\}) - \sum_{k \in S} \rho(J \setminus \{k\};k) (1 - \hat{x}_k) - \rho(J\setminus\{j\}; j)(1 - \hat{x}_j) + \sum_{k \in J \setminus S} \rho(S;k)\hat{x}_k -\rho(S;j)\hat{x}_j.
\end{align*}
Then,
\begin{align*}
H(S \cup \{j\}) - H(S) & = 
L(S \cup \{j\}) - L(S) - \rho(J\setminus\{j\}; j) (1 - \hat{x}_j) - \rho(S;j)\hat{x}_j \\
& = \bigl(\rho(S;j) - \rho(J\setminus\{j\}; j)\bigr) (1 - \hat{x}_j).
\end{align*}
It follows that
\begin{align*}
\bigl[H(S \cup \{j\}) - H(S)\bigr] - \bigl[H(R \cup \{j\}) - H(R)\bigr] & = \bigl(\rho(S;j) - \rho(R;j)\bigr) (1 - \hat{x}_j) \geq 0,
\end{align*}
where the inequality follows from Proposition \ref{p2}. Therefore, $H(S \cup \{j\}) - H(S)$ is non-increasing in $S$ and thus, $H$ is submodular by definition. This completes the proof. 
\Halmos
\endproof

\subsection{Proof of Proposition \ref{lemma:l2}} \label{apx-lemma-l2}
\proof{Proof:} 
For all $x,y \in \{0, 1\}^{|J|}$, it is easy to verify $x_j \vee y_j = (1- y_j)x_j+y_j$ and $x_j = - y_j x_j^2+(1+ y_j)x_j$. It follows that $\widehat L (x,y) = L(x, y)$ for all $x \in \{0, 1\}^{|J|}$.

It remains to show the concavity of $\widehat L(x,y)$. Since the sum of concave functions is concave, it suffices to show that $\widehat{L}_i (x) := \frac{U^{\tinyl}_i+\sum_{j \in J} {w_{ij} (-y_j x_j^2 +(1+ y_j) x_j)}}{U^{\tinyl}_i+U^{\tinyf}_i+\sum_{j \in J} {w_{ij}\bigl[(1-y_j)x_j+y_j\bigr]}}$ is concave for all $i \in I$. For ease of exposition, we denote its numerator $Q := U^{\tinyl}_i+\sum_{j \in J} {w_{ij} \bigl[-y_j x_j^2 +(1+y_j) x_j\bigr]}$, its denominator $P := U^{\tinyl}_i+U^{\tinyf}_i+\sum_{j \in J} {w_{ij}\bigl[(1-y_j)x_j+y_j\bigr]}$, and its Hessian $H := [h_{kl}]$. It follows that
\begin{eqnarray}
\frac{\partial \widehat{L}_i(x)}{\partial x_k} & \ = \ & \frac{-w_{ik} (1-y_k)Q}{P^2} + \frac{w_{ik}(-2y_k x_k+1+y_k)}{P}, \quad \forall k \in J, \nonumber \\[0.1cm]
h_{kk} & \ = \ & \frac{\partial^2 \widehat{L}_i(x)}{\partial x_k^2 }   = \frac{2w_{ik}^2 (1-y_k)^2 Q}{P^3} + \frac{-w_{ik}^2(1-y_k)(-2 y_k x_k+1+ y_k)}{P^2} \nonumber \\
& &\hspace{20mm} + \frac{-w_{ik}^2(1-y_k)(-2 y_k x_k+1+ y_k)}{P^2}+\frac{-2w_{ik} y_k}{P} \nonumber \\[0.1cm]
& \ = \ & \frac{2w_{ik}^2 (1-y_k)^2 Q}{P^3} - \frac{2w_{ik}^2(1-y_k)(-2 y_k x_k+1+ y_k)}{P^2}+\frac{-2w_{ik} y_k}{P}, \quad \forall k \in J, \nonumber \\[0.1cm]
\mathrm{and} \quad h_{k\ell}& \ = \ & \frac{\partial^2 \widehat{L}_i(x)}{\partial x_k \partial x_\ell } = \frac{2w_{ik} w_{i\ell} (1-y_k)(1-y_\ell) Q}{P^3} + \frac{-w_{ik}(1-y_k)w_{i\ell}(-2 y_\ell x_\ell+1+ y_\ell)}{P^2} \nonumber \\
& & \hspace{20mm} + \frac{-w_{i\ell}(1-y_\ell)w_{ik}(-2 y_k x_k+1+ y_k)}{P^2} \nonumber \\
& \ = \ & \frac{w_{ik}w_{i\ell}}{P^3} \Bigl[ (1-y_k)\Bigl((1-y_\ell) Q - (-2 y_\ell x_\ell +1+ y_\ell)P\Bigr) \nonumber\\
& & + (1-y_\ell)\Bigl((1-y_k) Q - (-2 y_k x_k +1+ y_k)P\Bigr) \Bigr], \quad \forall k, \ell \in J: k \neq \ell. \nonumber
\end{eqnarray}
Denote $J_1 := \{j \in J: y_j=1\}$ and  $J_0 := \{j \in J: y_j=0\}$. We simplify the expression of $H$ by examining the following four cases:\\
Case 1. If $y_k \in J_1, y_\ell \in J_1$ then $\displaystyle h_{kk} =-\frac{2}{P}w_{ik}$ and $h_{kl}=0$.\\[0.2cm]
Case 2. If $y_k\in J_0, y_\ell\in J_0$ then $\displaystyle h_{kk} = -\frac{2w_{ik}^2}{P^3}(P-Q)$ and $\displaystyle h_{kl} = -\frac{2w_{ik}w_{il}}{P^3}(P-Q)$.\\[0.2cm]
Case 3. If $y_k\in J_1, y_\ell\in J_0$ then $\displaystyle h_{kk}=-\frac{2}{P}w_{ik}$ and $\displaystyle h_{kl}= -\frac{2w_{ik}w_{il}}{P^2}(1-x_k)$.\\[0.2cm]
Case 4. If $y_k\in J_0, y_\ell\in J_1$ then $\displaystyle h_{kk} =-\frac{2w_{ik}^2}{P^3}(P-Q)$ and $\displaystyle h_{kl} = -\frac{2w_{ik}w_{il}}{P^2}(1-x_\ell)$.\\
To show that $H$ is negative semidefinite, we prove that $v^{\top} H v\le 0$ for all non-zero $v$ in $\mathbb{R}^{|J|}$. Indeed,
%\small
\begin{eqnarray}
v^{\top} Hv & = &\sum_{k \in J} {\sum_{\ell \in J} {h_{k\ell} v_k v_\ell}} \nonumber\\
& = &-\frac{2}{P^3} \Bigl[(P-Q)\sum_{k \in J_0} \sum_{\ell \in J_0} w_{ik} w_{i\ell} v_k v_\ell 
+ \sum_{k \in J_1} \sum_{\ell \in J_1, \ell \ne k} 0 + \sum_{k \in J_1} \bigl(- {P^2} w_{ik} v_k^2\bigr) \nonumber\\
& & \hspace{10mm} + \sum_{k\in J_1} \sum_{\ell \in J_0} \Bigl(-P(1-x_k)w_{ik} w_{i\ell} v_k v_\ell \Bigr)
+\sum_{k\in J_0} \sum_{\ell \in J_1} \Bigl(-P(1-x_\ell)w_{ik} w_{i\ell} v_k v_\ell \Bigr) \Bigr] \nonumber \\
&=&-\frac{2}{P^3}\left[(P-Q) \Bigl(\sum_{k \in J_0}w_{ik} v_k\Bigr)^2 + \sum_{k \in J_1} \bigl(P \sqrt{ w_{ik}} v_k\bigr)^2 + 2\sum_{k\in J_1} \sum_{\ell \in J_0} P(1-x_k)w_{ik}w_{i\ell} v_k v_\ell \right]. \label{eq:nsd:1}
\end{eqnarray}
%\normalsize
We notice that
\begin{eqnarray}
P-Q & = & U^{\tinyf}_i+\sum_{j\in J} w_{ij}\left\{\bigl[(1-y_j)x_j +y_j\bigr] -\bigl[-y_j x_j^2 +(1+y_j)x_j\bigr]\right\} \nonumber \\
& = &U^{\tinyf}_i+\sum_{j \in J} w_{ij} y_j (1-x_j)^2 =U^{\tinyf}_i+\sum_{j\in J_1} w_{ij} (1-x_j)^2. \label{eq:nsd:2}
\end{eqnarray}
Plugging \eqref{eq:nsd:2} into \eqref{eq:nsd:1} yields
%\footnotesize
\begin{eqnarray*}
v^{\top} H v &= & -\frac{2}{P^3}\Biggl\{\Bigl(\sum_{k\in J_1} \sqrt{w_{ik}} (1-x_k)\Bigr)^2 \sum_{k \in J_0}(w_{ik} v_k)^2 +  \sum_{k \in J_1} \bigl(P \sqrt{ w_{ik}} v_k\bigr)^2 \\
& & + 2\Bigl(\sum_{k\in J_1} P(1-x_k)w_{ik}v_k\Bigr)\sum_{\ell \in J_0}  w_{i\ell}v_\ell + U^{\tinyf}_i \Bigl(\sum_{k \in J_0} w_{ik} v_k\Bigr)^2 \Biggr\} \\
&=& -\frac{2}{P^3} \sum_{k\in J_1}\Biggl( \Bigl(\sqrt{w_{ik}}(1-x_k) \sum_{l\in J_0} w_{il}v_\ell \Bigr)^2 + \bigl(P\sqrt{w_{ik}}v_k\bigr)^2 
+ 2P(1-x_k)w_{ik}v_k\sum_{l \in J_0} w_{il}v_\ell \Biggr) \\
& & -\frac{2U^{\tinyf}_i}{P^3} \bigl(\sum_{k \in J_0} w_{ik} v_k\bigr)^2 \\
&=& -\frac{2}{P^3} \sum_{k\in J_1} \Bigl(\sqrt{w_{ik}}(1-x_k) \sum_{l\in J_0} w_{il}v_\ell+P\sqrt{w_{ik}}v_k\Bigr)^2 -\frac{2U^{\tinyf}_i}{P^3} \Bigl(\sum_{k \in J_0} w_{ik} v_k\Bigr)^2 \ \leq 0.
\end{eqnarray*}
%\normalsize
This completes the proof.
\Halmos
\endproof

\subsection{Features of $\theta \leq \widehat{L}(x, y)$ and A Proof}
\label{apx:prop:soc}
\begin{proposition} \label{prop:soc}
For fixed $y \in \{0, 1\}^{|J|}$, define $J_0 := \{j \in J: y_j = 0\}$, $J_1 := \{j \in J: y_j = 1\}$, and $\displaystyle s_i(y) := U^{\tinyl}_i + U^{\tinyf}_i + \sum_{j \in J}w_{ij}y_j$ for all $i \in I$. Then, the inequality $\theta \leq \widehat{L}(x, y)$ holds valid if and only if there exist $\{\theta_i\}_{i \in I}$ such that $\displaystyle \theta \leq \sum_{i \in I}h_i \theta_i$ and 
\begin{equation}
    \left\|\begin{bmatrix}
    2\sqrt{w_{ij}}(1 - x_j)_{j \in J_1}\\
    \displaystyle\sum_{j \in J} w_{ij}(1 - y_j)x_j + s_i(y) + \theta_i - 1\\
    2\sqrt{U^{\tinyf}_i}
    \end{bmatrix}\right\|_2 \ \leq \ \sum_{j \in J} w_{ij}(1 - y_j)x_j + s_i(y) - \theta_i + 1, \quad \forall i \in I. \label{prop:soc-note-0}
\end{equation}
\end{proposition}

\proof{Proof:} By definition, $\theta \leq \widehat{L}(x, y)$ holds valid if and only if there exist $\{\theta_i\}_{i \in I}$ such that $\displaystyle \theta \leq \sum_{i \in I}h_i \theta_i$ and 
\begin{equation}
    \theta_i \ \leq \ \frac{U^{\tinyl}_i+\sum_{j\in J}w_{ij}{\bigl[- y_j x_j^2+(1+ y_j)x_j\bigr]}}{U^{\tinyl}_i+U^{\tinyf}_i+\sum_{j\in J}w_{ij}\bigl[(1-y_j)x_j+y_j\bigr]}, \quad \forall i \in I. \label{prop:soc-note-1}
\end{equation}
To show that inequality \eqref{prop:soc-note-1} is second-order conic representable, we note that $\displaystyle \sum_{j\in J}w_{ij}{\bigl[- y_j x_j^2+(1+ y_j)x_j\bigr]} = \sum_{j \in J_0} w_{ij} x_j + \sum_{j \in J_1} w_{ij}(-x_j^2 + 2x_j)$ and $\displaystyle \sum_{j\in J}w_{ij}\bigl[(1-y_j)x_j+y_j\bigr] = \sum_{j \in J_0} w_{ij} x_j + \sum_{j \in J_1} w_{ij} \cdot 1 > 0$. We finish the proof by rewriting inequality \eqref{prop:soc-note-1} as follows.
\begin{align*}
    \mbox{\eqref{prop:soc-note-1}} \ \Longleftrightarrow & \ \theta_i  \left(U^{\tinyl}_i + U^{\tinyf}_i + \sum_{j \in J_0} w_{ij} x_j + \sum_{j \in J_1} w_{ij}\right) \leq U^{\tinyl}_i + U^{\tinyf}_i + \sum_{j \in J_0} w_{ij} x_j + \sum_{j \in J_1} w_{ij}(-x_j^2 + 2x_j) - U^{\tinyf}_i \\
    \Longleftrightarrow & \ \sum_{j \in J_1}w_{ij}(1 - x_j)^2 + \left(\sqrt{U^{\tinyf}_i}\right)^2 \leq \left(U^{\tinyl}_i + U^{\tinyf}_i + \sum_{j \in J_0} w_{ij} x_j + \sum_{j \in J_1} w_{ij}\right) (1 - \theta_i) \\
    \Longleftrightarrow & \ \sum_{j \in J_1}4w_{ij}(1 - x_j)^2 + \left(2\sqrt{U^{\tinyf}_i}\right)^2 + \left(\sum_{j \in J} w_{ij} (1 - y_j) x_j + s_i(y) + \theta_i - 1\right)^2 \\
    \leq & \ \left(\sum_{j \in J} w_{ij} (1 - y_j) x_j + s_i(y) - \theta_i + 1\right)^2 \ \Longleftrightarrow \ \mbox{\eqref{prop:soc-note-0}},
\end{align*}
where the second-to-last equivalence uses the equation that $xy = \frac{1}{4}(x+y)^2 - \frac{1}{4}(x-y)^2$ for any real numbers $x$ and $y$. \Halmos
\endproof

\subsection{Proof of Lemma \ref{lem:oc}} \label{apx:lem:oc}
\proof{Proof:} First, we recall that (S-CFLP) is equivalent to $\displaystyle\max_{x \in \mathcal{X}} \min_{y \in \mathcal{Y}}L(x, y)$. Suppose that $\bar{x}$ is an optimal solution to (S-CFLP) and $e^{\top}\bar{x} < p$. Then, there exists a $k \in J$ such that $\bar{x}_k = 0$. We construct a new solution $\bar{x}'$ such that $\bar{x}'_k = 1$ and $\bar{x}_j = \bar{x}'_j$ for all $j \neq k$. Then, $\bar{x}' \in \mathcal{X}$ because $e^{\top}\bar{x}' \leq p$. In addition, for any $y \in \mathcal{Y}$, we notice that $L(\bar{x}', y) \geq L(\bar{x}, y)$ by discussing the following two cases.
\begin{enumerate}[1.]
    \item If $y_k = 1$, then $(\bar{x}_j \vee y_j) = (\bar{x}'_j \vee y_j)$ for all $j \in J$. It follows that
\begin{align*}
    L(\bar{x}', y) \ = & \ \sum_{i\in I}h_{i}\left(\frac{U^{\tinyl}_i+\sum_{j\in J\setminus\{k\}}w_{ij} \bar{x}'_j + w_{ik}}{U^{\tinyl}_i+U^{\tinyf}_i+\sum_{j\in J}w_{ij}(\bar{x}'_{j}\vee y_{j})}\right) \\
    = & \ \sum_{i\in I}h_{i}\left(\frac{U^{\tinyl}_i+\sum_{j\in J\setminus\{k\}}w_{ij} \bar{x}_j + w_{ik}}{U^{\tinyl}_i+U^{\tinyf}_i+\sum_{j\in J}w_{ij}(\bar{x}_{j}\vee y_{j})}\right) \\
    \geq & \ \sum_{i\in I}h_{i}\left(\frac{U^{\tinyl}_i+\sum_{j\in J\setminus\{k\}}w_{ij} \bar{x}_j}{U^{\tinyl}_i+U^{\tinyf}_i+\sum_{j\in J}w_{ij}(\bar{x}_{j}\vee y_{j})}\right) \ = \ L(\bar{x}, y).
\end{align*}
    \item If $y_k = 0$, then $(\bar{x}_k\vee y_k) = 0$, $(\bar{x}'_k\vee y_k) = 1$, and $(\bar{x}_k\vee y_k) - \bar{x}_k = 0 = (\bar{x}'_k \vee y_k) - \bar{x}'_k$. It follows that
\begin{align*}
    L(\bar{x}', y) \ = & \ 1 - \sum_{i\in I}h_{i}\left(\frac{U^{\tinyf}_i+\sum_{j\in J}w_{ij} \left[(\bar{x}'_j \vee y_j) - \bar{x}'_j\right]}{U^{\tinyl}_i+U^{\tinyf}_i+\sum_{j\in J\setminus\{k\}}w_{ij}(\bar{x}'_{j}\vee y_{j}) + w_{ik}(\bar{x}'_k \vee y_k)}\right) \\
    = & \ 1 - \sum_{i\in I}h_{i}\left(\frac{U^{\tinyf}_i+\sum_{j\in J}w_{ij} \left[(\bar{x}_j \vee y_j) - \bar{x}_j\right]}{U^{\tinyl}_i+U^{\tinyf}_i+\sum_{j\in J\setminus\{k\}}w_{ij}(\bar{x}_{j}\vee y_{j}) + w_{ik}}\right) \\
    \geq & \ 1 - \sum_{i\in I}h_{i}\left(\frac{U^{\tinyf}_i+\sum_{j\in J}w_{ij} \left[(\bar{x}_j \vee y_j) - \bar{x}_j\right]}{U^{\tinyl}_i+U^{\tinyf}_i+\sum_{j\in J\setminus\{k\}}w_{ij}(\bar{x}_{j}\vee y_{j})}\right) \ = \ L(\bar{x}, y).
\end{align*}
\end{enumerate}
Therefore, $\displaystyle \min_{y \in \mathcal{Y}}L(\bar{x}', y) \geq \min_{y \in \mathcal{Y}}L(\bar{x}, y)$ and so $\bar{x}'$ is also optimal. Repeating this procedure yields an optimal solution $x^*$ to (S-CFLP) such that $e^{\top}x^* = p$.

Second, for any fixed $\hat{x} \in \mathcal{X}$, suppose that $\bar{y}$ is an optimal solution to the separation problem \eqref{sp} such that $e^{\top}\bar{y} < r$. Then, we can construct a $y^* \in \mathcal{Y}$ with $e^{\top}y^* = r$ by flipping sufficiently many entries of $\bar{y}$ from zero to one. Since $(\hat{x}_j \vee \bar{y}_j) \leq (\hat{x}_j \vee y^*_j)$ for all $j \in J$, we have $L(\hat{x}, y^*) \leq L(\hat{x}, \bar{y})$. Therefore, $y^*$ is also optimal to \eqref{sp}.
\Halmos
\endproof

\subsection{Proof of Proposition \ref{prop:approx}} \label{apx-prop-approx}
\proof{Proof:}
First, we rewrite $L(\hat{x}, y)$ as
\begin{align*}
    L(\hat{x}, y) = & \sum_{i\in I}h_{i}\left(\frac{U^{\tinyl}_i+\sum_{j\in J}w_{ij} \hat{x}_j}{U^{\tinyl}_i+U^{\tinyf}_i+\sum_{j\in J}w_{ij}(\hat{x}_{j} + y_{j} - \hat{x}_j\cdot y_j)}\right) \\
    = & \sum_{i\in I}h_{i}\left(\frac{a_i(\hat{x})}{a_i(\hat{x}) + w_i(\hat{x})}\right),
\end{align*}
where $w_i(\hat{x}) := U^{\tinyf}_i+\sum_{j\in J}w_{ij}(1 - \hat{x}_j)y_j$ for all $i \in I$. Without loss of optimality, we can assume that $w^{\tinyl}_i(\hat{x}) \leq w_i(\hat{x}) \leq w_i^{\tinyu}(\hat{x})$. Since $a_i(\hat{x}) > 0$ and $w^{\tinyl}_i(\hat{x}) \geq 0$, the function $\displaystyle \ell(w) := \frac{a_i(\hat{x})}{a_i(\hat{x}) + w}$ is convex in the interval $[w_i^{\tinyl}(\hat{x}), w_i^{\tinyu}(\hat{x})]$. As a result,
\begin{align*}
    L(\hat{x}, y) \equiv \ell(w_i(\hat{x})) \leq & \sum_{i\in I}h_{i} \left[ \ell(w^{\tinyl}_i(\hat{x})) + \frac{\ell(w^{\tinyu}_i(\hat{x})) - \ell(w^{\tinyl}_i(\hat{x}))}{w^{\tinyu}_i(\hat{x}) - w^{\tinyl}_i(\hat{x})} (w_i(\hat{x}) - w^{\tinyl}_i(\hat{x})) \right] \\
    = & \sum_{i\in I}h_{i} \left[ \frac{a_i(\hat{x})(a_i(\hat{x}) + w^{\tinyu}_i(\hat{x}) + w^{\tinyl}_i(\hat{x}))}{(a_i(\hat{x}) + w^{\tinyu}_i(\hat{x}))(a_i(\hat{x}) + w^{\tinyl}_i(\hat{x}))} - \left(\frac{a_i(\hat{x})}{(a_i(\hat{x}) + w^{\tinyu}_i(\hat{x}))(a_i(\hat{x}) + w^{\tinyl}_i(\hat{x}))}\right) w_i(\hat{x}) \right] \\
    = & \ \alpha(\hat{x}) - \beta(\hat{x})^{\top}y.
\end{align*}
Second, since $\alpha(\hat{x}) - \beta(\hat{x})^{\top}y$ is affine in $y$, a greedy algorithm solves the problem $\displaystyle \min_{y \in \overline{\mathcal{Y}}}\bigl\{\alpha(\hat{x}) - \beta(\hat{x})^{\top}y\bigr\}$ with an optimal solution $\hat{y}$ as described in the claim of this proposition.
\Halmos
\endproof

\subsection{Proof of Theorem \ref{thm:heuristic}} \label{apx-thm-heuristic}
\proof{Proof:}
First, by Lemma \ref{lem:oc}, (S-CFLP) can be represented as $\displaystyle \max_{x \in \overline{\mathcal{X}}} \min_{y \in \overline{\mathcal{Y}}} L(x, y)$. For all $x \in \overline{\mathcal{X}}$, $y \in \overline{\mathcal{Y}}$, and $i \in I$, we define $\displaystyle a_i(x) := U_i^{\tinyl}+\sum_{j \in J}w_{ij}x_j$ and $\displaystyle w_i(x,y) := U_i^{\tinyf} + \sum_{j \in J}w_{ij}(1-x_j)y_j$. Then, $L(x, y) = \sum_{i \in I}h_i\left(\frac{a_i(x)}{a_i(x)+w_i(x,y)}\right)$. We bound $L(x, y)$ from below by the harmonic-geometric-arithmetic mean inequality:
\begin{align}
L(x,y) \ \geq \ \frac{1}{\sum_{i \in I}h_i\left(\frac{a_i(x)+w_i(x,y)}{a_i(x)}\right)}, \quad \forall x \in \overline{\mathcal{X}}, \ \forall y \in \overline{\mathcal{Y}}. \label{heu-note-5}
\end{align}
It follows that
\begin{subequations}
\begin{align*}
\max_{x \in \overline{\mathcal{X}}} \min_{y \in \overline{\mathcal{Y}}} L(x,y) \ & \geq \  \max_{x \in \overline{\mathcal{X}}} \min_{y \in \overline{\mathcal{Y}}} \frac{1}{\sum_{i \in I}h_i\left(\frac{a_i(x)+w_i(x,y)}{a_i(x)}\right)}\\
& = \ \frac{1}{\displaystyle 1+\min_{x \in \overline{\mathcal{X}}} \max_{y \in \overline{\mathcal{Y}}} \sum_{i \in I}\left(\frac{h_i w_i(x,y)}{a_i(x)}\right)},
\end{align*}
\end{subequations}
where the equality uses the facts that $\sum_{i \in I}h_i = 1$ and that function $1/x$ is decreasing in $x$ on $\mathbb{R}_+$. Hence,  solving problem
\begin{equation}
\min_{x \in \overline{\mathcal{X}}} \max_{y \in \overline{\mathcal{Y}}}\sum_{i \in I}\left(\frac{h_i w_i(x,y)}{a_i(x)}\right) \label{heu-note-1}
\end{equation}
produces a feasible solution and a lower bound of (S-CFLP).

Second, we recast formulation \eqref{heu-note-1} as the MISOCP \eqref{heu}. To do this, we derive
\begin{subequations}
\begin{align}
& \ \min_{x \in \overline{\mathcal{X}}} \max_{y \in \overline{\mathcal{Y}}}\sum_{i \in I}\left(\frac{h_i w_i(x,y)}{a_i(x)}\right) \nonumber \\
= & \ \min_{x \in \overline{\mathcal{X}}} \max_{y \in \overline{\mathcal{Y}}}\left\{\sum_{i \in I}\left(\frac{h_i U^{\tinyf}_i}{U^{\tinyl}_i+\sum_{k \in J} w_{ik}x_k}\right)+\sum_{j \in J}\left(\sum_{i \in I}\frac{h_i w_{ij}(1-x_j)}{U^{\tinyl}_i+\sum_{k \in J}w_{ik}x_k}\right)y_j\right\} \label{heu-note-2} \\
= & \ \min_{x \in \overline{\mathcal{X}}, \ \mu \geq 0} \ \ \sum_{i \in I}\left(\frac{h_i U^{\tinyf}_i}{U^{\tinyl}_i+\sum_{k \in J}w_{ik}x_k}\right) + r\lambda + e^{\top}\mu \label{heu-note-3} \\
& \qquad \mbox{s.t.} \quad \ \ \lambda + \mu_j \geq \sum_{i \in I} \left(\frac{h_i w_{ij}(1-x_j)}{U^{\tinyl}_i+\sum_{k \in J}w_{ik}x_k}\right), \quad \forall j \in J \nonumber \\
= & \ \min_{\substack{x \in \overline{\mathcal{X}}, \ \mu \geq 0,\\ s \geq 0, \ t \geq 0}} \ \ \sum_{i \in I} h_i U^{\tinyf}_i s_i + r\lambda + e^{\top}\mu \label{heu-note-4} \\
& \qquad \mbox{s.t.} \quad \ \ \lambda + \mu_j \geq \sum_{i \in I} h_i w_{ij}t_{ij}, \quad \forall j \in J, \\
& \qquad \qquad \quad s_i \geq \frac{1}{U_i^{\tinyl}+\sum_{k \in J} w_{ik}x_k}, \quad \forall i \in I,  \label{heu-1c}\\
& \qquad \qquad \quad t_{ij} \geq \frac{1-x_j}{U_i^{\tinyl}+\sum_{k \in J} w_{ik}x_k}, \quad \forall i \in I, \forall j \in J,\label{heu-1d}
\end{align}
\end{subequations}
where equality \eqref{heu-note-3} is due to the strong duality of linear programming. Indeed, since the inner maximization problem in \eqref{heu-note-2} has an objective function linear in $y$, the feasible region $\overline{\mathcal{Y}}$ of $y$ can be replaced by its convex hull $\mbox{conv}(\overline{\mathcal{Y}}) = \{y \in [0, 1]^{|J|}: e^{\top}y = r\}$, because $e^{\top}y = r$ produces a totally unimodular constraint matrix. Then, this inner maximization problem is equivalent to a linear program with a feasible region $y \in \{y \in \mathbb{R}^{|J|}_+: e^{\top}y = r, \ y \leq e\}$. Taking the dual of this linear program yields formulation \eqref{heu-note-4}--\eqref{heu-1d}, where dual variables $\lambda$ and $\mu$ are associated with primal constraints $e^{\top}y = r$ and $y \leq e$, respectively. We represent constraints \eqref{heu-1c} as
\begin{subequations}
\begin{align*}
\mbox{\eqref{heu-1c}} \ \iff \ & \left(s_i+U^{\tinyl}_i+\sum_{k \in J}w_{ik}x_k\right)^2 - \left(U^{\tinyl}_i+\sum_{k \in J}w_{ik}x_k-s_i\right)^2 \geq 4 \\
\iff \ & 
\begin{Vmatrix}
\begin{bmatrix} 
2 \\
U^{\tinyl}_i+\sum_{k \in J}w_{ik}x_k - s_i 
\end{bmatrix}
\end{Vmatrix} \leq s_i+U^{\tinyl}_i+\sum_{k \in J}w_{ik}x_k \ \iff \ \mbox{\eqref{heu-1c-ref}}.
\end{align*}
\end{subequations}
Likewise, we recast constraints \eqref{heu-1d} as
\begin{align*}
\mbox{\eqref{heu-1d}} \ \iff \ & t_{ij} \geq \frac{(1-x_j)^2}{U_i^{\tinyl}+\sum_{k \in J} w_{ik}x_k} \\
\iff \ & \begin{Vmatrix}
\begin{bmatrix} 
2(1-x_j)  \\
U^{\tinyl}_i+\sum_{k \in J} w_{ik}x_k-t_{ij} 
\end{bmatrix}
\end{Vmatrix} \leq t_{ij}+U^{\tinyl}_i+\sum_{k \in J} w_{ik}x_k \ \iff \ \mbox{\eqref{heu-1d-ref}},
\end{align*}
where the first equivalence is because $x_j$ is binary-valued.

Third, let $\ell(x, y) := \displaystyle \left[\sum_{i \in I}h_i\left(\frac{a_i(x)+w_i(x,y)}{a_i(x)}\right)\right]^{-1}$. Then, $\displaystyle x^{\mbox{\tiny H}} \in \argmax_{x \in \overline{\mathcal{X}}}\left\{\min_{y \in \overline{\mathcal{Y}}} \ell(x, y)\right\}$ by the equivalence between \eqref{heu-note-1} and \eqref{heu}. In addition, $L(x, y) \geq \ell(x, y)$ by inequality \eqref{heu-note-5}. We bound $L(x, y)$ from above by deriving
\begin{align*}
L(x,y) \ & \equiv \ \sum_{i \in I}h_i\left(\frac{a_i(x)}{a_i(x)+w_i(x,y)}\right) \\
& \leq \ \frac{(\gamma_M(x, y)+\gamma_m(x, y))^2}{4\gamma_M(x, y) \gamma_m(x, y)}\left[\sum_{i \in I}h_i\left(\frac{a_i(x)+w_i(x,y)}{a_i(x)}\right)\right]^{-1} \\
& \leq \ \frac{(\gamma_M+\gamma_m)^2}{4\gamma_M \gamma_m} \ell (x, y),
\end{align*}
where $\displaystyle \gamma_M(x, y) := \max_{i \in I}\left\{\frac{a_i(x)}{a_i(x)+w_i(x,y)}\right\}$ and $\displaystyle \gamma_m(x, y) := \min_{i \in I}\left\{\frac{a_i(x)}{a_i(x)+w_i(x,y)}\right\}$. The first inequality above follows from the Kantorovich inequality~\citep[see][]{kantorovich1945functional,henrici1961two}, and the second inequality is because
\begin{align*}
\frac{(\gamma_M(x, y)+\gamma_m(x, y))^2}{4\gamma_M(x, y) \gamma_m(x, y)} \ & = \ \frac{1}{4}\left(\frac{\gamma_M(x, y)}{\gamma_m(x, y)} + \frac{\gamma_m(x, y)}{\gamma_M(x, y)} + 2\right) \\
& \leq \ \frac{1}{4}\left(\frac{\gamma_M}{\gamma_m} + \frac{\gamma_m}{\gamma_M} + 2\right) \ = \ \frac{(\gamma_M+\gamma_m)^2}{4\gamma_M \gamma_m},
\end{align*}
where the inequality is because $1 \leq \frac{\gamma_M(x, y)}{\gamma_m(x, y)} \leq \frac{\gamma_M}{\gamma_m}$ and function $f(z) := z + \frac{1}{z} + 2$ is increasing on the interval $[1, \infty)$. 

Now pick any $\displaystyle x^* \in \argmax_{x \in \overline{\mathcal{X}}}\left\{\min_{y \in \overline{\mathcal{Y}}}L(x, y)\right\}$. Then, $\displaystyle z^{\mbox{\tiny H}} \equiv \min_{y \in \overline{\mathcal{Y}}}L(x^{\mbox{\tiny H}}, y) \leq \min_{y \in \overline{\mathcal{Y}}}L(x^*, y) = z^*$ by definition of $x^*$. In addition,
\begin{align*}
z^{\mbox{\tiny H}} \equiv \min_{y \in \overline{\mathcal{Y}}}L(x^{\mbox{\tiny H}}, y) \ \geq \ \min_{y \in \overline{\mathcal{Y}}}\ell(x^{\mbox{\tiny H}}, y) \ \geq \ \min_{y \in \overline{\mathcal{Y}}}\ell(x^*, y) \ \geq \ \min_{y \in \overline{\mathcal{Y}}} \left\{\frac{4\gamma_M \gamma_m}{(\gamma_M+\gamma_m)^2} L(x^*, y)\right\} = \frac{4\gamma_M \gamma_m}{(\gamma_M+\gamma_m)^2} z^*,
\end{align*}
where the first and the third inequalities are because $\ell(x, y) \leq L(x, y) \leq \frac{(\gamma_M+\gamma_m)^2}{4\gamma_M \gamma_m} \ell(x, y)$, and the second inequality is by definition of $x^{\mbox{\tiny H}}$. This finishes the proof.
\Halmos
\endproof

\subsection{RO Approximations in Section~\ref{sec:extension:outside}} \label{apx-ext-outside}
\begin{proposition} \label{prop:ext-outside}
It holds that
\begin{equation*}
\max_{x \in \mathcal{X}} \min_{y \in \mathcal{Y}(x)} L^+(x, y) \leq z^{\tinyo} \leq \max_{x \in \mathcal{X}} \min_{y \in \mathcal{Y}(x)} \big\{1-F^+(x, y)\big\}.
\end{equation*}
In addition, it holds that
\begin{align*}
\max_{x \in \mathcal{X}} \min_{y \in \mathcal{Y}(x)} L^+(x, y) \ & = \ \max_{x \in \mathcal{X}} \min_{y \in \mathcal{Y}} \sum_{i\in I} h_i \left(\frac{U^{\tinyl}_i + \sum_{j \in J}w_{ij}x_j}{U^{\tinyl}_i+U^{\tinyf}_i+\sum_{j \in J}w_{ij}(x_{j} \vee y_{j}) + U^{\tinyo}_i}\right) \\[1em]
\text{and} \quad \max_{x \in \mathcal{X}} \min_{y \in \mathcal{Y}(x)}\big\{1 - F^+(x,y)\big\} \ & = \ \max_{x \in \mathcal{X}} \min_{y \in \mathcal{Y}} \sum_{i\in I} h_i \left(\frac{U^{\tinyl}_i + \sum_{j \in J}w_{ij}x_j + U^{\tinyo}_i}{U^{\tinyl}_i+U^{\tinyf}_i+\sum_{j \in J}w_{ij}(x_{j} \vee y_{j}) + U^{\tinyo}_i}\right).
\end{align*}
\end{proposition}
\proof{Proof of Proposition~\ref{prop:ext-outside}:}
First, pick any $x \in \mathcal{X}$. By definition, each $y^* \in \argmax_{y \in \mathcal{Y}(x)}F^+(x,y)$ is a feasible solution for the formulation $\min_{y \in \mathcal{Y}(x)}L^+(x,y)$. Hence, it holds that $\min_{y \in \mathcal{Y}(x)}L^+(x,y) \leq L^+(x,y^*)$ for all $x \in \mathcal{X}$ and it follows that $\max_{x \in \mathcal{X}}\min_{y \in \mathcal{Y}(x)} L^+(x, y) \leq z^{\tinyo}$.

Second, we notice that $L^+(x,y) \leq 1 - F^+(x,y)$. Since $\argmax_{y \in \mathcal{Y}(x)}F^+(x,y) = \argmin_{y \in \mathcal{Y}(x)}\big\{1-F^+(x,y)\big\}$ by definition, it holds that $L^+(x,y^*) \leq \min_{y \in \mathcal{Y}(x)} \big\{1-F^+(x,y)\big\}$ for any $x \in \mathcal{X}$ and $y^* \in \argmax_{y \in \mathcal{Y}(x)}F^+(x,y)$. It follows that $z^{\tinyo} \leq \max_{x \in \mathcal{X}} \min_{y \in \mathcal{Y}(x)} \big\{1-F^+(x, y)\big\}$.

Third, it follows from a similar proof to that of Theorem~\ref{lemma:l1} that replacing $(x_j + y_j)$ with $(x_j \vee y_j)$ in definitions of $L^+(x,y)$ and $1 - F^+(x,y)$ allows us to remove the decision dependency (i.e., replacing $\mathcal{Y}(x)$ with $\mathcal{Y}$) in the two RO formulations without loss of optimality.
\Halmos
\endproof

\subsection{Reformulation of (S-CFLP-U) in Section~\ref{sec:extension:change}} \label{apx-prop:extension:change}
\begin{proposition} \label{prop:ext:change}
It holds that
\begin{equation*}
\max_{x \in \mathcal{X}} \min_{y \in \mathcal{Y}(x)} L^{\textrm{U}}(x, y) \ = \ \max_{x \in \mathcal{X}} \min_{y \in \mathcal{Y}} L^{\textrm{U}}(x, y).
\end{equation*}
In addition, for any $y \in \mathcal{Y}$, define $t_i(x) := U^{\tinyl}_i + U^{\tinyf}_i +  \sum_{j \in J} (v^{\tinyl}_{ij} + v^{\tinyf}_{ij} + w_{ij})\big[(1 - y_j)x_j + y_j\big]$. Then, the hypograph of the objective function, i.e., the inequality $\theta \leq L^{\textrm{U}}(x, y)$, holds if and only if there exist $\{\theta_i\}_{i \in I}$ such that $\theta \leq \sum_{i \in I}h_i \theta_i$ and
\begin{align}
    & \ \left\|\begin{bmatrix}
    2\sqrt{w_{ij}}(1 - x_j)_{j \in J_1}\\
    \left(2\sqrt{v^{\tinyf}_{ij}}x_j\right)_{j \in J_0} \\
    \left(2\sqrt{v^{\tinyf}_{ij}}\right)_{j \in J_1} \\
    2\sqrt{U^{\tinyf}_i} \\
    t_i(x) + \theta_i - 1
    \end{bmatrix}\right\|_2 \ \leq \ t_i(x) - \theta_i + 1, \quad \forall i \in I, \label{eq:prop:ext:change}
\end{align}
where $J_0 = \{j \in J: y_j = 0\}$ and $J_1 = \{j \in J: y_j = 1\}$.
\end{proposition}
\proof{Proof of Proposition~\ref{prop:ext:change}:}
First, every $y_j$ in the definition of $L^{\textrm{U}}(x, y)$ appears in the form $(x_j \vee y_j)$. As a result, for any $x \in \mathcal{X}$ and for any $j \in J$ such that $x_j = 1$, the value of the corresponding $y_j$ does not affect the value of $L^{\textrm{U}}(x, y)$. That is because in this case $(x_j \vee y_j) = (1 \vee y_j) = 1$. It follows that we can relax the constraints $y_j \leq 1 - x_j$, for all $j \in J$, from the set $\mathcal{Y}(x)$ without loss of optimality. Therefore, it holds that $\min_{y \in \mathcal{Y}(x)}L^{\textrm{U}}(x, y) = \min_{y \in \mathcal{Y}}L^{\textrm{U}}(x, y)$.

Second, by definition, $\theta \leq L^{\textrm{U}}(x, y)$ holds if and only if there exist $\{\theta_i\}_{i \in I}$ such that $\theta \leq \sum_{i \in I}h_i \theta_i$ and
\begin{align}
    \theta_i \ \leq & \ \frac{U^{\tinyl}_i + \sum_{j \in J} v^{\tinyl}_{ij}(x_j \vee y_j) + \sum_{j \in J}w_{ij}x_j}{U^{\tinyl}_i+U^{\tinyf}_i+\sum_{j \in J}(v^{\tinyl}_{ij} + v^{\tinyf}_{ij} + w_{ij})(x_{j} \vee y_{j})} \nonumber \\[1em]
    = & \ 1 - \frac{U^{\tinyf}_i + \sum_{j \in J} v^{\tinyf}_{ij}(x_j \vee y_j) + \sum_{j \in J}w_{ij}(x_j \vee y_j - x_j)}{U^{\tinyl}_i+U^{\tinyf}_i+\sum_{j \in J}(v^{\tinyl}_{ij} + v^{\tinyf}_{ij} + w_{ij})(x_j \vee y_j)} \nonumber \\[1em]
    = & \ 1 - \frac{U^{\tinyf}_i + \sum_{j \in J_0} v^{\tinyf}_{ij}x_j^2 + \sum_{j \in J_1} v^{\tinyf}_{ij} + \sum_{j \in J_1}w_{ij}(1 - x_j)^2}{U^{\tinyl}_i+U^{\tinyf}_i+\sum_{j \in J}(v^{\tinyl}_{ij} + v^{\tinyf}_{ij} + w_{ij})\big[(1 - y_j)x_j + y_j\big]}, \quad \forall i \in I, \label{prop:ext:change-note-1}
\end{align}
where the last equality uses the facts that $x_j = x_j^2$ and $(1 - x_j) = (1 - x_j)^2$ since $x_j \in \{0,1\}$. We finish the proof by rewriting inequality \eqref{prop:ext:change-note-1} as follows.
\begin{align*}
\mbox{\eqref{prop:ext:change-note-1}} \ \Longleftrightarrow & \ U^{\tinyf}_i + \sum_{j \in J_0} v^{\tinyf}_{ij}x_j^2 + \sum_{j \in J_1} v^{\tinyf}_{ij} + \sum_{j \in J_1}w_{ij}(1 - x_j)^2 \leq (1 - \theta_i) t_i(x) \\
\Longleftrightarrow & \ \Big(2\sqrt{U^{\tinyf}_i}\Big)^2 + \sum_{j \in J_0} \Big(2\sqrt{v^{\tinyf}_{ij}}x_j\Big)^2 + \sum_{j \in J_1} \Big(2\sqrt{v^{\tinyf}_{ij}}\Big)^2 + \sum_{j \in J_1}\Big(2\sqrt{w_{ij}}(1 - x_j)\Big)^2 \\
\leq & \ \big(t_i(x) - \theta_i + 1 \big)^2 - \big(t_i(x) + \theta_i - 1\big)^2 \\
\Longleftrightarrow & \ \mbox{\eqref{eq:prop:ext:change}},
\end{align*}
where the second-to-last equivalence uses the equation that $xy = \frac{1}{4}(x+y)^2 - \frac{1}{4}(x-y)^2$ for any real numbers $x$ and $y$.
\Halmos
\endproof

\subsection{Exponential Increase of the (S-CFLP) Search Space in $p$ and $r$} \label{apx-prop:exponential}
\begin{proposition} \label{prop:exponential}
The search space of the (S-CFLP) problem has a cardinality $C(|J|,p) \times C(|J|-p,r)=\dbinom{|J|}{p} \times \dbinom{|J|-p}{r}$. In addition, when $p\le |J|/2$ and $r\le(|J|-p)/2$, it holds that $C(|J|,p) \times C(|J|-p,r) \geq 2^{p+r}$, that is, (S-CFLP)'s search space increases exponentially in $p$ and $r$.
\end{proposition}
\proof{Proof of Proposition~\ref{prop:exponential}:} The cardinality of the (S-CFLP) search space follows from constraints~\eqref{bl-b}, \eqref{bl-e} and the definition of combinatorial numbers $C(|J|,p)$, $C(|J|-p,r)$. In addition,
\begin{align*}
C(|J|,p) \times C(|J|-p,r) \ = & \ \frac{|J|!}{p!(|J|-p)!} \times \frac{(|J|-p)!}{r!(|J|-p-r)!} \\
= & \ \frac{|J|(|J|-1)\cdots(|J|-p-r+1)}{[p(p-1)\cdots 1][r(r-1)\cdots 1]} \\
= & \ \frac{|J|(|J|-1)\cdots(|J|-p+1)}{p(p-1)\cdots 1} \times \frac{(|J|-p)(|J|-p-1)\cdots(|J|-p-r+1)}{r(r-1)\cdots 1} \\
\geq & \ \left(\frac{|J|}{p}\right)^p \left(\frac{|J|-p}{r}\right)^r \ \geq \ 2^{p+r},
\end{align*}
where the first inequality is because $|J|/p \leq (|J|-k)/(p-k)$ for all $k \in [0, p-1]$ and $(|J|-p)/r \leq (|J|-p-k)/(r-k)$ for all $k \in [0, r-1]$, and the second inequality follows from the assumption that $p\le |J|/2$ and $r\le(|J|-p)/2$.
\Halmos
\endproof

\section{Valid Inequalities for MISOCP \eqref{heu}} \label{apx-vi-misocp}
We derive valid linear inequalities with regard to the second-order conic constraints \eqref{heu-1c-ref}--\eqref{heu-1d-ref} (more specifically, their representation \eqref{heu-1c}--\eqref{heu-1d}) to further improve the efficacy of solving MISOCP \eqref{heu}. These valid inequalities exploit the convexity of the right-hand side of inequality \eqref{heu-1c}, and although that of \eqref{heu-1d} is \emph{non}-convex, we recover the convexity by resorting to perspective function.
\begin{proposition} \label{apx-prop-vi}
Given any $\hat{x} \in \mathcal{X}$, the following linear inequalities are valid for MISOCP \eqref{heu}:
\begin{subequations}
\begin{align}
s_i \ \geq \ & \frac{U_i^{\tinyl} + \sum_{k \in J}w_{ik}(2\hat{x}_k - x_k)}{(U_i^{\tinyl}+\sum_{k \in J} w_{ik}\hat{x}_k)^2}, \quad \forall i \in I, \label{apx-vi-s} \\
t_{ij} \ \geq \ & \frac{(1 - \hat{x}_j)\left( (U_i^{\tinyl}+2\sum_{k \neq j} w_{ik}\hat{x}_k)(1 - x_j) - \sum_{k \neq j} w_{ik} x_k \right)}{(U_i^{\tinyl}+\sum_{k \neq j} w_{ik}\hat{x}_k)^2}, \quad \forall i \in I, \forall j \in J. \label{apx-vi-t}
\end{align}
\end{subequations}
\end{proposition}
\proof{Proof of Proposition \ref{apx-prop-vi}:}
From the proof of Theorem \ref{thm:heuristic}, we observe that the second-order conic constraints \eqref{heu-1c-ref} and \eqref{heu-1d-ref} are equivalent to inequalities \eqref{heu-1c} and \eqref{heu-1d}, respectively.

First, we notice that the right-hand side of inequality \eqref{heu-1c} is convex in variables $x$ because it is the composite of function $c(z) := 1/z$ on $z > 0$ and an affine function $z(x) := U_i^{\tinyl}+\sum_{k \in J} w_{ik} x_k$. Then, the supporting hyperplane of this convex function yields inequalities \eqref{apx-vi-s}.

Second, we notice that the right-hand side of inequality \eqref{heu-1d} is not convex in variables $x$. Nevertheless, for each $j \in J$, this right-hand side is convex in variables $(x_1, \ldots, x_{j-1}, x_{j+1}, \ldots, x_{|J|})$ when fixing $x_j = 0$. In contrast, if we fix $x_j = 1$ then this right-hand side becomes zero. This motivates us to replace the right-hand side of \eqref{heu-1d} with the perspective function of its restriction generated by fixing $x_j = 0$. More specifically, we consider inequalities
\begin{equation*}
t_{ij} \ \geq \ t_{ij}(x) := \frac{(1-x_j)^2}{U_i^{\tinyl}(1 - x_j)+\sum_{k \neq j} w_{ik}x_k}, \quad \forall i \in I, \forall j \in J.
\end{equation*}
We notice that this inequality is equivalent to \eqref{heu-1d} whenever $x \in \mathcal{X} \subseteq \{0,1\}^{|J|}$, but its right-hand side $t_{ij}(x)$ is now convex in variables $x$ because it is the perspective function of $1/(U_i^{\tinyl}+\sum_{k \neq j} w_{ik}x_k)$, which is convex in $x$. Then, the supporting hyperplane of $t_{ij}(x)$ yields inequalities \eqref{apx-vi-t}. This finishes the proof.
\Halmos
\endproof

\section{Additional Computational Results} \label{apx-compu}

\subsection{Detailed Results of Optimality Gap Improvements} 
\label{apx-fig-gaps}

In Table~\ref{tab:fig-gaps}, we report the optimality gap data used to depict Figure~\ref{fig:gaps}.

% Table generated by Excel2LaTeX from sheet 'Gap data'
\begin{table}[ht!]
\fontsize{9}{11}\selectfont
  \centering
  \caption{Detailed optimality gap improvements after the first, third, and tenth rounds of cut generation}
  \resizebox{0.87\textwidth}{!}{
    \begin{tabular}{lrrrrrrrrr}
    \hline
    \multicolumn{1}{c}{\multirow{2}[2]{*}{Instance}} & \multicolumn{3}{c}{SC} & \multicolumn{3}{c}{BI} & \multicolumn{3}{c}{SCBI} \bigstrut[t]\\
          & \multicolumn{1}{c}{\texttt{Gap$_1$}} & \multicolumn{1}{c}{\texttt{Gap$_3$}} & \multicolumn{1}{c}{\texttt{Gap$_{10}$}} & \multicolumn{1}{c}{\texttt{Gap$_1$}} & \multicolumn{1}{c}{\texttt{Gap$_3$}} & \multicolumn{1}{c}{\texttt{Gap$_{10}$}} & \multicolumn{1}{c}{\texttt{Gap$_1$}} & \multicolumn{1}{c}{\texttt{Gap$_3$}} & \multicolumn{1}{c}{\texttt{Gap$_{10}$}} \\
    \hline
    20-20-2-2 & 5.84\% & 5.84\% & 0.00\% & 36.13\% & 0.00\% & 0.00\% & 5.84\% & 0.00\% & 0.00\% \\
    20-20-3-2 & 18.45\% & 2.96\% & 2.96\% & 18.71\% & 0.86\% & 0.00\% & 18.45\% & 2.96\% & 0.00\% \\
    20-20-2-3 & 0.00\% & 0.00\% & 0.00\% & 37.38\% & 0.00\% & 0.00\% & 0.00\% & 0.00\% & 0.00\% \\
    40-40-2-2 & 0.00\% & 0.00\% & 0.00\% & 19.91\% & 14.78\% & 0.00\% & 0.00\% & 0.00\% & 0.00\% \\
    40-40-3-2 & 17.82\% & 6.75\% & 2.19\% & 21.84\% & 2.19\% & 0.47\% & 17.82\% & 0.47\% & 0.47\% \\
    40-40-2-3 & 3.65\% & 3.49\% & 0.00\% & 20.73\% & 3.05\% & 1.99\% & 3.65\% & 0.06\% & 0.00\% \\
    60-60-2-2 & 1.20\% & 0.00\% & 0.00\% & 24.48\% & 9.18\% & 0.93\% & 1.20\% & 1.20\% & 0.00\% \\
    60-60-3-2 & 100.00\% & 1.13\% & 0.37\% & 11.26\% & 0.15\% & 0.00\% & 11.26\% & 0.00\% & 0.00\% \\
    60-60-2-3 & 0.21\% & 0.00\% & 0.00\% & 25.16\% & 7.58\% & 0.00\% & 0.21\% & 0.21\% & 0.00\% \\
    80-80-2-2 & 0.10\% & 0.10\% & 0.10\% & 5.39\% & 2.96\% & 2.41\% & 0.10\% & 0.10\% & 0.05\% \\
    80-80-3-2 & 100.00\% & 100.00\% & 0.00\% & 3.02\% & 0.89\% & 0.89\% & 3.02\% & 1.27\% & 1.27\% \\
    80-80-2-3 & 1.14\% & 1.14\% & 0.15\% & 5.53\% & 1.50\% & 0.37\% & 1.14\% & 0.37\% & 0.00\% \\
    100-100-2-2 & 0.00\% & 0.00\% & 0.00\% & 9.94\% & 2.18\% & 0.15\% & 0.00\% & 0.00\% & 0.00\% \\
    100-100-3-2 & 100.00\% & 100.00\% & 0.00\% & 7.30\% & 2.80\% & 1.03\% & 7.30\% & 0.69\% & 0.60\% \\
    100-100-2-3 & 6.02\% & 6.02\% & 3.02\% & 10.07\% & 7.81\% & 0.17\% & 6.02\% & 2.49\% & 0.55\% \\
    \hline
    \textbf{Average} & 23.63\% & 15.16\% & 0.59\% & 17.12\% & 3.73\% & 0.56\% & 5.07\% & 0.65\% & 0.20\% \\
    \hline
    \end{tabular}%
    }
  \label{tab:fig-gaps}%
\end{table}%

\subsection{Effects of Approximate Separation}
\label{sec:compu:appr}
We evaluate the approximate separation method proposed in Section \ref{sec:cuts-approx}. That is, in line 6 of Algorithm \ref{algo-b&c}, we always solve the approximate separation problem first and, only when this fails to find a cut, we invoke the ``exact'' separation to guarantee optimality of the algorithm. We denote these implementations as SC-AS, BI-AS, and SCBI-AS corresponding to SC, BI, and SCBI in Section \ref{sec:perform}. Table \ref{Table:AS1} reports the results, where \texttt{Imp(\%)} represents the relative improvement on solution time over the benchmark implementations. From this table, we observe that approximate separation speeds up Algorithm \ref{algo-b&c} significantly. For example, even for the most effective SCBI implementation, approximate separation provides a $2\times$ speedup on average. We repeat this experiment on the more challenging instances with larger $p$ and $r$ values. The results are reported in Tables \ref{Table:AS2}--\ref{Table:AS3} and the observations are similar. We also observe that the advantage of SCBI over BI is lost after the approximate separation speeds up both implementations. This can be explained by that the approximate separation weakens the depth of SC and BI cuts, which in turn weakens their complementarity in SCBI-AS. 
\begin{table}[ht!]
\fontsize{9}{11}\selectfont
  \centering
  \caption{Effectiveness of the approximate separation procedure for speeding up cut generation in Algorithm \ref{algo-b&c}}
   \resizebox{\textwidth}{!}{
   \begin{tabular}{lrrrrrrrrr}
    \toprule
    \multicolumn{1}{c}{\multirow{2}[4]{*}{Instance}} & \multicolumn{1}{c}{SC} & \multicolumn{2}{c}{SC-AS} & \multicolumn{1}{c}{BI} & \multicolumn{2}{c}{BI-AS} & \multicolumn{1}{c}{SCBI} & \multicolumn{2}{c}{SCBI-AS} \\
\cmidrule{2-10}          & \multicolumn{1}{c}{\texttt{Time(s)}} & \multicolumn{1}{c}{\texttt{Time(s)}} & \multicolumn{1}{c}{\texttt{Imp(\%)}} & \multicolumn{1}{c}{\texttt{Time(s)}} & \multicolumn{1}{c}{\texttt{Time(s)}} & \multicolumn{1}{c}{\texttt{Imp(\%)}} & \multicolumn{1}{c}{\texttt{Time(s)}} & \multicolumn{1}{c}{\texttt{Time(s)}} & \multicolumn{1}{c}{\texttt{Imp(\%)}} \\
    \midrule
    20-20-2-2 & 0.94  & 0.13  & 650.40 & 2.14  & 0.11 & 1864.22 & 0.75  & 0.13  & 500.00 \\
    20-20-3-2 & 1.66  & 0.20  & 715.76 & 3.09  & 0.23 & 1221.79 & 0.97 & 0.27  & 264.29 \\
    20-20-2-3 & 2.00  & 0.47  & 327.35 & 2.42  & 0.28 & 761.92 & 1.23 & 0.38  & 245.87 \\
    40-40-2-2 & 13.23  & 5.28  & 150.60 & 4.06  & 2.84 & 42.86 & 3.75  & 2.05  & 83.19 \\
    40-40-3-2 & 68.09  & 22.89  & 197.47 & 15.28 & 4.89 & 212.43 & 11.16 & 7.59  & 46.92 \\
    40-40-2-3 & 64.44  & 17.80  & 262.07 & 20.88 & 5.81 & 259.17 & 11.02 & 6.66  & 65.50 \\
    60-60-2-2 & 69.95  & 36.25  & 92.97 & 29.36 & 7.97 & 268.43 & 9.67 & 5.53  & 74.87 \\
    60-60-3-2 & 777.56  & 563.89  & 37.89 & 79.88 & 25.06 & 218.70 & 39.33 & 29.44 & 33.60 \\
    60-60-2-3 & 640.11  & 148.75  & 330.33 & 211.22 & 40.27 & 424.57 & 94.92 & 33.03 & 187.37 \\
    80-80-2-2 & 353.55  & 244.41  & 44.66 & 65.49 & 30.55 & 114.37 & 25.75 & 24.05 & 7.09 \\
    80-80-3-2 & 13655.10  & 6366.08  & 114.50 & 147.78 & 57.36 & 157.64 & 146.78 & 83.17 & 76.48 \\
    80-80-2-3 & 5181.42  & 2848.00  & 81.93 & 384.99 & 238.14 & 61.66 & 228.08 & 187.55 & 21.61 \\
    100-100-2-2 & 636.63  & 356.59  & 78.53 & 57.97 & 50.42 & 14.97 & 44.95 & 47.77 & -5.89 \\
    100-100-3-2 & 13418.00  & 14942.20  & -10.20 & 233.02 & 102.11 & 128.20 & 190.53 & 175.80 & 8.38 \\
    100-100-2-3 & 5469.91  & 4593.98  & 19.07 & 384.00 & 299.56 & 28.19 & 273.86 & 260.83 & 5.00 \\
    \hline
    \textbf{Average} & 2690.17  & 2009.79  & \textbf{206.22} & 109.44  & 57.71  & \textbf{385.28} & 72.19  & 57.61  & \textbf{107.62} \\
    \bottomrule
    \end{tabular}%
    }
  \label{Table:AS1}%
\end{table}%

\begin{table}[ht!]
\fontsize{9}{11}\selectfont
  \centering
  \caption{Effectiveness of approximation separation on instances with $|I|=|J|=20$ and varying $p$-, $r$-values}
  \resizebox{0.9\textwidth}{!}{
  \begin{tabular}{lrrrrrrrrr}
    \toprule
    \multicolumn{1}{c}{\multirow{2}[4]{*}{Instance}} & \multicolumn{1}{c}{SC} & \multicolumn{2}{c}{SC-AS} & \multicolumn{1}{c}{BI} & \multicolumn{2}{c}{BI-AS} & \multicolumn{1}{c}{SCBI} & \multicolumn{2}{c}{SCBI-AS} \\
\cmidrule{2-10}          & \multicolumn{1}{c}{\texttt{Time(s)}} & \multicolumn{1}{c}{\texttt{Time(s)}} & \multicolumn{1}{c}{\texttt{Imp(\%)}} & \multicolumn{1}{c}{\texttt{Time(s)}} & \multicolumn{1}{c}{\texttt{Time(s)}} & \multicolumn{1}{c}{\texttt{Imp(\%)}} & \multicolumn{1}{c}{\texttt{Time(s)}} & \multicolumn{1}{c}{\texttt{Time(s)}} & \multicolumn{1}{c}{\texttt{Imp(\%)}} \\
    \midrule
    20-20-2-2 & 0.97  & 0.16 & 521.15 & 0.88  & 0.13 & 600.00 & 0.56  & 0.13 & 349.60 \\
    20-20-4-2 & 2.70  & 0.89 & 203.37 & 2.03  & 0.28 & 622.78 & 1.16  & 0.75  & 54.13 \\
    20-20-6-2 & 6.61  & 3.63 & 82.34 & 1.69  & 0.67 & 151.19 & 1.75  & 1.25  & 40.00 \\
    20-20-8-2 & 3.58  & 2.30 & 55.77 & 2.45  & 0.61  & 302.13 & 1.39  & 0.81 & 71.18 \\
    20-20-10-2 & 2.67  & 1.36 & 96.62 & 1.91  & 0.80 & 139.15 & 1.44  & 0.81 & 76.88 \\
    20-20-2-4 & 3.09  & 0.99 & 214.11 & 1.67  & 0.44 & 281.74 & 1.34  & 0.50   & 168.80 \\
    20-20-4-4 & 9.14  & 3.84 & 137.77 & 4.19  & 0.56 & 643.87 & 3.88  & 1.48 & 161.12 \\
    20-20-6-4 & 13.27  & 5.08 & 161.24 & 6.83  & 1.50   & 355.20 & 6.05  & 3.02 & 100.50 \\
    20-20-8-4 & 17.11  & 8.13 & 110.57 & 12.70 & 2.56 & 395.82 & 9.75  & 4.00     & 143.75 \\
    20-20-10-4 & 17.45  & 5.50   & 217.33 & 21.53 & 2.73 & 687.53 & 34.64 & 5.17 & 569.78 \\
    20-20-2-6 & 2.73  & 0.91 & 201.77 & 5.70  & 0.17 & 3215.70 & 1.75  & 0.23 & 647.86 \\
    20-20-4-6 & 11.36  & 5.86 & 93.87 & 17.39 & 1.28 & 1257.61 & 8.52  & 3.44 & 147.75 \\
    20-20-6-6 & 24.50  & 10.67 & 129.57 & 28.88 & 4.08 & 608.07 & 27.06 & 8.34 & 224.33 \\
    20-20-8-6 & 43.13  & 14.39 & 199.67 & 80.33 & 10.63 & 656.03 & 57.06 & 16.95 & 236.60 \\
    20-20-10-6 & 50.69  & 24.03 & 110.92 & 137.94 & 34.19 & 303.47 & 69.58 & 32.34 & 115.12 \\
    20-20-2-8 & 2.64  & 0.89  & 196.63 & 2.06  & 0.57 & 277.66 & 1.58  & 0.78 & 101.79 \\
    20-20-4-8 & 12.70  & 5.44 & 133.60 & 18.14 & 2.16 & 741.42 & 8.30  & 4.30 & 93.09 \\
    20-20-6-8 & 43.05  & 16.56 & 159.91 & 87.31 & 10.30 & 747.94 & 33.94 & 16.75 & 102.61 \\
    20-20-8-8 & 68.17  & 34.70 & 96.44 & 303.06 & 65.92 & 359.73 & 69.72 & 56.19 & 24.08 \\
    20-20-10-8 & 61.42  & 30.91 & 98.74 & 1370.64 & 1943.80 & -29.49 & 81.05 & 53.59 & 51.22 \\
    20-20-2-10 & 2.24  & 0.86 & 160.19 & 1.69  & 0.44 & 286.04 & 1.66  & 0.67 & 146.58 \\
    20-20-4-10 & 15.17  & 4.78 & 217.27 & 19.59 & 2.45 & 698.78 & 11.14 & 4.47 & 149.30 \\
    20-20-6-10 & 45.11  & 20.83 & 116.58 & 134.98 & 32.11 & 320.39 & 40.56 & 36.89 & 9.95 \\
    20-20-8-10 & 79.94  & 43.22 & 84.96 & 1619.81 & 1997.02 & -18.89 & 92.58 & 71.72 & 29.08 \\
    20-20-10-10 & 74.02  & 36.63 & 102.09 & 8493.41  & \multicolumn{1}{r}{N/A} & \multicolumn{1}{r}{N/A} & 87.84 & 65.69 & 33.73 \\
    \hline
    \textbf{Average} & 24.54  & 11.30  & \textbf{156.10} & 495.07  & 171.47  & \textbf{566.83} & 26.17  & 15.61  & \textbf{153.95} \\
    \bottomrule
    \end{tabular}%
    }
  \label{Table:AS2}%
\end{table}%

\begin{table}[ht!]
\fontsize{9}{11}\selectfont
  \centering
  \caption{Effectiveness of approximation separation on instances with $|I|=|J|=30$ and varying $p$-, $r$-values}
  \resizebox{0.9\textwidth}{!}{
  \begin{tabular}{lrrrrrrrrr}
    \toprule
    \multicolumn{1}{c}{\multirow{2}[4]{*}{Instance}} & \multicolumn{1}{c}{SC} & \multicolumn{2}{c}{SC-AS} & \multicolumn{1}{c}{BI} & \multicolumn{2}{c}{BI-AS} & \multicolumn{1}{c}{SCBI} & \multicolumn{2}{c}{SCBI-AS} \\
\cmidrule{2-10}          & \multicolumn{1}{c}{\texttt{Time(s)}} & \multicolumn{1}{c}{\texttt{Time(s)}} & \multicolumn{1}{c}{\texttt{Imp(\%)}} & \multicolumn{1}{c}{\texttt{Time(s)}} & \multicolumn{1}{c}{\texttt{Time(s)}} & \multicolumn{1}{c}{\texttt{Imp(\%)}} & \multicolumn{1}{c}{\texttt{Time(s)}} & \multicolumn{1}{c}{\texttt{Time(s)}} & \multicolumn{1}{c}{\texttt{Imp(\%)}} \\
    \midrule
    \multicolumn{1}{l}{30-30-3-3} & \multicolumn{1}{r}{34.11} & 25.88 & 31.82  & 8.33  & 2.77  & 201.19 & 6.06  & 4.88  & 24.37 \\
    \multicolumn{1}{l}{30-30-6-3} & \multicolumn{1}{r}{450.97} & 414.45 & 8.81  & 21.84 & 10.50 & 108.04 & 19.16 & 19.00 & 0.81 \\
    \multicolumn{1}{l}{30-30-9-3} & \multicolumn{1}{r}{LIMIT} & \multicolumn{1}{r}{LIMIT} & \multicolumn{1}{r}{N/A} & 36.02 & 17.55 & 105.25 & 34.19 & 25.19 & 26.33 \\
    \multicolumn{1}{l}{30-30-12-3} & \multicolumn{1}{r}{2198.05} & 2145.91 & 2.43  & 27.42 & 16.16 & 69.73 & 20.20 & 15.33 & 24.13 \\
    \multicolumn{1}{l}{30-30-15-3} & \multicolumn{1}{r}{88.42} & 62.84 & 40.70  & 28.59 & 10.80 & 164.83 & 20.49 & 17.22 & 15.94 \\
    \multicolumn{1}{l}{30-30-3-6} & \multicolumn{1}{r}{143.83} & 94.08 & 52.88  & 33.52 & 13.42 & 149.72 & 34.53 & 18.81 & 45.52 \\
    \multicolumn{1}{l}{30-30-6-6} & \multicolumn{1}{r}{3259.80} & 3196.03 & 2.00  & 140.97 & 54.94 & 156.60 & 130.28 & 107.39 & 17.57 \\
    \multicolumn{1}{l}{30-30-9-6} & \multicolumn{1}{r}{LIMIT} & \multicolumn{1}{r}{LIMIT} & \multicolumn{1}{r}{N/A} & 269.45 & 132.45 & 103.43 & 400.83 & 374.53 & 6.56 \\
    \multicolumn{1}{l}{30-30-12-6} & \multicolumn{1}{r}{LIMIT} & \multicolumn{1}{r}{LIMIT} & \multicolumn{1}{r}{N/A} & 1837.03 & 648.66 & 183.21 & 1981.98 & 2249.34 & -13.49 \\
    \multicolumn{1}{l}{30-30-3-9} & \multicolumn{1}{r}{92.06} & 89.02 & 3.42  & 38.48 & 13.64 & 182.14 & 31.13 & 19.89 & 36.09 \\
    \multicolumn{1}{l}{30-30-6-9} & \multicolumn{1}{r}{2098.56} & 2294.66 & -8.55  & 314.88 & 85.20 & 269.55 & 304.11 & 237.67 & 21.85 \\
    \multicolumn{1}{l}{30-30-3-12} & \multicolumn{1}{r}{55.72} & 42.30 & 31.73  & 25.36 & 4.74  & 435.56 & 20.95 & 16.86 & 19.53 \\
    \multicolumn{1}{l}{30-30-6-12} & \multicolumn{1}{r}{2199.31} & 1954.39 & 12.53  & 1026.81 & 353.95 & 190.10 & 965.25 & 1093.69 & -13.31 \\
    \multicolumn{1}{l}{30-30-3-15} & \multicolumn{1}{r}{39.92} & 36.64 & 8.96  & 22.73 & 5.91  & 284.87 & 21.39 & 22.58 & -5.55 \\
    \multicolumn{1}{l}{30-30-6-15} & \multicolumn{1}{r}{2609.08} & 2430.55 & 7.35  & \multicolumn{1}{r}{LIMIT} & \multicolumn{1}{r}{LIMIT} & \multicolumn{1}{r}{N/A} & 2234.25 & 2410.95 & -7.91 \\
    \hline
    \textbf{Average} & 969.16  & 941.47  & \textbf{16.98} & 153.54  & 52.00  & \textbf{201.12} & 143.05  & 143.03  & \textbf{17.00} \\
    \bottomrule
    \end{tabular}%
    }
  \label{Table:AS3}%
\end{table}%

\subsection{Performance of the Approximation Algorithm}
\label{sec:compu:socp}
We evaluate the effectiveness of the MISOCP approximation algorithm proposed in Section \ref{sec:heuristic}. To do this, we solve the same instances as in Section \ref{sec:perform} by the SCBI implementation and by the MISOCP formulation \eqref{heu}, respectively. We evaluate solution quality of the approximate solution by $\texttt{Obj-GAP} := 1 - z^H/z^*$ and the relative saving on CPU seconds by $\texttt{T-Gap}$ in Tables \ref{table socp-1}--\ref{table socp-2}. From these tables, we observe that the proposed approximation algorithm is able to find good-quality solutions in most instances. For example, the approximate algorithm produces a \texttt{Obj-GAP} below 5\% in 78\% of the instances (43/55). In addition, the algorithm achieves so within an average solution time roughly 45\% of that of SCBI across all instances. Notably, the saving on solution time is particularly significant on challenging instances (e.g., 30-30-12-6 and 30-30-6-15). This demonstrates the applicability of the approximation algorithm in larger-sized instances.

\begin{table}[ht!]
\fontsize{9}{11}\selectfont
  \centering
  \caption{Performance of the approximation algorithm for solving diverse (S-CFLP) instances}
    \label{table socp-1}
    \resizebox{1.0\textwidth}{!}{
    \begin{tabular}{lrrrrrr|lrrrrrr}
    \hline
    \multicolumn{1}{c}{\multirow{2}[2]{*}{Instance}} & \multicolumn{2}{c}{SCBI} & \multicolumn{4}{c|}{AA}       & \multicolumn{1}{c}{\multirow{2}[2]{*}{Instance}} & \multicolumn{2}{c}{SCBI} & \multicolumn{4}{c}{AA} \bigstrut[t]\\
          & \multicolumn{1}{c}{\texttt{Time(s)}} & \multicolumn{1}{c}{\texttt{Obj}} & \multicolumn{1}{c}{\texttt{Time(s)}} & \multicolumn{1}{c}{\texttt{Obj}} & \multicolumn{1}{c}{\texttt{T-GAP}} & \multicolumn{1}{c|}{\texttt{Obj-GAP}} &       & \multicolumn{1}{c}{\texttt{Time(s)}} & \multicolumn{1}{c}{\texttt{Obj}} & \multicolumn{1}{c}{\texttt{Time(s)}} & \multicolumn{1}{c}{\texttt{Obj}} & \multicolumn{1}{c}{\texttt{T-GAP}} & \multicolumn{1}{c}{\texttt{Obj-GAP}} \\
    \hline
    20-20-2-2 & 0.75  & 0.5195 & 0.67  & 0.4768 & 10.40\% & 8.23\% & 30-30-3-3 & 6.06  & 0.5077 & 38.19 & 0.5068 & -529.85\% & 0.17\% \bigstrut[t]\\
    20-20-3-2 & 0.97  & 0.6256 & 1.08  & 0.5917 & -11.25\% & 5.42\% & 30-30-6-3 & 19.16 & 0.6812 & 68.73 & 0.6634 & -258.81\% & 2.61\% \\
    20-20-2-3 & 1.30  & 0.4136 & 0.86  & 0.3797 & 33.69\% & 8.20\% & 30-30-9-3 & 34.19 & 0.7637 & 58.20 & 0.7597 & -70.24\% & 0.52\% \\
    40-40-2-2 & 3.75  & 0.5003 & 4.63  & 0.5000 & -23.33\% & 0.07\% & 30-30-12-3 & 20.20 & 0.8222 & 84.94 & 0.8167 & -320.42\% & 0.66\% \\
    40-40-3-2 & 11.16 & 0.6084 & 6.92  & 0.5878 & 37.95\% & 3.39\% & 30-30-15-3 & 20.49 & 0.8597 & 83.16 & 0.8557 & -305.94\% & 0.47\% \\
    40-40-2-3 & 11.02 & 0.3996 & 4.69  & 0.3993 & 57.44\% & 0.06\% & 30-30-3-6 & 34.53 & 0.3456 & 39.59 & 0.3450 & -14.66\% & 0.17\% \\
    60-60-2-2 & 9.67  & 0.5054 & 19.31 & 0.4889 & -99.67\% & 3.26\% & 30-30-6-6 & 130.28 & 0.5180 & 123.80 & 0.5180 & 4.98\% & 0.00\% \\
    60-60-3-2 & 39.33 & 0.6029 & 69.88 & 0.5978 & -77.67\% & 0.85\% & 30-30-9-6 & 400.83 & 0.6252 & 126.88 & 0.6197 & 68.35\% & 0.88\% \\
    60-60-2-3 & 94.92 & 0.4031 & 21.86 & 0.3949 & 76.97\% & 2.05\% & 30-30-12-6 & 1981.98 & 0.7015 & 383.98 & 0.6919 & 80.63\% & 1.37\% \\
    80-80-2-2 & 25.75 & 0.5003 & 25.89 & 0.4941 & -0.55\% & 1.23\% & 30-30-3-9 & 31.13 & 0.2646 & 45.94 & 0.2643 & -47.59\% & 0.12\% \\
    80-80-3-2 & 146.78 & 0.6060 & 165.27 & 0.5964 & -12.59\% & 1.60\% & 30-30-6-9 & 304.11 & 0.4265 & 214.67 & 0.4250 & 29.41\% & 0.35\% \\
    80-80-2-3 & 228.08 & 0.3954 & 28.53 & 0.3901 & 87.49\% & 1.33\% & 30-30-3-12 & 20.95 & 0.2195 & 53.22 & 0.2171 & -153.99\% & 1.09\% \\
    100-100-2-2 & 44.95 & 0.5014 & 115.77 & 0.5014 & -157.52\% & 0.00\% & 30-30-6-12 & 965.25 & 0.3672 & 502.05 & 0.3648 & 47.99\% & 0.65\% \\
    100-100-3-2 & 190.53 & 0.6040 & 207.72 & 0.6015 & -9.02\% & 0.41\% & 30-30-3-15 & 21.39 & 0.1894 & 60.02 & 0.1745 & -180.58\% & 7.87\% \\
    100-100-2-3 & 273.86 & 0.3970 & 130.64 & 0.3965 & 52.30\% & 0.12\% & 30-30-6-15 & 2234.25 & 0.3263 & 583.34 & 0.3018 & 73.89\% & 7.50\% \\
    \hline
    \textbf{Average} & 72.19 & 0.5055 & 53.58 & 0.4931 & 25.78\% & 2.42\% & \textbf{Average} & 414.99 & 0.5079 & 164.45 & 0.5016 & 60.37\% & 1.23\% \\
    \hline
    \end{tabular}%
  }
\end{table}%

% Table generated by Excel2LaTeX from sheet 'results - 20 20'
\begin{table}[ht!]
\fontsize{9}{11}\selectfont
  \centering
  \caption{Performance of the approximation algorithm on instances with $|I|=|J|=20$ and varying $p$-, $r$-values}
    \label{table socp-2}
    \resizebox{1.0\textwidth}{!}{
    \begin{tabular}{lrrrrrr|lrrrrrr}
    \hline
    \multicolumn{1}{c}{\multirow{2}[2]{*}{Instance}} & \multicolumn{2}{c}{SCBI} & \multicolumn{4}{c|}{AA}       & \multicolumn{1}{c}{\multirow{2}[2]{*}{Instance}} & \multicolumn{2}{c}{SCBI} & \multicolumn{4}{c}{AA} \bigstrut[t]\\
          & \multicolumn{1}{c}{\texttt{Time(s)}} & \multicolumn{1}{c}{\texttt{Obj}} & \multicolumn{1}{c}{\texttt{Time(s)}} & \multicolumn{1}{c}{\texttt{Obj}} & \multicolumn{1}{c}{\texttt{T-GAP}} & \multicolumn{1}{c|}{\texttt{Obj-GAP}} &       & \multicolumn{1}{c}{\texttt{Time(s)}} & \multicolumn{1}{c}{\texttt{Obj}} & \multicolumn{1}{c}{\texttt{Time(s)}} & \multicolumn{1}{c}{\texttt{Obj}} & \multicolumn{1}{c}{\texttt{T-GAP}} & \multicolumn{1}{c}{\texttt{Obj-GAP}} \\
    \hline
    20-20-2-2 & 0.56  & 0.5195  & 6.72  & 0.4768 & -1095.55\% & 8.23\% & 20-20-8-6 & 57.06 & 0.6198  & 14.47 & 0.6001 & 74.65\% & 3.17\% \bigstrut[t]\\
    20-20-4-2 & 1.16  & 0.6963  & 5.52  & 0.6584 & -377.08\% & 5.45\% & 20-20-10-6 & 69.58 & 0.6855  & 16.86 & 0.6773 & 75.77\% & 1.20\% \\
    20-20-6-2 & 1.75  & 0.7693  & 6.49  & 0.7585 & -270.57\% & 1.40\% & 20-20-2-8 & 1.58  & 0.2196  & 8.98  & 0.2159 & -469.33\% & 1.68\% \\
    20-20-8-2 & 1.39  & 0.8254  & 6.70  & 0.8160 & -382.23\% & 1.14\% & 20-20-4-8 & 8.30  & 0.3722  & 6.27  & 0.3476 & 24.48\% & 6.61\% \\
    20-20-10-2 & 1.44  & 0.8592  & 7.25  & 0.8514 & -404.17\% & 0.91\% & 20-20-6-8 & 33.94 & 0.4766  & 12.22 & 0.4606 & 64.00\% & 3.36\% \\
    20-20-2-4 & 1.34  & 0.3400  & 8.86  & 0.3192 & -559.15\% & 6.10\% & 20-20-8-8 & 69.72 & 0.5646  & 18.19 & 0.5449 & 73.91\% & 3.48\% \\
    20-20-4-4 & 3.88  & 0.5248  & 5.91  & 0.5006 & -52.41\% & 4.60\% & 20-20-10-8 & 81.05 & 0.6414  & 35.44 & 0.6303 & 56.28\% & 1.73\% \\
    20-20-6-4 & 6.05  & 0.6344  & 7.61  & 0.6149 & -25.85\% & 3.07\% & 20-20-2-10 & 1.66  & 0.1896  & 7.42  & 0.1865 & -347.92\% & 1.63\% \\
    20-20-8-4 & 9.75  & 0.7024  & 8.13  & 0.6860 & 16.67\% & 2.33\% & 20-20-4-10 & 11.14 & 0.3290  & 8.36  & 0.3076 & 24.97\% & 6.49\% \\
    20-20-10-4 & 34.64 & 0.7541  & 9.64  & 0.7470 & 72.17\% & 0.94\% & 20-20-6-10 & 40.56 & 0.4363  & 16.59 & 0.4132 & 59.09\% & 5.30\% \\
    20-20-2-6 & 1.75  & 0.2675  & 8.58  & 0.2361 & -390.17\% & 11.75\% & 20-20-8-10 & 92.58 & 0.5263  & 29.55 & 0.5062 & 68.08\% & 3.83\% \\
    20-20-4-6 & 8.52  & 0.4355  & 6.92  & 0.4079 & 18.71\% & 6.33\% & 20-20-10-10 & 87.84 & 0.6059  & 49.08 & 0.5794 & 44.13\% & 4.38\% \\
    \cline{8-14}
    20-20-6-6 & 27.06 & 0.5414  & 10.61 & 0.5177 & 60.80\% & 4.37\% & \textbf{Average} & 26.17 & 0.5415  & 12.89 & 0.5224 & 50.73\% & 3.52\% \\
    \hline
    \end{tabular}%
  }
\end{table}%

\subsection{Varying Size of Customers $|I|$} \label{apx-sec:vary-customer-size-I}
In (S-CFLP), a set $I$ of discrete nodes are selected to represent a group of customers, and each node can represent a residential community or a business center located in the considered area. Intuitively, the larger $I$ is, the more representative the model would become, which may also be more challenging to solve. We examine the trade-off between model fidelity and computational burden by varying the size of customers $|I|$. In specific, we evaluate the solution time (\texttt{Time(s)}) and optimal objective value (\texttt{Obj}) of instances $|I|$-$100$-$2$-$2$ with $|I|$ ranging from 20 to 2000 and report the results in Table \ref{table 6}. From this table, we observe that the solution time remains under 5 minutes even if we increase $|I|$ to 2000. Relatively, the solution time increases by roughly 15 times while $|I|$ increases by 100 times. This demonstrates the scalability of Algorithm \ref{algo-b&c}. In addition, we observe that the optimal objective value remains stable when $|I| \geq 100$. This demonstrates that setting $|I| = 100$, as in previous sections, produces sufficiently representative instances.

\begin{table}[ht!]
\fontsize{9}{12}\selectfont
  \centering
  \caption{Computational results of instances $|I|$-$100$-$2$-$2$ with varying $|I|$} 
  \label{table 6}
     \resizebox{\textwidth}{!}{
     \begin{tabular}{lrrrrrrrrrrr}
    \toprule
    $|I|$   & 20 & 40 & 60 & 80 & 100   & 200   & 400   & 800  & 1200  & 1600  & 2000 \\
    \midrule
    \texttt{Time(s)} & 19.00 & 49.52 & 29.23 & 54.80 & 44.95 & 65.69 & 85.70  & 130.50 & 284.80 & 254.91 & 279.97 \\
    \texttt{Obj} & 0.5021 & 0.5004 & 0.5013 & 0.5003 & 0.5014 & 0.5011 & 0.5007 & 0.5009 & 0.5001 & 0.5002 & 0.5072 \\
    \texttt{\#Cuts}  & 161   & 224   & 135   & 225   & 176   & 197   & 184   & 176   & 228   & 187   & 164 \\
    \bottomrule
    \end{tabular}%
    }
\end{table}%

\end{appendices}

%%%%%%%%%%%%%%%%%
\end{document}